\documentclass[11pt]{amsart}
\usepackage{amsmath}
\usepackage{a4wide}
\usepackage[utf8]{inputenc}
\usepackage{amssymb}
\usepackage{amsopn}
\usepackage{epsfig}
\usepackage{amsfonts}
\usepackage{latexsym}
\usepackage{amsthm}
\usepackage{enumerate}
\usepackage[UKenglish]{babel}
\usepackage{verbatim}
\usepackage{color}
\usepackage{calrsfs}
\usepackage{pgf}
\usepackage{tikz}
\usepackage{subfig}
\usetikzlibrary{math}
\usetikzlibrary{positioning} 
\usepackage{xcolor} 

\usepackage{mathrsfs}   
\usepackage{bm}

 \usetikzlibrary{arrows,automata}

\DeclareMathAlphabet{\pazocal}{OMS}{zplm}{m}{n}
\newtheorem{theorem}{Theorem}[section]
\newtheorem{lemma}[theorem]{Lemma}
\newtheorem{proposition}[theorem]{Proposition}

\theoremstyle{definition}
\newtheorem{definition}[theorem]{Definition}
\newtheorem{example}[theorem]{Example}

\theoremstyle{remark}
\newtheorem{remark}[theorem]{Remark}

\numberwithin{equation}{section}

\newcommand{\R}{\ensuremath{\mathbb{R}}}
\newcommand{\N}{\ensuremath{\mathbb{N}}}

\newcommand{\+}{ {\ensuremath{\oplus}}}
\renewcommand{\b}{ {\mathbf{b}}}
\renewcommand{\c}{ {\mathbf{c}}}
\renewcommand{\d}{ {\mathbf{d}}}

\renewcommand{\t}{ {\mathbf{t}}}

\newcommand{\re}{\ensuremath{\phi}}

\newcommand{\U} {{\mathbf{U}}}

\newcommand{\su}{ {\mathbf{u}}}
\newcommand{\sv}{ {\mathbf{v}}}

\newcommand{\Om}{ {\Omega}}

\newcommand{\set}[1]{\left\{#1\right\}}

\newcommand{\la}{\lambda}
\newcommand{\ga}{\gamma}

\newcommand{\f}{\infty}
\newcommand{\de}{\delta}
\newcommand{\om}{\omega}
\newcommand{\al}{\alpha}
\newcommand{\lle}{\preccurlyeq}
\newcommand{\lge}{\succcurlyeq}
\renewcommand{\a}{ \mathbf{a}}
\newcommand{\si}{\sigma}

\newcommand{\ra}{\rightarrow}

\begin{document}

\title[Phase transitions for unique codings of fat Sierpinski Gaskets]{Phase transitions for unique codings of fat Sierpinski Gaskets}

\author{Yi Cai}
\address[Y. Cai]{College of Science, Shanghai Institute of Technology
Shanghai 201418, People's Republic of China}
\email{caiyi@sit.edu.cn}

\author{Derong Kong}
\address[D. Kong]{College of Mathematics and Statistics, Center of Mathematics, Chongqing University, Chongqing, 401331, P.R.China.}
\email{derongkong@126.com}

\author{Wenxia Li}
\address[W. Li]{School of Mathematical Sciences, Key Laboratory of MEA (Ministry of Education) \& Shanghai Key Laboratory of PMMP, East China Normal University, Shanghai 200241, People's Republic of China}
\email{wxli@math.ecnu.edu.cn}

\author{Yuhan Zhang}
\address[Y. Zhang]{College of Mathematics and Statistics,   Chongqing University, Chongqing, 401331, P.R.China}
\email{zyh240207@163.com}

\dedicatory{}



\begin{abstract}
Given an integer  $M\ge 1$ and $\beta\in(1, M+1)$, let $S_{\beta, M}$ be the fat Sierpinski gasket in $\mathbb R^2$ generated by the iterated function system $\left\{f_d(x)=\frac{x+d}{\beta}: d\in\Omega_M\right\}$, where $\Omega_M=\{(i,j)\in\mathbb Z_{\ge 0}^2: i+j\le M\}$. Then each $x\in S_{\beta, M}$ may be represented as a series  $x=\sum_{i=1}^\infty\frac{d_i}{\beta^i}=:\Pi_\beta((d_i))$, and the infinite sequence $(d_i)\in\Omega_M^{\mathbb N}$ is called a \emph{coding} of $x$. Since $\beta<M+1$, a point in $S_{\beta, M}$ may have multiple codings. Let $U_{\beta, M}$ be the set of $x\in S_{\beta, M}$ having a unique coding, that is
\[
U_{\beta, M}=\set{x\in S_{\beta, M}: \#\Pi_\beta^{-1}(x)=1}.
\]
When $M=1$, Kong and Li [2020, Nonlinearity] described  two critical bases for the phase transitions of the intrinsic univoque set $\widetilde U_{\beta, 1}$, which is a subset of $U_{\beta, 1}$. In this paper we consider $M\ge 2$, and characterize the two critical bases $\beta_G(M)$ and $\beta_c(M)$ for the phase transitions of $U_{\beta, M}$: (i) if $\beta\in(1, \beta_G(M)]$, then   $U_{\beta, M}$ is finite; (ii) if $\beta\in(\beta_G(M), \beta_c(M))$ then $U_{\beta, M}$ is countably infinite; (iii) if $\beta=\beta_c(M)$ then $U_{\beta, M}$ is uncountable and has zero Hausdorff dimension; (iv) if $\beta>\beta_c(M)$ then $U_{\beta, M}$ has positive Hausdorff dimension.  Moreover, we show that the first critical base $\beta_G(M)$ is a Perron number, while the second critical base $\beta_c(M)$ is a transcendental number.  
  \end{abstract}

\maketitle

\section{Introduction}
Given $\beta>1$ and $M\in\N$, let $S_{\beta, M}$ be the {\em Sierpinski gasket} in $\R^2$ generated by the iterated function system (IFS) $\set{f_{\beta, d}(x)=\frac{x+d}{\beta}: d\in\Omega_M}$, where 
\[
\Omega_M:=\set{\al_{ij}=(i,j): i, j\in\set{0,1,\ldots, M}\textrm{ and } i+j\le M}.
\]
Note that each digit $\alpha_{ij}=(i,j)$ is a vector in $\R^2$ which determines the transition of $f_{\beta,\al_{ij}}$.
It is well-known that $S_{\beta, M}$ is the unique nonempty compact set  satisfying the equation $S_{\beta, M}=\bigcup_{d\in{\color{black}\Omega_M}}f_{\beta, d}(S_{\beta, M})$ (cf.~\cite{Hutchinson_1981}). To study   $S_{\beta, M}$ it is useful to introduce the \emph{coding} map
\begin{equation}\label{eq:coding-map}
	\Pi_\beta: \Omega_M^\N\to S_{\beta, M};\quad (d_i)\mapsto\sum_{i=1}^{\f}\frac{d_i}{\beta^i},
\end{equation}
which is always continuous and  surjective. If $\beta>M+1$, the IFS $\mathcal F_{\beta, M}=\set{f_{\beta, d}: d\in\Omega_M}$ satisfies the strong separation condition, and then $\Pi_\beta$ is bijective (cf.~\cite{Falconer_1990}). If $\beta=M+1$, the IFS $\mathcal F_{\beta,M}$ fails the strong separation condition but still admits  the open set condition, and thus  $\Pi_\beta$ is nearly bijective up to a countable set. However, for $\beta<M+1$ the IFS $\mathcal F_{\beta, M}$ fails the open set condition and it always has overlaps. In particular, if $\beta\in(1, \frac{M+2}{2}]$ the attractor $S_{\beta, M}$ is an isosceles triangle with vertices $(0,0), (\frac{M}{\beta-1}, 0)$ and $(0, \frac{M}{\beta-1})$ (see Proposition \ref{prop:U-beta-coincidence}). In this case, for Lebesgue almost every $x\in S_{\beta, M}$ the preimage $\Pi_\beta^{-1}(x)$ is uncountable (cf.~\cite{Sidorov_2007}). 
When $\beta\in(\frac{M+2}{2}, M+1)$ the structure of $S_{\beta, M}$ is getting more complicated. Broomhead et al.~\cite{Bro-Mon-Sid-04} showed that for $M=1$ and $\beta\le \beta_*$ the set $S_{\beta, M}$ has non-empty interior, where $\beta_*\approx 1.54369$ is a root of $x^3-2x^2+2x-2=0$. 

When $\beta\in(1, M+1)$ the \emph{fat} Sierpinski gasket has attracted much attention in the past twenty years. Simon and Solomyak \cite{Simon-Solomyak-02} showed that for $M=1$ there exists a dense set of $\beta\in(1, 2)$ such that $\dim_H S_{\beta, 1}<\log 3/\log\beta$, where $\log 3/\log\beta$ is the self-similar dimension of $S_{\beta, 1}$. 
On the other hand, Jordan and Pollicott \cite{Jordan-Pollicott-06} proved that for $M=1$  the fat Sierpinski gasket $S_{\beta, 1}$ has positive Lebesgue measure for Lebesgue almost every $\beta\in(1, 1.70853]$.  
Recently, Hochman \cite{Hochman-2015} showed that $\dim_H S_{\beta, M}=\min\set{2, \dim_S S_{\beta, M}}$ for $\beta\in(1,M+1)$ up to a set of zero packing dimension, where $\dim_S S_{\beta, M}=\frac{\log(M+1)+\log(M+2)-\log 2}{\log \beta}$ is the similarity dimension of $S_{\beta, M}$.

\begin{figure}[h!]
	\begin{center}
		\begin{tikzpicture}[
			scale=6,  
			axis/.style={very thick, ->},
			important line/.style={thick},
			dashed line/.style={dashed, thin},
			pile/.style={thick, ->, >=stealth', shorten <=2pt, shorten
				>=2pt},
			every node/.style={color=black}
			]
			
			\fill[black!20]({0},0)--({4/(4-1)-4/4},0)--(0,{4/(4-1)-4/4})--cycle; 
			\fill[black!20]({1/4},0)--({4/(4-1)-3/4},0)--({1/4},{4/(4-1)-4/4})--cycle; 
			\fill[black!20]({2/4},0)--({4/(4-1)-2/4},0)--({2/4},{4/(4-1)-4/4})--cycle; 
			\fill[black!20]({3/4},0)--({4/(4-1)-1/4},0)--({3/4},{4/(4-1)-4/4})--cycle; 
			\fill[black!20]({4/4},0)--({4/(4-1)},0)--({4/4},{4/(4-1)-4/4})--cycle; 
			
			\fill[black!20]({0},{1/4})--({4/(4-1)-4/4},{1/4})--(0,{4/(4-1)-3/4})--cycle; 
			\fill[black!20]({1/4},{1/4})--({4/(4-1)-3/4},{1/4})--({1/4},{4/(4-1)-3/4})--cycle; 
			\fill[black!20]({2/4},{1/4})--({4/(4-1)-2/4},{1/4})--({2/4},{4/(4-1)-3/4})--cycle; 
			\fill[black!20]({3/4},{1/4})--({4/(4-1)-1/4},{1/4})--({3/4},{4/(4-1)-3/4})--cycle; 
			
			\fill[black!20]({0},{2/4})--({4/(4-1)-4/4},{2/4})--(0,{4/(4-1)-2/4})--cycle; 
			\fill[black!20]({1/4},{2/4})--({4/(4-1)-3/4},{2/4})--({1/4},{4/(4-1)-2/4})--cycle; 
			\fill[black!20]({2/4},{2/4})--({4/(4-1)-2/4},{2/4})--({2/4},{4/(4-1)-2/4})--cycle;

			\fill[black!20]({0},{3/4})--({4/(4-1)-4/4},{3/4})--(0,{4/(4-1)-1/4})--cycle; 
			\fill[black!20]({1/4},{3/4})--({4/(4-1)-3/4},{3/4})--({1/4},{4/(4-1)-1/4})--cycle; 
			
			\fill[black!20]({0},{4/4})--({4/(4-1)-4/4},{4/4})--(0,{4/(4-1)})--cycle; 
			
			\fill[black!40]({1/4},0)--({4/(4-1)-4/4},0)--({1/4},{4/(4-1)-5/4})--cycle;    
			\fill[black!40]({2/4},0)--({4/(4-1)-3/4},0)--({2/4},{4/(4-1)-5/4})--cycle; 
			\fill[black!40]({3/4},0)--({4/(4-1)-2/4},0)--({3/4},{4/(4-1)-5/4})--cycle; 
			\fill[black!40]({4/4},0)--({4/(4-1)-1/4},0)--({4/4},{4/(4-1)-5/4})--cycle; 
			
			\fill[black!40](0,{1/4})--(0,{4/(4-1)-4/4})--({4/(4-1)-5/4},{1/4})--cycle;    
			\fill[black!40](0,{2/4})--(0,{4/(4-1)-3/4})--({4/(4-1)-5/4},{2/4})--cycle; 
			\fill[black!40](0,{3/4})--(0,{4/(4-1)-2/4})--({4/(4-1)-5/4},{3/4})--cycle; 
			\fill[black!40](0,{4/4})--(0,{4/(4-1)-1/4})--({4/(4-1)-5/4},{4/4})--cycle; 
			
			\fill[black!40]({1/4},{1/4})--({4/(4-1)-4/4},{1/4})--({1/4},{4/(4-1)-4/4})--cycle;    
			\fill[black!40]({2/4},{1/4})--({4/(4-1)-3/4},{1/4})--({2/4},{4/(4-1)-4/4})--cycle; 
			\fill[black!40]({3/4},{1/4})--({4/(4-1)-2/4},{1/4})--({3/4},{4/(4-1)-4/4})--cycle; 
			\fill[black!40]({4/4},{1/4})--({4/(4-1)-1/4},{1/4})--({4/4},{4/(4-1)-4/4})--cycle; 
			
			\fill[black!40]({1/4},{2/4})--({4/(4-1)-4/4},{2/4})--({1/4},{4/(4-1)-3/4})--cycle;    
			\fill[black!40]({2/4},{2/4})--({4/(4-1)-3/4},{2/4})--({2/4},{4/(4-1)-3/4})--cycle; 
			\fill[black!40]({3/4},{2/4})--({4/(4-1)-2/4},{2/4})--({3/4},{4/(4-1)-3/4})--cycle;

			\fill[black!40]({1/4},{3/4})--({4/(4-1)-4/4},{3/4})--({1/4},{4/(4-1)-2/4})--cycle;    
			\fill[black!40]({2/4},{3/4})--({4/(4-1)-3/4},{3/4})--({2/4},{4/(4-1)-2/4})--cycle; 
			
			\fill[black!40]({1/4},{4/4})--({4/(4-1)-4/4},{4/4})--({1/4},{4/(4-1)-1/4})--cycle;  
			
			\draw[axis] (-0.12,0)  -- ({(4/(4-1))+0.1},0) node(xline)[right]
			{$x$};
			\draw[axis] (0,-0.12) -- (0,{4/(4-1)+0.1}) node(yline)[above] {$y$};
			\node[] at (-0.07,-0.07){$0$};
			\draw[important line] (0,{4/(4-1)})--({4/(4-1)},0);
			
			\draw[important line] ({1/4},0)--({1/4},{4/(4-1)-1/4});
			\draw[important line] ({2/4},0)--({2/4},{4/(4-1)-2/4});
			\draw[important line] ({3/4},0)--({3/4},{4/(4-1)-3/4});
			\draw[important line] ({4/4},0)--({4/4},{4/(4-1)-4/4});
			
			\draw[important line] (0,{1/4})--({4/(4-1)-1/4},{1/4});
			\draw[important line] (0,{2/4})--({4/(4-1)-2/4},{2/4});
			\draw[important line] (0,{3/4})--({4/(4-1)-3/4},{3/4});
			\draw[important line] (0,{4/4})--({4/(4-1)-4/4},{4/4});
			
			\draw[important line] ({4/(4-1)-2/4},0)--(0,{4/(4-1)-2/4});
			\draw[important line] ({4/(4-1)-1/4},0)--(0,{4/(4-1)-1/4});
			\draw[important line] ({4/(4-1)-3/4},0)--(0,{4/(4-1)-3/4});
			\draw[important line] ({4/(4-1)-4/4},0)--(0,{4/(4-1)-4/4});
			
			\node[] at ({1/(2*4)}, {1/(2*4)}){$\al_{00}$};
			
			\node[] at ({3/(2*4)}, {1/(2*4)}){$\al_{10}$};
			
			\node[] at ({5/(2*4)}, {1/(2*4)}){$\al_{20}$};
			\node[] at ({7/(2*4)}, {1/(2*4)}){$\al_{30}$};
			\node[] at ({9/(2*4)}, {1/(2*4)}){$\al_{40}$};

			\node[] at ({1/(2*4)}, {3/(2*4)}){$\al_{01}$};
			
			\node[] at ({1/(2*4)}, {5/(2*4)}){$\al_{02}$};
			\node[] at ({1/(2*4)}, {7/(2*4)}){$\al_{03}$};
			\node[] at ({1/(2*4)}, {9/(2*4)}){$\al_{04}$};

			\node[] at ({3/(2*4)}, {3/(2*4)}){$\al_{11}$};
			\node[] at ({3/(2*4)}, {5/(2*4)}){$\al_{12}$};
			\node[] at ({3/(2*4)}, {7/(2*4)}){$\al_{13}$};

			\node[] at ({5/(2*4)}, {3/(2*4)}){$\al_{21}$};
			\node[] at ({7/(2*4)}, {3/(2*4)}){$\al_{31}$};

			\node[] at ({5/(2*4)}, {5/(2*4)}){$\al_{22}$};

			\node[] at ({1/(4)}, {-0.08}){$\frac{1}{\beta}$};
			
			\node[] at ({2/(4)}, {-0.08}){$\frac{2}{\beta}$};
			\node[] at ({3/(4)}, {-0.08}){$\frac{3}{\beta}$};
			\node[] at ({4/(4)}, {-0.08}){$\frac{4}{\beta}$};
			\node[] at ({4/(4-1)}, {-0.08}){$\frac{4}{\beta-1}$};
			
			\node[] at ( {-0.08}, {1/(4)}){$\frac{1}{\beta}$};
			
			\node[] at ( {-0.08}, {2/(4)}){$\frac{2}{\beta}$};
			\node[] at ( {-0.08}, {3/(4)}){$\frac{3}{\beta}$};
			\node[] at ( {-0.08}, {4/(4)}){$\frac{4}{\beta}$};
			\node[] at ( {-0.08}, {4/(4-1)}){$\frac{4}{\beta-1}$};
			%
			
			%
			%
			%

		\end{tikzpicture}
	\end{center}
	\caption{The graph for the first generation of $S_{\beta, M}$ with $M=4$ and $\beta=4$.   The convex hull $\Delta_{\beta, M}$ is the triangle with vertices $(0, 0), (4/(\beta-1), 0)$ and $(0, 4/(\beta-1))$. Each light grey triangle corresponds to a $f_{\beta, \al_{ij}}(\Delta_{\beta,M})$ for some $\al_{ij}\in\Om_M$, and the overlap region $O_{\beta, M}$ is the union of   eighteen  small dark grey triangles.}\label{Fig:S-4}
\end{figure}

Observe that for $\beta\in(1, M+1)$ the overlap region
\[
O_{\beta, M}=\bigcup_{c, d\in\Om_{M}, c\ne d}f_{\beta, c}(\Delta_{\beta,M})\cap f_{\beta, d}(\Delta_{\beta, M})
\]
is nonempty, where $\Delta_{\beta, M}$ is the convex hull of $S_{\beta, M}$ (see Figure \ref{Fig:S-4}). Motivated by the work of  Kong and Li \cite{Kong-Li-2020}, we considered the \emph{intrinsic univoque set}
\begin{equation}\label{eq:U-beta}
	\widetilde{U}_{\beta, M}:=\set{\Pi_\beta((d_i))\in S_{\beta, M}: \Pi_\beta(d_{n+1}d_{n+2}\ldots)\notin O_{\beta, M}~\forall n\ge 0},
\end{equation}
where $\Pi_\beta$ is defined in (\ref{eq:coding-map}). Note that $\widetilde{U}_{\beta, M}$ is the survivor set of some open dynamical systems (cf.~\cite{Pianigiani-Yorke-1979}). Indeed, let $T_{\beta, M}$ be the expanding map from the isosceles triangle $\Delta_{\beta, M}$ to itself such that $T_{\beta, M}(x)=f_{\beta, d}^{-1}(x)$ if $x\in f_{\beta, d}(\Delta_{\beta, M})$ for any $d\in\Om_M$. Then $S_{\beta, M}$ is invariant under $T_{\beta, M}$, and the intrinsic univoque set $\widetilde{U}_{\beta, M}$ is the survivor set consisting of all $x\in\Delta_{\beta, M}$ whose orbit under $T_{\beta, M}$ never enters  the hole $O_{\beta, M}$, i.e., $\widetilde{U}_{\beta, M}=\set{x\in\Delta_{\beta, M}: T_{\beta, M}^n(x)\notin O_{\beta, M}~\forall n\ge 0}$. 
For more information on open dynamical systems we refer to the papers \cite{Urbanski_1986, Demmers-2005, Demers-Young-2006, Demers-Wright-Young-2010, Glendinning-Sidorov-2015, Lyndsey-2016, Bunimovich-Yurchenko-2011, Kalle-Kong-Langeveld-Li-18, Allaart-Kong-2023} and the references therein. 

Since $\beta\in(1, M+1)$, $S_{\beta, M}$ is a self-similar set with overlaps. Inspired by the work of unique $\beta$-expansions (cf.~\cite{Erdos_Joo_Komornik_1990, DeVries_Komornik_2008, Komornik-Kong-Li-17}), we  study the \emph{univoque set} $U_{\beta, M}$ in $S_{\beta, M}$, i.e.,  
\begin{equation}\label{eq:U-beta-intrinsic}
	U_{\beta, M}=\set{x\in S_{\beta, M}: \#\Pi_\beta^{-1}(x)=1},
\end{equation}
where $\# A$ denotes the cardinality of a set $A$. Note that each point in $\widetilde{U}_{\beta, M}$ has a unique coding under $\Pi_\beta$. Then $\widetilde{U}_{\beta, M}\subset U_{\beta, M}$ for any $\beta\in(1, M+1)$. 

When $\beta\in(1,M+1)$ is getting larger, the overlap region $O_{\beta, M}$ is getting relatively smaller, and it turns out that  $\widetilde{U}_{\beta, M}$ is getting larger. In this paper we will describe the phase transition phenomenon for $\widetilde{U}_{\beta, M}$ and $U_{\beta, M}$. More precisely, we will describe   two critical bases $\beta_G(M)$ and $\beta_c(M)$ with $\beta_G(M)<\beta_c(M)$. The first critical base $\beta_G(M)$, called the \emph{generalized golden ratio}, is the unique base such that the cardinality of    $\widetilde{U}_{\beta, M}$ increases  from finite to infinite when $\beta$ passes $\beta_G(M)$. The second critical base $\beta_c(M)$, called the \emph{generalized Komornik-Loreti constant}, is the unique base in which the Hausdorff dimension  of   $\widetilde{U}_{\beta, M}$ changes from zero to positive when $\beta$ crosses $\beta_c(M)$. When $M\ge 2$, we will show that the two critical bases $\beta_G(M)$ and $\beta_c(M)$ can  also be applied to $U_{\beta, M}$.

{\color{black}Write $\N:=\set{1,2,3,\ldots}$, $\N_0:=\N\cup\set{0}$.} Let $(\tau_i)_{i=0}^\f=01101001\ldots\in\set{0,1}^\N$ be the classical Thue-Morse sequence defined by (cf.~\cite{Allouche_Shallit_1999})
\[
\tau_{0}=0,\quad \tau_{2i+1}=1-\tau_i\quad \textrm{and}\quad \tau_{2i}=\tau_i\quad\forall i\ge 0.
\] Now we are ready to define the two critical bases $\beta_G(M)$ and $\beta_c(M)$. 
\begin{definition}
	\label{def:beta-G}
	The \emph{generalized golden ratio} $\beta_G(M)$ for $M\ge 1$ is defined as follows. If $M=3N+1$ for some $N\in\N_{0}:=\set{0,1,2,\ldots}$, then $\beta_G(3N+1)$ is the unique root in $(N+1, N+2)$ of the equation
	\[
	x^3-(N+1)x^2-N x-(N+1)=0.
	\]
	Similarly, if $M=3N+2$ with $N\in\N_0$, then $\beta_G(3N+2)$ is the unique root in $(N+1, N+2)$ of the equation
	\[
	x^3-(N+1)(x^2+x+1)=0.
	\]
	Finally, if $M=3N+3$ with $N\in\N_0$, then   $\beta_G(3N+3)=N+2$. 
\end{definition}

\begin{figure}[h!]
	\centering
	\includegraphics[width=13cm]{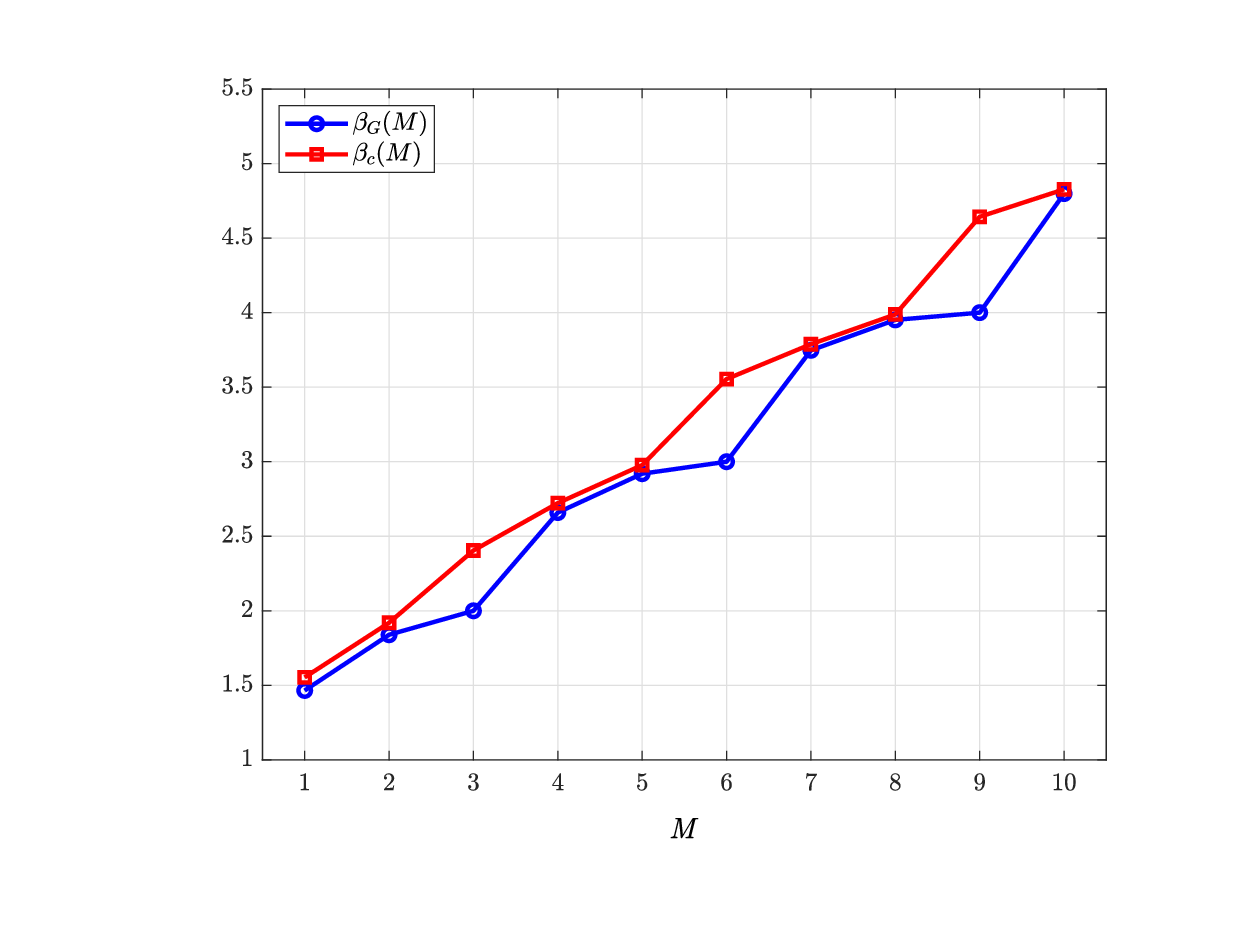}
	\caption{The graph of the two critical bases $\beta_G(M)$ and $\beta_c(M)$ with $M=1,2,\ldots, 10$. }\label{fig:critical-values}
\end{figure}

\begin{definition}
	\label{def:beta-c}The \emph{generalized Komornik-Loreti constant} $\beta_c(M)$ for $M\ge 1$ is defined as follows. If $M=3N+1$ for some $N\in\N_0$, then   $\beta_c(3N+1)$ is the unique root in $(N+1, N+2)$ of the equation
	\[
	\sum_{i=0}^{\f}\left(\frac{N+\tau_{2i+1}}{x^{3i+1}}+\frac{N}{x^{3i+2}}+\frac{N+\tau_{2i+2}}{x^{3i+3}}\right)=1,
	\]
	where $(\tau_i)_{i=0}^\f=01101001\ldots$ is the classical Thue-Morse sequence. Similarly, if $M=3N+2$ with $N\in\N_0$, then $\beta_c(3N+2)$ is   the unique root in $(N+1, N+2)$ of the equation
	\[
	\sum_{i=0}^{\f}\left(\frac{N+\tau_{2i+1}}{x^{3i+1}}+\frac{N+1}{x^{3i+2}}+\frac{N+\tau_{2i+2}}{x^{3i+3}}\right)=1.
	\]
	Finally, if $M=3N+3$ with $N\in\N_0$, then   $\beta_c(3N+3)$ is the unique root in $(N+2, N+3)$ of the equation
	\[
	\sum_{i=0}^{\f}\left(\frac{N+2\tau_{2i+1}}{x^{3i+1}}+\frac{N+\tau_{2i}}{x^{3i+2}}+\frac{N+1+\tau_{2i+2}}{x^{3i+3}}\right)=1.
	\] 
\end{definition}
By Definitions \ref{def:beta-G} and \ref{def:beta-c} it follows that 
\[{\color{black}\left[\frac{M}{3}\right]}+1\le \beta_G(M)<\beta_c(M)<{\color{black}\left[\frac{M}{3}\right]}+2,\]
where $[r]$ denotes the integer part of a real number $r$. Then $\lim_{M\to\f}\frac{\beta_G(M)}{M}=\lim_{M\to\f}\frac{\beta_c(M)}{M}=\frac{1}{3}$.
In the following table we numerically calculate the two critical bases $\beta_G(M)$ and $\beta_c(M)$ for $M=1,2,\ldots, 10$; and they are also  illustrated in Figure \ref{fig:critical-values}.
\begin{table}[htbp] 
	\centering  
	\resizebox{\linewidth}{!}{ 
		\begin{tabular}{c|c|c|c|c|c|c|c|c|c|c}
			\hline
			$M$&1&2&3&4&5&6&7&8&9&10\\\hline
			$\beta_G(M)\approx$&1.46557&1.83929&2&2.65897&2.91964&3&3.74735&3.95137&4&4.79885\\\hline
			$\beta_c(M)\approx$&1.55356&1.91988&2.40462&2.72236&2.97737&3.55447&3.78826&3.98799&4.64302&4.82717\\\hline
	\end{tabular}}
	\caption{ {\color{black}$\beta_G(M)$ and $\beta_c(M)$ with $M=1, 2,\ldots, 10$.}} 
	\label{tab1}  
\end{table}

We will show that each $\beta_G(M)$ is a Perron number, and each $\beta_c(M)$ is a transcendental number. When $M=1$ the following phase transition phenomenon  for $\widetilde{U}_{\beta,1}$ was established by 
Sidorov \cite[Theorem 4.1]{Sidorov_2007} for the cardinality transition,  and by Kong and Li \cite[Theorem 2]{Kong-Li-2020} for the dimension transition.
\begin{theorem}[\cite{Sidorov_2007, Kong-Li-2020}]\label{th:sidorov-kong-li}
	Let $\beta_G(1)\approx 1.46557$ and $\beta_c(1)\approx 1.55356$ be defined as in Definitions \ref{def:beta-G} and \ref{def:beta-c}, respectively.   Then $\beta_G(1)$ is a Perron number, and  $\beta_c(1)$ is a transcendental number. 
	\begin{enumerate}[{\rm(i)}]
		\item If $\beta\in(1, \beta_G(1)]$, then $\#\widetilde{U}_{\beta, 1}=\#U_{\beta, 1}=3$;
		\item If $\beta\in(\beta_G(1), \beta_c(1))$, then   $\widetilde{U}_{\beta, 1}$ is countably infinite;
		\item If $\beta=\beta_c(1)$, then   $\widetilde{U}_{\beta, 1}$   is uncountable and has zero Hausdorff dimension;
		\item If $\beta\in(\beta_c(1), 2)$, then $\widetilde{U}_{\beta, 1}$ has positive Hausdorff dimension.
	\end{enumerate}
\end{theorem}
Since $\widetilde{U}_{\beta, 1}\subset U_{\beta, 1}$, by Theorem \ref{th:sidorov-kong-li} it follows that  $U_{\beta, 1}$ is infinite for $\beta>\beta_G(1)$. So,   Theorem \ref{th:sidorov-kong-li} implies that $\beta_G(1)$ is also a critical base for $U_{\beta, 1}$, i.e., $U_{\beta, 1}$ is finite if and only if $\beta\le \beta_G(1)$.      However, we don't know whether $\beta_c(1)$ is still the critical base for $U_{\beta, 1}$ in the sense that $U_{\beta, 1}$ has positive Hausdorff dimension if and only if $\beta>\beta_c(1)$.

In this paper we extend Theorem \ref{th:sidorov-kong-li} to   $M\ge 2$, and show that (i) $\widetilde{U}_{\beta, M}$ and $U_{\beta, M}$ are finite if and only if $\beta\le\beta_G(M)$;  (ii) $\widetilde{U}_{\beta, M}$ and $U_{\beta, M}$ have positive Hausdorff dimension if and only if $\beta>\beta_c(M)$. 
{The new challenge is that for $M=3N+1$ with $N\ge 1$ the digit $N-1$ might occur  in $d_{n}^1$ or $d_n^2$ for a sequence $(d_i)\in\widetilde{\U}_{\beta_c, M}$, while for $M=1$ this can never happen. Fortunately, we are able to show  that this digit $N-1$ can occur at most 3 times in any sequence of $\widetilde{\U}_{\beta_c, M}$ (see Lemma \ref{lem:exceptional-digits}). Furthermore, completely different from the case for $M=1$, we need more effort to find the critical base $\beta_c(M)$ for $M=3N+3$. In other words, we need to construct a completely new type of sequences to represent $\beta_c(M)${\color{black};} see more in Section \ref{sec:3N}.}

Now we state the main result of the paper.
\begin{theorem}
	\label{th:main-result}
	Given $M\ge 2$,  let $\beta_G(M)$ and $\beta_c(M)$ be defined as in Definitions \ref{def:beta-G} and \ref{def:beta-c} respectively. Then 
	$\beta_G(M)$ is a Perron number, and $\beta_c(M)$ is a transcendental number. 
	\begin{enumerate}[{\rm(i)}]
		\item If $\beta\in(1, \beta_G(M)]$, then $\# U_{\beta, M}=\#\widetilde{U}_{\beta, M}=3$;
		\item If $\beta\in(\beta_G(M),\beta_c(M))$, then $U_{\beta, M}=\widetilde{U}_{\beta, M}$ is countably infinite;
		\item If $\beta=\beta_c(M)$, then $U_{\beta, M}=\widetilde{U}_{\beta, M}$ is uncountable and has zero Hausdorff dimension;
		\item If $\beta>\beta_c(M)$, then both $U_{\beta, M}$ and $\widetilde{U}_{\beta, M}$ have positive Hausdorff {\color{black}dimension}.
	\end{enumerate}
\end{theorem} 
When $\beta>\beta_c(M)$ and $\beta$ is a generalized multinacci number, we are able to give the explicit formula for the Hausdorff dimension of $U_{\beta, M}$ (see Example \ref{ex:1} below).

The rest of the paper is organized as follows. In the next section we give the characterization of $\widetilde{U}_{\beta, M}$. In Section \ref{sec:golden-ratio} we show that the generalized golden ratio $\beta_G(M)$ defined in Definition \ref{def:beta-G} is the first critical base for $\widetilde{U}_{\beta, M}$ {\color{black}in the sense that}  $\widetilde{U}_{\beta, M}$ is finite if and only if $\beta\le \beta_G(M)$. Furthermore, we show that the generalized Komornik-Loreti constant $\beta_c(M)$ in Definition \ref{def:beta-c} is the second critical base for $\widetilde{U}_{\beta, M}$. This will be proved in Sections \ref{sec:3N+1}-- \ref{sec:3N} for $M=3N+1, M=3N+2$ and $M=3N+3$, respectively. We point out that the two critical bases $\beta_G(M)$ and $\beta_c(M)$ are also applied to $U_{\beta, M}$ by Proposition  \ref{prop:U-beta-coincidence}.

\section{Characterization of $\widetilde{U}_{\beta, M}$}\label{sec:characterization-U}
Given $M\in\N$, we  recall  
$\Om_M=\set{\al_{ij}=(i,j)\in\N^2_0:  i+j\leq M}.$
For $d=(d^1, d^2)\in\Omega_M$ let   $d^\+:=d^1+d^2$ and $\overline{d^\+}:=M-d^\+$. Then 
\begin{equation}\label{eq:sum-equal}
	d^1+d^2+\overline{d^\+}=M\quad\textrm{for any }d\in\Omega_M.
\end{equation}
First we introduce  some terminology from symbolic dynamics (cf.~\cite{Lind_Marcus_1995}). For $n\in\N$ let $\Om_M^n$ be the set of all length $n$ words with each digit from $\Om_M$. Denote by $\Om_M^*:=\bigcup_{n=1}^\f\Om_M^n$  the set of all finite words over the alphabet  $\Om_M$. Moreover, let $\Om_M^\N$ be the set of all infinite sequences  with each digit from $\Om_M$. For a word $d_1d_2\ldots d_n\in\Om_M^*$ or   a sequence $(d_i)\in\Om_M^\N$ it is useful to represent it in a matrix form: 
\[
d_1d_2\ldots d_n\sim\left(\begin{array}{c}
	d_1^1d_2^1\ldots d_n^1 \\
	d_1^2d_2^2\ldots d_n^2 \\
	\overline{d_1^\+}\,\overline{d_2^\+}\ldots \overline{d_n^\+} 
\end{array}\right), \quad 
d_1d_2\ldots \sim\left(\begin{array}{c}
	d_1^1d_2^1\ldots   \\
	d_1^2d_2^2\ldots  \\
	\overline{d_1^\+}\,\overline{d_2^\+}\ldots 
\end{array}\right).
\]
This representation plays an important role in our proof of Theorem \ref{th:main-result}; and it will be clear after we {\color{black}prove} Proposition \ref{prop:lexicographical-characterization-U}.
Note by (\ref{eq:sum-equal}) that each row in the above presentation is a length $n$ word or an infinite sequence over the alphabet  $\set{0,1,\ldots, M}$.

Let  $\theta: \Om_M\to\Om_M$ be a cyclic permutation map defined as follows. For $d=(d^1, d^2)\in\Omega_M$ set 
\begin{equation}\label{eq:theta}
	\theta(d)=\theta((d^1, d^2))=(M-(d^1+d^2), d^1)=(\overline{d^\+}, d^1).
\end{equation}
Then $\theta(\Om_M)=\Om_M$. Furthermore, 
$\theta^2(d)=(d^2,\overline{d^\+})$ and $\theta^3(d)=d$. In other words, by using the matrix form we have
\[\left(
\begin{array}{c}
	d^1 \\
	d^2\\
	\overline{d^\+} 
\end{array}\right)\stackrel{\theta}{\longrightarrow}
\left(
\begin{array}{c}
	\overline{d^\+} \\
	d^1\\
	d^2 
\end{array}\right)\stackrel{\theta^2}{\longrightarrow}
\left(
\begin{array}{c}
	d^2 \\
	\overline{d^\+}\\
	d^1
\end{array}\right)\stackrel{\theta^3}{\longrightarrow}
\left(
\begin{array}{c}
	d^1 \\
	d^2\\
	\overline{d^\+} 
\end{array}\right).
\]
For a block $d_1\ldots d_n\in\Om_M^*$ we define $\theta(d_1\ldots d_n)=\theta(d_1)\ldots \theta(d_n)\in\Om_M^*$; and for a sequence $(d_i)\in\Om_M^\N$ we set $\theta((d_i))=\theta(d_1)\theta(d_2)\ldots\in\Om_M^\N$.

Given $\beta\in(1,M+1)$, for an infinite sequence $(d_i)\in\Omega_M^\N$ we write
\[
\Pi_\beta(d_1d_2\ldots)=\sum_{i=1}^{\f}\frac{d_i}{\beta^i}=\left(\sum_{i=1}^{\f}\frac{d_i^1}{\beta^i}, \sum_{i=1}^{\f}\frac{d_i^2}{\beta^i}\right)=:(\pi_\beta((d_i^1)), \pi_\beta((d_i^2))).
\]
Note that the two induced sequences $(d_i^1)$ and $(d_i^2)$ both have digits in $\set{0,1,\ldots, M}$. Then we also need to introduce the one-dimensional symbolic dynamics $(\set{0,1,\ldots, M}^\N, \si)$, where  $\set{0,1,\ldots, M}^\N$ is the set of all infinite sequences with each digit from $\set{0,1,\ldots, M}$, and $\si$ is the {\color{black}left-shift} map.  For $n\in\N$ let 
$\set{0,1,\ldots, M}^n$ be the set of all length $n$ words with each digit from $\set{0,1,\ldots, M}$. Set $\set{0,1,\ldots, M}^*:=\bigcup_{n=1}^\f\set{0,1,\ldots, M}^n$. For a word $w=w_1\ldots w_n\in\set{0,1,\ldots, M}^*$, if $w_n<M$ then we set $w^+:=w_1\ldots w_{n-1}(w_n+1)$; if $w_n>0$ then we set $w^-:=w_1\ldots w_{n-1}(w_n-1)$.  We also write $w^\f$ for the periodic sequence with {\color{black}the} periodic block $w$. For a sequence $(c_i)\in\set{0,1,\ldots, M}^\N$ we denote its \emph{reflection} by $\overline{(c_i)}=(M-c_1)(M-c_2)\ldots$. Furthermore, we use \emph{lexicographical order} between sequences in $\set{0,1,\ldots, M}^\N$. For two infinite sequences $(c_i), (d_i)\in\set{0,1,\ldots, M}^\N$ we write $(c_i)\prec (d_i)$ if there exists $N\in\N$ such that $c_1\ldots c_{N-1}=d_1\ldots d_{N-1}$ and $c_N<d_N$. We also write $(c_i)\lle (d_i)$ if $(c_i)\prec (d_i)$ or $(c_i)=(d_i)$. Similarly, we write $(c_i)\succ (d_i)$   if $(d_i)\prec (c_i)$; and write $(c_i)\lge (d_i)$ if $(d_i)\lle(c_i)$.

For $\beta\in(1,M+1)$  
let $\de(\beta)=\de_1(\beta)\de_2(\beta)\ldots \in\set{0,1,\ldots, M}^\N$ be the \emph{quasi-greedy} $\beta$-expansion of $1$, i.e., the lexicographically largest sequence $(\de_i)$ in $\set{0,1,\ldots, M}^\N$ not ending with $0^\f$ and satisfying $\pi_\beta((\de_i))=1$ (cf.~\cite{Daroczy_Katai_1993}). The following characterization of $\de(\beta)$ was proven in \cite{Baiocchi_Komornik_2007}.  
\begin{lemma}
	\label{lem:char-quasi-greedy-expansion}
	The map $\beta\mapsto \de(\beta)$ is strictly increasing and bijective from $(1,M+1]$ to the set of sequences $(\de_i)\in\{0,1,\ldots, M\}^\N$ satisfying
	\[
	0^\f\prec\si^n((\de_i))\lle (\de_i)\quad\forall n\ge 0.
	\]
	Furthermore, the inverse map $\de^{-1}$ is continuous with respect to the order topology. 
\end{lemma}

In general, for $\beta\in(1,M+1)$ and $x\in(0,1]$ let   $a(x,\beta)$ be the \emph{quasi-greedy} $\beta$-expansion of $x$, i.e., the lexicographically largest sequence $(a_i) \in\set{0,1,\ldots, M}^\N$ not ending with $0^\f$ and satisfying $\pi_\beta((a_i))=x$. The following result for quasi-greedy expansions was established by Parry \cite{Parry_1960} and de Vries et al.~\cite[Lemma 2.2 and Proposition 2.4]{Vries-Komornik-Loreti-2022}. 
\begin{lemma}
	\label{lem:char-quasi-greedy-expansion-x}
	Let $\beta\in(1, M+1]$. The map $x\mapsto a(x,\beta)$ is strictly  increasing and  bijective from $(0, 1]$ to the set 
	\[
	\set{(a_i)\in\set{0,1,\ldots, M}^\N: 0^\f\prec \si^n((a_i))\lle \de(\beta)~\forall n\ge 0}.
	\]
	Furthermore, the map $x\mapsto a(x,\beta)$ is left-continuous with respect to the order topology. 
\end{lemma}

For $\beta\in(1, M+1)$ we recall from (\ref{eq:U-beta}) and (\ref{eq:U-beta-intrinsic}) the definitions of the intrinsic  univoque set $\widetilde{U}_{\beta, M}$ and   univoque set ${U}_{\beta,M}$, respectively. 
According to the coding map $\Pi_\beta$ in (\ref{eq:coding-map}), we define the symbolic intrinsic univoque set $\widetilde{\U}_{\beta, M}\subset \Om_M^\N$ and the symbolic univoque set $\U_{\beta, M}\subset\Om_M^\N$ by
\[
\widetilde{\U}_{\beta, M}:=\Pi_\beta^{-1}(\widetilde{U}_{\beta, M})\quad\textrm{and}\quad  \U_{\beta, M}:=\Pi_\beta^{-1}(U_{\beta, M}).
\]Note that $\widetilde{U}_{\beta, M}\subset U_{\beta, M}$ for any $\beta\in(1, M+1)$, and each $x\in U_{\beta, M}$ has a unique coding under $\Pi_\beta$.
Then  $\Pi_\beta$ is  bijective  from $\widetilde{\U}_{\beta, M}$ to $\widetilde{U}_{\beta, M}$, and from $\U_{\beta, M}$ to $U_{\beta, M}$  as well.

The following   characterization of $\widetilde{\U}_{\beta, M}$   plays an important role in the proof of Theorem \ref{th:main-result}. When  $M=1$, it   was proven in \cite[Proposition 2.2]{Kong-Li-2020}.

\begin{proposition}
	\label{prop:lexicographical-characterization-U}
	Let $M\in\N$ and $\beta\in(1,M+1)$. Then     $(d_i)\in\widetilde{\U}_{\beta, M}$ if and only if   $(d_i)\in\Omega_M^\N$ satisfies
	\begin{equation}\label{eq:chara-U-beta}
		\left\{
		\begin{array}{ccc}
			d_{n+1}^1d_{n+2}^1 \ldots \prec \de(\beta)& \textrm{whenever} & d_n^1<M, \\
			d_{n+1}^2d_{n+2}^2 \ldots \prec \de(\beta)& \textrm{whenever} & d_n^2<M, \\
			\overline{d_{n+1}^\+}\overline{d_{n+2}^\+} \ldots \prec \de(\beta)& \textrm{whenever} & \overline{d_n^\+}<M.
		\end{array}
		\right.
	\end{equation} 
\end{proposition}

Before proving the proposition  we clarify our main idea.  Note that $d_n^1<M$ implies that the first coordinate of $\Pi_\beta(d_nd_{n+1}\ldots)\in U_{\beta,M}$ is strictly smaller than $\frac{d_n^1+1}{\beta}$. Similar to the one-dimensional unique $\beta$-expansions we can deduce that this is equivalent to $d_{n+1}^1d_{n+2}^1\ldots\prec \de(\beta)$. Similar ideas can be applied to study the case for $d_n^2<M$. However, the case for $\overline{d_n^\+}<M$ is more involved. Note that $\overline{d_n^\+}=M-(d_n^1+d_n^2)$. Then $\overline{d_n^\+}<M$ is equivalent to $d_n^1+d_n^2>0$. Observe that $d_n\in\Omega_M$, i.e., $d_n^1, d_n^2\in\set{0,1,\ldots, M}$ and $d_n^1+d_n^2\le M$. Similar to the previous two cases, we split the proof of the case for $\overline{d_n^\+}<M$ into four subcases: (i) $d_n^1<M, d_n^2<M$ and $d_n^1+d_n^2>0$; (ii) $d_n^1+d_n^2=0$, i.e., $d_n^1=d_n^2=0$; (iii) $d_n^1=0$ and $d_n^2=M$; (iv) $d_n^1=M$ and $d_n^2=0$. Each subcase can be proved similar to the one-dimensional unique $\beta$-expansions.

\begin{proof}
	First  we prove the necessity. Take  $(d_i)\in \widetilde{\U}_{\beta, M}$.  We will prove (\ref{eq:chara-U-beta}) by considering  the following three cases: (I) $d_n^1<M$; (II) $d_n^2<M$; (III) $\overline{d_n^\+}<M$ for some $n\in\N$.

	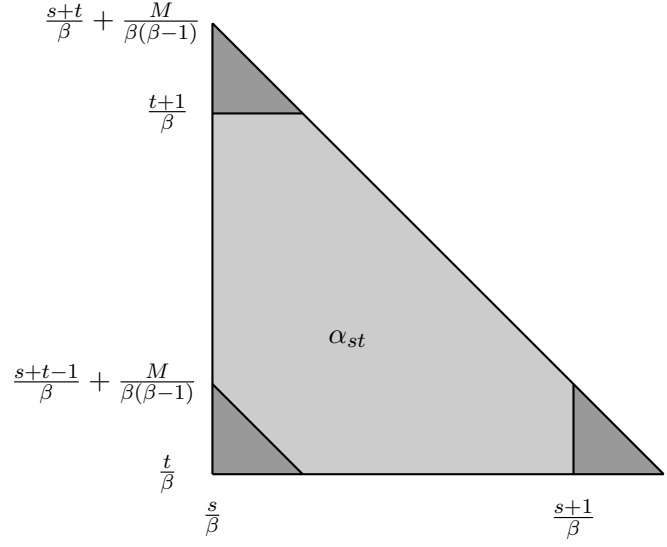
\begin{figure}[h!]
		\begin{center}
			\begin{tikzpicture}[
				scale=12,
				axis/.style={very thick, ->},
				important line/.style={thick},
				dashed line/.style={dashed, thin},
				pile/.style={thick, ->, >=stealth', shorten <=2pt, shorten
					>=2pt},
				every node/.style={color=black}
				]
				\fill[black!20](0,0)--(0.5,0)--(0,0.5)--cycle;
				\fill[black!40](0,0)--(0.1,0)--(0,0.1)--cycle;
				
				\fill[black!40](0.4,0)--(0.5,0)--(0.4,0.1)--cycle;
				
				\fill[black!40](0,0.4)--(0.1,0.4)--(0,0.5)--cycle;
				
				\draw[important line] (0,0)--(0.5,0);
				
				\draw[important line] (0, 0)--(0,0.5);
				
				\draw[important line] (0.5, 0)--(0,0.5);
				
				\draw[important line] (0.1,0)--(0,0.1);
				\draw[important line] (0.4,0)--(0.4,0.1);
				\draw[important line] (0,0.4)--(0.1,0.4);
				
				\node[] at (0,-0.05){$\frac{s}{\beta}$};
				\node[] at (-0.12, 0.1){$\frac{s+t-1}{\beta}+\frac{M}{\beta(\beta-1)}$};
				\node[] at (0.4,-0.05){$\frac{s+1}{\beta}$};
				\node[] at (-0.1, 0.5){$\frac{s+t}{\beta}+\frac{M}{\beta(\beta-1)}$};
				\node[] at (-0.05,0){$\frac{t}{\beta}$};
				
				\node[] at (-0.05,0.4){$\frac{t+1}{\beta}$};
				
				\node[] at (0.15, 0.15){$\al_{st}$};

			\end{tikzpicture} 
		\end{center}
		\caption{The typical  pattern  of $f_{\beta, d}(\Delta_{\beta,M})\setminus O_{\beta,M}$ with   $d=(d^1, d^2)\in\Omega_M$ satisfying $d^1<M, d^2<M$ and $\overline{d^\+}=M-(d^1+d^2)<M$.}\label{fig:1}
	\end{figure}
	Case (I). Suppose $d_n^1=s<M$  for some $n$. Then $d_n=\al_{sj}$ for some $j\in\set{0,1,\ldots,M}$. Since $(d_i)\in\widetilde{\U}_{\beta, M}=\Pi_\beta^{-1}(\widetilde{U}_{\beta, M})$,  by the definition of $\widetilde{U}_{\beta, M}$ in (\ref{eq:U-beta-intrinsic}) it follows that
	\[\Pi_\beta(d_nd_{n+1}\ldots) \in\bigcup_{j=0}^M f_{\beta,\al_{s j}}(\Delta_{\beta,M})\setminus {O}_{\beta,M}.\]
	In view of Figure \ref{Fig:S-4}, there are four possible patterns for $f_{\beta, d}(\Delta_{\beta,M})\setminus O_{\beta,M}$ with $d\in\Omega_M$: see  Figure \ref{fig:1}  for the typical pattern, and see Figure  \ref{fig:2} for the patterns in the three corners. Since $d_n^1=s<M$, this implies $\pi_\beta(d_n^1d_{n+1}^1\ldots)<\frac{s+1}{\beta}$. Note that $d_n^1=s$. Then 
	\begin{equation}\label{eq:dec31-1}\pi_\beta(d_{n+1}^1d_{n+2}^1\ldots)<1=\pi_\beta(\delta(\beta)).\end{equation}
	So, it suffices to prove $d_{n+1}^1d_{n+2}^1\ldots \prec \de(\beta)$. 
	
	Suppose on the contrary that $d_{n+1}^1d_{n+2}^1\ldots \lge\de(\beta)$. If the equality holds, then we have $\pi_\beta(d_{n+1}^1d_{n+2}^1\ldots)=1$, leading to a contradiction to (\ref{eq:dec31-1}). If strict inequality holds, then there exists $k\in\N$ such that $d_{n+1}^1\ldots d_{n+k-1}^1=\de_1(\beta)\ldots \de_{k-1}(\beta)$ and $d_{n+k}^1>\de_{k}(\beta)$. By Lemmas \ref{lem:char-quasi-greedy-expansion} and \ref{lem:char-quasi-greedy-expansion-x} it follows that
	\begin{align*}
		1 = \pi_\beta(\de(\beta))    &=\sum_{i=1}^{k}\frac{\de_i(\beta)}{\beta^i}+\frac{1}{\beta^k}\sum_{i=1}^{\f}\frac{\de_{k+i}(\beta)}{\beta^i}\\
		&\le\sum_{i=1}^{k}\frac{\de_i(\beta)}{\beta^i}+\frac{1}{\beta^k}\\
		&\le\sum_{i=1}^{k}\frac{d_{n+i}^1}{\beta^i}\le \pi_\beta(d_{n+1}^1d_{n+2}^1\ldots),
	\end{align*}
	again leading to a contradiction to (\ref{eq:dec31-1}). This proves $d_{n+1}^1d_{n+2}^1\ldots\prec \de(\beta)$ as required.

	Case (II). Suppose $d_n^2=t<M$  for some $n$. Similar to Case (I) we can prove that $d_{n+1}^2d_{n+2}^2\ldots \prec \de(\beta)$.
	
	Case (III). Suppose $\overline{d_n^\+}=M-(s+t)<M$  for some $n$, where $d_n^1=s$ and $d_n^2=t$.  Then $d_n^\+=d_n^1+d_n^2=M-\overline{d_n^\+}=s+t>0$. Since $(d_i)\in\widetilde{\U}_{\beta, M}=\Pi_\beta^{-1}(\widetilde{U}_{\beta, M})$,  by (\ref{eq:U-beta-intrinsic}) we obtain that (see Figure     \ref{fig:1})
	\[
	\pi_\beta(d_n^\+d_{n+1}^\+\ldots)=\pi_\beta(d_n^1d_{n+1}^1\ldots) +\pi_\beta(d_n^2d_{n+1}^2\ldots) >\frac{s+t-1}{\beta}+\frac{M}{\beta(\beta-1)}.
	\] 
	Note that $d_n^\+=d_n^1+ d_n^2=s+t$. Then 
	$
	\pi_\beta(d_{n+1}^\+d_{n+2}^\+\ldots) >\frac{M}{\beta-1}-1,
	$
	which implies 
	\[
	\pi_\beta(\overline{d_{n+1}^\+ d_{n+2}^\+\ldots}) <1=\pi_\beta(\de(\beta)).
	\]
	By the same argument as in Case (I) we conclude  $\overline{d_{n+1}^\+d_{n+2}^\+}\prec \de(\beta)$.

	Next we prove the sufficiency. Let $(d_i)\in\Omega_M^\N$ be an infinite  sequence satisfying (\ref{eq:chara-U-beta}). Then it suffices to prove 
	\begin{equation}\label{eq:char-1}
		\Pi_\beta(d_nd_{n+1}\ldots)\in f_{\beta, d_n}(\Delta_{\beta,M})\setminus  O_{\beta, M}\quad\forall n\ge 1.
	\end{equation}
	Take   $d_n=(s,t)\in\Om_M$ for some $n\in\N$. In view of Figures \ref{Fig:S-4} and \ref{fig:1}, we split the proof of (\ref{eq:char-1}) into the following four cases.
	
	Case A. $s<M, t<M$ and $s+t>0$. Then $d_n^1=s<M, d_n^2=t<M$ and $\overline{d_n^\+}=M-(s+t)<M$. First we claim that 
	\begin{equation}\label{eq:jan14-1}
		\pi_\beta(d_n^1d_{n+1}^1\ldots)<\frac{s+1}{\beta}.
	\end{equation}
	Note that $d_n^1=s<M$. Starting with $k_0=n$ we define by (\ref{eq:chara-U-beta}) an infinite sequence of indices $k_0<k_1<k_2<\cdots$ satisfying the conditions
		\[
		d_{k_{j-1}+1}^1\ldots d_{k_j-1}^1=\de_1(\beta)\ldots\de_{k_j-k_{j-1}-1}(\beta)\quad\textrm{and}\quad d_{k_j}^1<\de_{k_j-k_{j-1}}(\beta)
		\]
		for $j=1,2,\ldots$.
	Then by Lemma \ref{lem:char-quasi-greedy-expansion} it follows that
	\begin{align*}
		\pi_\beta(d_n^1d_{n+1}^1\ldots)&=\frac{s}{\beta}+\beta^{n-1}\sum_{i=n+1}^{\f}\frac{d_i^1}{\beta^i} =\frac{s}{\beta}+\beta^{n-1}\sum_{j=1}^{\f}\sum_{i=1}^{k_j-k_{j-1}}\frac{d_{k_{j-1}+i}^1}{\beta^{k_{j-1}+i}}  \\
		\le& \frac{s}{\beta}+\beta^{n-1}\sum_{j=1}^{\f}\left(\sum_{i=1}^{k_j-k_{j-1}}\frac{\de_i(\beta)}{\beta^{k_{j-1}+i}}-\frac{1}{\beta^{k_j}}\right) \\
		< & \frac{s}{\beta}+\beta^{n-1}\sum_{j=1}^{\f}\left(\frac{1}{\beta^{k_{j-1}}}-\frac{1}{\beta^{k_j}}\right)=\frac{s}{\beta}+\frac{\beta^{n-1}}{\beta^{k_0}}=\frac{s+1}{\beta}.
	\end{align*}
	This proves (\ref{eq:jan14-1}). Similarly, since $d_{n}^2=t<M$ and $\overline{d_n^\+}=M-(s+t)<M$, by (\ref{eq:chara-U-beta}) and the same argument as above we can prove 
	\[
	\pi_\beta(d_n^2d_{n+1}^2\ldots)<\frac{t+1}{\beta}\quad\textrm{and}\quad \pi_\beta(\overline{d_n^\+d_{n+1}^\+\ldots})<\frac{s+t}{\beta};
	\]
	and the second inequality implies $\pi_\beta(d_n^\+d_{n+1}^\+\ldots)>\frac{s+t-1}{\beta}+\frac{M}{\beta(\beta-1)}$.
	Therefore, in view of Figure  \ref{fig:1}, this proves  (\ref{eq:char-1}).
	
	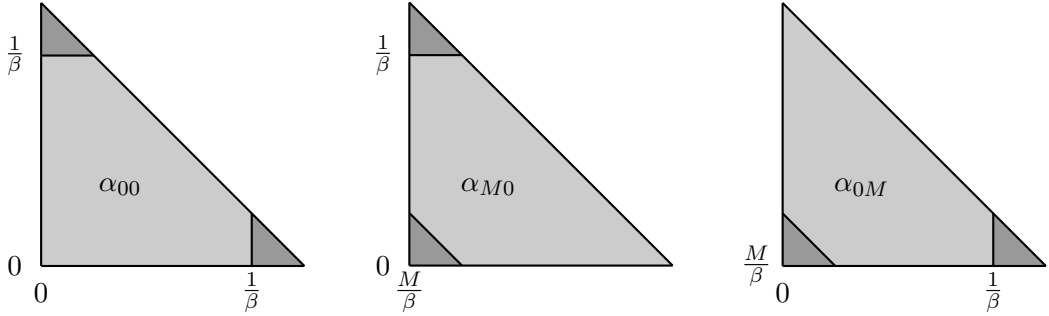
\begin{figure}[h!]
		\centering
		\subfloat[$d=(0, 0)$ \label{fig:sub1}]{
			\begin{minipage}[b]{0.3\textwidth} 
				\centering
				
				\begin{tikzpicture}[
					scale=7,
					axis/.style={very thick, ->},
					important line/.style={thick},
					dashed line/.style={dashed, thin},
					pile/.style={thick, ->, >=stealth', shorten <=2pt, shorten
						>=2pt},
					every node/.style={color=black}
					]
					\fill[black!20](0,0)--(0.5,0)--(0,0.5)--cycle;
					\fill[black!40](0.4,0)--(0.5,0)--(0.4,0.1)--cycle;
					
					\fill[black!40](0,0.4)--(0.1,0.4)--(0,0.5)--cycle;            
					\draw[important line] (0,0)--(0.5,0);
					
					\draw[important line] (0, 0)--(0,0.5);
					
					\draw[important line] (0.5, 0)--(0,0.5);
					
					\draw[important line] (0.4,0)--(0.4,0.1);
					\draw[important line] (0,0.4)--(0.1,0.4);
					
					\node[] at (0,-0.05){$0$};
					\node[] at (0.4,-0.05){$\frac{1}{\beta}$};
					
					\node[] at (-0.05,0){$0$};
					
					\node[] at (-0.05,0.4){$\frac{1}{\beta}$};
					
					\node[] at (0.15, 0.15){$\al_{00}$};
				\end{tikzpicture}
			\end{minipage}
		}
		\subfloat[$d=(M, 0)$ \label{fig:sub2}]{
			\begin{minipage}[b]{0.3\textwidth}
				\centering
				\begin{tikzpicture}[
					scale=7,
					axis/.style={very thick, ->},
					important line/.style={thick},
					dashed line/.style={dashed, thin},
					pile/.style={thick, ->, >=stealth', shorten <=2pt, shorten
						>=2pt},
					every node/.style={color=black}
					]
					\fill[black!20](0,0)--(0.5,0)--(0,0.5)--cycle;
					\fill[black!40](0,0)--(0.1,0)--(0,0.1)--cycle;
					
					
					\fill[black!40](0,0.4)--(0.1,0.4)--(0,0.5)--cycle;                
					\draw[important line] (0,0)--(0.5,0);
					
					\draw[important line] (0, 0)--(0,0.5);
					
					\draw[important line] (0.5, 0)--(0,0.5);
					
					\draw[important line] (0.1,0)--(0,0.1);
					\draw[important line] (0,0.4)--(0.1,0.4);
					
					\node[] at (0,-0.05){$\frac{M}{\beta}$};
					
					\node[] at (-0.05,0){$0$};
					
					\node[] at (-0.05,0.4){$\frac{1}{\beta}$};
					
					\node[] at (0.15, 0.15){$\al_{M0}$};
				\end{tikzpicture}
			\end{minipage}
		}
		\subfloat[$d=(0, M)$ \label{fig:sub3}]{
			\begin{minipage}[b]{0.3\textwidth}
				\centering
				\begin{tikzpicture}[
					scale=7,
					axis/.style={very thick, ->},
					important line/.style={thick},
					dashed line/.style={dashed, thin},
					pile/.style={thick, ->, >=stealth', shorten <=2pt, shorten
						>=2pt},
					every node/.style={color=black}
					]
					\fill[black!20](0,0)--(0.5,0)--(0,0.5)--cycle;
					\fill[black!40](0,0)--(0.1,0)--(0,0.1)--cycle;
					
					\fill[black!40](0.4,0)--(0.5,0)--(0.4,0.1)--cycle;
					
					
					\draw[important line] (0,0)--(0.5,0);
					
					\draw[important line] (0, 0)--(0,0.5);
					
					\draw[important line] (0.5, 0)--(0,0.5);
					
					\draw[important line] (0.1,0)--(0,0.1);
					\draw[important line] (0.4,0)--(0.4,0.1);
					
					\node[] at (0,-0.05){$0$};
					\node[] at (0.4,-0.05){$\frac{1}{\beta}$};
					
					\node[] at (-0.05,0){$\frac{M}{\beta}$};
					
					
					\node[] at (0.15, 0.15){$\al_{0M}$};
				\end{tikzpicture}
			\end{minipage}
		}
		\caption{The three patterns of $f_{\beta, d}(\Delta_{\beta,M})\setminus O_{\beta,M}$ with $d=(0,0), d=(M,0)$ and $d=(0,M)$  respectively.}
		\label{fig:2}
	\end{figure}

	Case B. $s+t=0$. Then $d_n^1=d_n^2=0$ and $\overline{d_n^\+}=M$. By the same argument as in Case A we can prove that 
	\[
	\pi_\beta(d_n^1d_{n+1}^1\ldots) <\frac{1}{\beta},\quad \pi_\beta(d_n^2d_{n+1}^2\ldots) <\frac{1}{\beta}.
	\]	
	In view of Figure  \ref{fig:sub1} ,  we can also prove (\ref{eq:char-1}).
	
	Case C. $s=M$. Then $d_n^1=M$, $d_n^2=0$ and $\overline{d_n^\+}=0$. Similar to Case A, we have 
	\[
	\pi_\beta(d_n^2d_{n+1}^2\ldots) <\frac{1}{\beta},\quad \pi_\beta(d_n^\+d_{n+1}^\+\ldots) >\frac{M-1}{\beta}+\frac{M}{\beta(\beta-1)}.
	\] 
	In view of Figure  \ref{fig:sub2}, this yields (\ref{eq:char-1}).
	
	Case D. $t=M$. Then $d_n^1=0$, $d_n^2=M$ and $\overline{d_n^\+}=0$. By the same argument as in Case A we obtain that 
	\[
	\pi_\beta(d_n^1d_{n+1}^1\ldots) <\frac{1}{\beta},\quad \pi_\beta(d_n^\+d_{n+1}^\+\ldots) >\frac{M-1}{\beta}+\frac{M}{\beta(\beta-1)}.
	\]
	In view of Figure  \ref{fig:sub3}, this also proves (\ref{eq:char-1}). 
\end{proof}
\begin{remark}
	\label{rem:chara-U-beta}
	Note by (\ref{eq:theta}) that $\theta$ is the cyclic permutation. Then 
	by Proposition \ref{prop:lexicographical-characterization-U} it follows that 
	\begin{equation}\label{eq:cyclic-permutation}
		(d_i)\in\widetilde{\U}_{\beta, M}\quad\Longleftrightarrow\quad \theta((d_i))\in\widetilde{\U}_{\beta, M}\quad\Longleftrightarrow\quad\theta^2((d_i))\in\widetilde{\U}_{\beta, M}.
	\end{equation}
\end{remark}

At the end of this section we show that for $\beta\le (M+2)/2$ the univoque set $\U_{\beta, M}$ coincides with the intrinsic univoque set $\widetilde{\U}_{\beta,M}$.
\begin{proposition}
	\label{prop:U-beta-coincidence}
	Let $M\in\N$ and $\beta\in(1, M+1)$. Then  $S_{\beta, M}=\Delta_{\beta, M}$ if and only if $\beta\in(1, \frac{M+2}{2}]$. Hence, for any $\beta\in(1, \frac{M+2}{2}]$ we have 
	\[
	U_{\beta, M}=\widetilde{U}_{\beta, M}\quad\textrm{and}\quad \U_{\beta,M}=\widetilde{\U}_{\beta, M}.
	\]
\end{proposition}
\begin{proof}
	In view of Figures \ref{Fig:S-4} and \ref{fig:1}, the fat Sierpinski gasket $S_{\beta, M}$ is equal to its convex hull $\Delta_{\beta, M}$ if and only if for any $d=(s, t)\in\N_0^2$ with $s+t\le M$ we must have 
	\[
	\frac{s+1}{\beta}+\frac{t+1}{\beta}\le \frac{s+t}{\beta}+\frac{M}{\beta(\beta-1)},
	\]
	which is equivalent to $\beta\le \frac{M+2}{2}$. This proves the first statement. 
	
	Suppose $\beta\in(1, \frac{M+2}{2}]$. Then $S_{\beta, M}=\Delta_{\beta, M}$, and in this case each point in the overlap region $O_{\beta, M}$ has at least two codings (cf.~\cite[Proposition 2.4]{Kong-Li-2020}). So, by (\ref{eq:U-beta}) and (\ref{eq:U-beta-intrinsic}) it follows  that   $\widetilde{U}_{\beta,M}=U_{\beta, M}$, and also  $\widetilde{\U}_{\beta, M}=\U_{\beta, M}$   since $\Pi_\beta$ is  bijective.    
\end{proof}
We emphasize that when $M\ge 2$ the critical bases $\beta_G(M)$ and $\beta_c(M)$ defined in Definitions \ref{def:beta-G} and \ref{def:beta-c} are both strictly smaller than $\frac{M+2}{2}$ (see Figure \ref{fig:critical-values} for $M\le 10$). So, Proposition \ref{prop:U-beta-coincidence} implies that the critical bases for the phase transitions of $\widetilde{U}_{\beta, M}$ and $U_{\beta, M}$ are the same. Hence, in the following we only prove Theorem \ref{th:main-result} for $\widetilde{U}_{\beta, M}$.

\section{Generalized golden ratio $\beta_G(M)$ for a general $M$}\label{sec:golden-ratio}
Recall from Definition \ref{def:beta-G} the generalized golden ratio $\beta_G(M)$. Then by Lemma \ref{lem:char-quasi-greedy-expansion} it follows that  the quasi-greedy $\beta_G(M)$-expansion of $1$ is given by
\begin{equation}\label{eq:dec31-2}
	\de(\beta_G(M))=\left\{
	\begin{array}{lll}
		((N+1)NN)^\f & \textrm{if} & M=3N+1,\\
		((N+1)(N+1)N)^\f&\textrm{if}& M=3N+2,\\
		(N+1)^\f & \textrm{if} & M=3N+3.
	\end{array}
	\right.
\end{equation} The following result shows that $\beta_G(M)$ is the first critical base for $\widetilde{\U}_{\beta, M}$, i.e., $\widetilde{\U}_{\beta, M}$ is finite if and only if $\beta\le \beta_G(M)$. When $M=1$, this was proven by  Sidorov \cite[Theorem 4.1]{Sidorov_2007} (see also, \cite{Kong-Li-2020}).

\begin{proposition}\label{prop:beta-G}
	Let $\beta_G(M)$ be defined as in Definition \ref{def:beta-G}. Then $\beta_G(M)$ is a Perron number, and 
	\[\#\widetilde{\U}_{\beta, M}<+\f\quad\Longleftrightarrow\quad \beta\le \beta_G(M).\]
	Furthermore, for $\beta\in(1, \beta_G(M)]$ we have 
	\[
	\widetilde{\U}_{\beta, M}=\bigcup_{j=0}^2\set{\theta^j(\al_{00}^\f)}=\set{\al_{00}^\f, \al_{0M}^\f, \al_{M0}^\f}.
	\]
\end{proposition} 
\begin{proof}
	By (\ref{eq:dec31-2}) it follows that $\beta_G(M)$ is a Perron number. For the remaining statements we consider three cases: (I) $M=3N+1$; (II) $M=3N+2$;  and (III) $M=3N+3$ for some $N\in\N_0$.
	
	Case (I).   $M=3N+1$ with $N\in\N_0$. Then by (\ref{eq:dec31-2}), Lemma \ref{lem:char-quasi-greedy-expansion} and Proposition \ref{prop:lexicographical-characterization-U} we only need to prove that $\widetilde{\U}_{\beta, M}=\set{\al_{00}^\f, \al_{M0}^\f, \al_{0M}^\f}$ if $\de(\beta)\lle ((N+1)NN)^\f$; and $\#\widetilde{\U}_{\beta, M}=\f$ if $\de(\beta)\succ ((N+1)NN)^\f$.
	
	Take $\beta\in(1, M+1]$ so that  $\de(\beta)\lle ((N+1)NN)^\f$. Let $(d_i)\in\widetilde{\U}_{\beta, M}$. Recall that each $d_i=(d_i^1, d_i^2)$ and $\overline{d_i^\+}=M-d_i^\+$ with $d_i^\+=d_i^1+d_i^2$. Suppose $d_1^1, d_1^2, \overline{d_1^\+}<M$. Then by Proposition \ref{prop:lexicographical-characterization-U} and (\ref{eq:sum-equal}) it follows that 
	\begin{equation}\label{eq:3n+1-case2}
		d_2\in\set{\al_{(N+1)(N+1)}, \al_{(N-1)(N+1)}, \al_{(N+1)(N-1)}}=\set{\theta^j(\al_{(N+1)(N+1)}): j=0,1,2}
	\end{equation}
	or
	\begin{equation}\label{eq:3n+1-case1}
		d_2\in\set{\al_{(N+1)N}, \al_{N(N+1)}, \al_{NN}}=\set{\theta^j(\al_{(N+1)N}): j=0,1,2}.
	\end{equation}
	For (\ref{eq:3n+1-case2}) we assume by (\ref{eq:cyclic-permutation}) that $d_2=\al_{(N+1)(N+1)}$. Then by Proposition \ref{prop:lexicographical-characterization-U} and (\ref{eq:sum-equal}) we obtain that 
	\[
	d_2d_3d_4 =\al_{(N+1)(N+1)}\al_{NN}\al_{NN},
	\]
	which yields $\overline{d_2^\+d_3^\+d_4^\+}=(N-1)(N+1)(N+1)$. This
	leads to a contradiction to $(d_i)\in\widetilde{\U}_{\beta, M}$.
	For (\ref{eq:3n+1-case1}) we   assume by (\ref{eq:cyclic-permutation}) that $d_2=\al_{(N+1)N}$. Then by Proposition \ref{prop:lexicographical-characterization-U} and (\ref{eq:sum-equal}) it follows that 
	\[
	d_3\in\set{\al_{N(N+1)}, \al_{NN}, \al_{(N-1)(N+1)}}.
	\] 
	If $d_3=\al_{(N-1)(N+1)}$, then by a similar argument as in (\ref{eq:3n+1-case2}) we will obtain that $d_3d_4d_5=\al_{(N-1)(N+1)}\al_{(N+1)N}\al_{(N+1)N}$, leading to a contradiction to $(d_i)\in\widetilde{\U}_{\beta, M}$. So, $$d_2d_3\in\set{\al_{(N+1)N}\al_{N(N+1)}, \al_{(N+1)N}\al_{NN}}.$$ By Proposition \ref{prop:lexicographical-characterization-U} and the same argument as above we can deduce that 
	\[
	\left(\begin{array}{c}
		d_2^1d_3^1\ldots \\
		d_2^2 d_3^2\ldots\\
		\overline{d_2^\+}\overline{d_3^\+}\ldots
	\end{array}\right)\in\set{\left(\begin{array}{c}
			((N+1)NN)^\f\\
			(N(N+1)N)^\f \\
			(NN(N+1))^\f
		\end{array}\right), 
		\left(\begin{array}{c}
			((N+1)NN)^\f\\
			(N N(N+1))^\f \\
			(N(N+1)N)^\f
		\end{array}\right)
	},\]
	which again leads to a contradiction to $(d_i)\in\widetilde{\U}_{\beta, M}$. 
	
	Therefore, $d_1^1=M, d_1^2=M$ or $\overline{d_1^\+}=M$.  By (\ref{eq:cyclic-permutation}) we may assume $d_1^1=M$. Then by (\ref{eq:sum-equal}) we have $d_1^2=\overline{d_1^\+}=0$. By Proposition \ref{prop:lexicographical-characterization-U} and the same argument as above we can prove that $(d_i)=\al_{M0}^\f$. Therefore, by (\ref{eq:cyclic-permutation}) this proves $\widetilde{\U}_{\beta, M}=\set{\al_{00}^\f, \al_{0M}^\f, \al_{M0}^\f}$. 
	
	On the other hand, if $\beta\in(1, M+1]$ with $\de(\beta)\succ ((N+1)NN)^\f$ then  we can prove that 
	\[\al_{M0}^k(\al_{(N+1)N}\al_{N(N+1)}\al_{NN})^\f\in\widetilde{\U}_{\beta, M}\quad\textrm{for any }k\in\N_0.\]
	So $\#\widetilde{\U}_{\beta, M}=\f$.  
	
	Case (II).   $M=3N+2$ with $N\in\N_0$. By (\ref{eq:dec31-2}), Lemma \ref{lem:char-quasi-greedy-expansion} and Proposition \ref{prop:lexicographical-characterization-U} it suffices to prove that $\widetilde{\U}_{\beta, M}=\set{\al_{00}^\f, \al_{M0}^\f, \al_{0M}^\f}$ if $\de(\beta)\lle ((N+1)(N+1)N)^\f$; and $\#\widetilde{\U}_{\beta, M}=\f$ if $\de(\beta)\succ ((N+1)(N+1)N)^\f$.
	
	Take $\beta\in(1, M+1]$ with $\de(\beta)\lle ((N+1)(N+1)N)^\f$, and let $(d_i)\in\widetilde{\U}_{\beta, M}$. If $d_1^1, d_1^2, \overline{d_1^\+}<M$, then by Proposition \ref{prop:lexicographical-characterization-U} and (\ref{eq:sum-equal}) it follows that  
	\[
	d_2\in\set{\al_{(N+1)(N+1)}, \al_{(N+1)N}, \al_{N(N+1)}}=\set{\theta^j(\al_{(N+1)(N+1)}): j=0,1,2}.
	\]
	Then by (\ref{eq:cyclic-permutation}) we may assume $d_2=\al_{(N+1)(N+1)}$. By Proposition \ref{prop:lexicographical-characterization-U} and (\ref{eq:sum-equal}) we obtain that 
	\[
	d_2d_3\ldots=(\al_{(N+1)(N+1)}\al_{(N+1)N}\al_{N(N+1)})^\f\quad\textrm{or}\quad (\al_{(N+1)(N+1)}\al_{N(N+1)}\al_{(N+1)N})^\f.
	\]
	Both cases will lead to a contradiction to $(d_i)\in\widetilde{\U}_{\beta, M}$. Similar to  Case $M=3N+1$, we can deduce that  $\widetilde{\U}_{\beta, M}=\set{\al_{00}^\f, \al_{M0}^\f, \al_{0M}^\f}$. 
	On the other hand, if $\beta\in(1, M+1]$ with $\de(\beta)\succ ((N+1)(N+1)N)^\f$, then we can prove that 
	\[\al_{M0}^k(\al_{(N+1)(N+1)}\al_{N(N+1)}\al_{(N+1)N})^\f\in\widetilde{\U}_{\beta, M}\quad\textrm{for any }k\in\N_0.\]
	So $\#\widetilde{\U}_{\beta, M}=\f$.  
	
	Case (III).      $M=3N+3$ with $N\in\N_0$. By (\ref{eq:dec31-2}), Lemma \ref{lem:char-quasi-greedy-expansion} and Proposition \ref{prop:lexicographical-characterization-U} it suffices to prove that $\widetilde{\U}_{\beta, M}=\set{\al_{00}^\f, \al_{0M}^\f, \al_{M0}^\f}$ if $\de(\beta)\lle (N+1)^\f$, and $\#\widetilde{\U}_{\beta, M}=\f$ if $\de(\beta)\succ (N+1)^\f$.  
	
	Take $\beta\in(1, M+1]$ such that  $\de(\beta)\lle (N+1)^\f$, and take $(d_i)\in\widetilde{\U}_{\beta, M}$. Note that $d_1\in\Omega_M$. If $d_1^1, d_1^2, \overline{d_1^\+}<M=3(N+1)$, then by Proposition \ref{prop:lexicographical-characterization-U} and (\ref{eq:sum-equal}) it follows that $d_2^1=d_2^2=\overline{d_2^\+}=N+1$. Since $N+1<M$, by the same argument we can deduce that 
	\[
	d_2^1d_3^1\ldots=d_2^2d_3^2\ldots=\overline{d_2^\+}\overline{d_3^\+}\ldots=(N+1)^\f\lge \de(\beta),
	\]
	leading to a contradiction to $(d_i)\in\widetilde{\U}_{\beta, M}$. So, $d_1^1=M, d_1^2=M$ or $\overline{d_1^\+}=M$. By the same argument as that for $M=3N+1$ we can prove  $\widetilde{\U}_{\beta, M}=\set{\al_{00}^\f, \al_{0M}^\f, \al_{M0}^\f}$.
	On the other hand, take $\beta\in(1, M+1]$ so that $\de(\beta)\succ (N+1)^\f$. Then one can verify by Proposition \ref{prop:lexicographical-characterization-U} that 
	\[
	\al_{M0}^k\al_{(N+1)(N+1)}^\f\in\widetilde{\U}_{\beta, M}\quad\forall k\in\N_0,
	\]
	which implies $\#\widetilde{\U}_{\beta, M}=\f$. This completes the proof.
\end{proof}

\section{Generalized Komornik-Loreti constant $\beta_c(M)$ for $M=3N+1$}\label{sec:3N+1}
Let $M=3N+1$ with $N\in\N_0$. Recall from Definition \ref{def:beta-c} the generalized Komornik-Loreti constant $\beta_c(M)$. In this section we will show that $\dim_H\widetilde{\U}_{\beta, M}>0$ if and only if $\beta>\beta_c(M)$; and  prove   Theorem \ref{th:main-result} for $M=3N+1$.
When $M=1$, Theorem \ref{th:main-result} was established by Kong and Li \cite{Kong-Li-2020}. {In this section  we may assume $M=3N+1$ with $N\in\N$. The new challenge is that for $M=3N+1$ with $N\ge 1$ the digit $N-1$ might occur  in $d_{n}^1, d_n^2$ or $\overline{d_n^\+}$ for a sequence $(d_i)\in\widetilde{\U}_{\beta_c, M}$, while for $M=1$ this can never happen. Remember that each digit $d_i=(d_i^1, d_i^2)$ and $\overline{d_i^\+}=M-d_i^\+$ with $d_i^\+=d_i^1+d_i^2$. Fortunately, we show in Lemma \ref{lem:exceptional-digits} that this digit $N-1$ can occur at most 3 times in any sequence of $\widetilde{\U}_{\beta_c, M}$.}

\subsection{Properties of $\beta_c(M)$}
Recall that $(\tau_i)_{i=0}^\f=0110100110010110\ldots\in\set{0,1}^\N$ is the classical Thue-Morse sequence (cf.~\cite{Allouche_Shallit_1999}). Then $\tau_{2i+1}=1-\tau_i$ and $\tau_{2i}=\tau_i$ for all $i\ge 0$. Now we introduce a
Thue-Morse  type sequence $(\la_i)\in\set{N, N+1}^\N$ such that  
\begin{equation}\label{eq:relation-lambda-tau}
	\la_{3k+1}=\tau_{2k+1}+N,\quad \la_{3k+2}=N, \quad \la_{3k+3}=\tau_{2k+2}+N
\end{equation}
for all $k\ge 0$. Then $(\la_i)$ begins with $(N+1)N(N+1)NN(N+1)NNN(N+1)N(N+1)$. 

The Thue-Morse type sequence  $(\la_i)$ can also be obtained  by  substitutions. 

\begin{lemma}\label{lem-obs-(3N+1)}
	Let $(\tau_i)_{i=0}^\f$ be the classical Thue-Morse sequence. Then 
	\[(\la_i)_{i=1}^\infty=\si(\varphi_1(\tau_0)\varphi_1(\tau_1)\varphi_1(\tau_2)\cdots),\]
	where $\si$ is the left-shift map, and $
	\varphi_1$ is the substitution defined by $ 0\mapsto N(N+1)N;~ 1\mapsto (N+1)NN.
	$ 
\end{lemma}
\begin{proof}
	Set $\la_0=N$. Then $\la_0=N+\tau_0$. By (\ref{eq:relation-lambda-tau}) and the definition of $(\tau_i)$ it follows that 
	\begin{equation}\label{eq:jan11-1}
		\la_{3k}=\tau_{k}+N,\quad \la_{3k+1}=1-\tau_{k}+N,\quad \la_{3k+2}=N
	\end{equation}
	for all $k\ge 0$. Note that $(\tau_i)_{i=0}^\f$ is the fixed point of the substitution $\phi: 0\mapsto 01;~1\mapsto 10$ (cf.~\cite[Ch.~6]{Allouche_Shallit_2003}). That is 
	\[
	(\tau_i)_{i=0}^\f=\phi(\tau_0)\phi(\tau_1)\phi(\tau_2)\cdots.
	\]
	Therefore, comparing  (\ref{eq:jan11-1}) with the definition of $\varphi_1: 0\mapsto (0+N)(1+N)N;~1\mapsto (1+N)(0+N)N$ we obtain that 
	\[
	(\la_i)_{i=0}^\f=\varphi_1(\tau_0)\varphi_1(\tau_1)\varphi_1(\tau_2)\cdots,
	\]
	completing the proof.
\end{proof}

Next we show that $(\la_i)$ is the quasi-greedy expansion of $1$ in the base $\beta_c(M)$.

\begin{proposition}
	\label{prop:beta_c-3N+1}
	Let $M=3N+1$. The quasi-greedy $\beta_c(M)$-expansion of $1$  is given by 
	\[
	\de(\beta_c(M))=\la_1\la_2\la_3\ldots=(N+1)N(N+1)NN(N+1)NNN\ldots,
	\]
	where   $(\la_i)$ is defined   in (\ref{eq:relation-lambda-tau}). Furthermore, $\beta_c(M)$ is a transcendental number. 
\end{proposition}
Let $(\overline{\tau_i})_{i=0}^\f=(1-\tau_i)=10010110\ldots\in\set{0,1}^\N$ be the reflection  of $(\tau_i)$. Set   
\begin{equation}\label{eq:Sigma-1}
	{\Sigma}_1:=\set{(N+1)N(N+1), NNN,(N+1)NN,  NN(N+1)}.
\end{equation}
We define  a block map  $\re_1: {\Sigma}_1  \ra {\Sigma}_1$ by
\begin{equation}\label{eq:substitution-3N+1}
	\begin{split}
		(N+1)N(N+1) \quad&\longleftrightarrow\quad NNN,   \qquad 
		(N+1)NN \quad\longleftrightarrow\quad  NN(N+1).
	\end{split}
\end{equation}
Then for a sequence $\a=a_1a_2\ldots \in{\Sigma}_1^\N$ we set $\re_1(\a):=\re_1(a_1)\re_1(a_2)\cdots\in\Sigma_1^\N$. We emphasize that each digit $a_i$ is a block of length $3$ from ${\Sigma}_1$. Note that $(\la_i)\in\Sigma_1^\N$. We define its \emph{conjugate}   by 
\[
(\hat\la_i)_{i=1}^\f:=\re_1(\la_1\la_2\la_3)\re_1(\la_4\la_5\la_6)\ldots\re_1(\la_{3k+1}\la_{3k+2}\la_{3k+3})\ldots\in\Sigma_1^\N.
\]
Then $(\hat\la_i)$ begins with $NNN(N+1)NN(N+1)N(N+1)NNN$. Comparing the definitions of $(\overline{\tau_i})$ and $(\hat\la_i)$ we obtain that for any $k\in\N_0$,
\begin{equation}\label{eq:relation-lambda-tau-conjugacy}
	\hat\la_{3k+1}=\overline{\tau_{2k+1}}+N,\quad \hat\la_{3k+2}=N,\quad \hat\la_{3k+3}=\overline{\tau_{2k+2}}+N.
\end{equation}

\begin{lemma}
	\label{lem:inequality-lambda-i}
	For any $n\in\N_0$ and any $0\le i<3\cdot 2^n$, we have 
	\begin{align}
		\hat\la_1\ldots \hat\la_{3\cdot 2^n-i}&\prec \la_{i+1}\ldots \la_{3\cdot 2^n}\lle \la_1\ldots\la_{3\cdot 2^n-i},\label{eq:inequ-1}\\
		\hat\la_1\ldots\hat\la_{3\cdot 2^n-i}&\lle \hat\la_{i+1}\ldots\hat\la_{3\cdot 2^n}\prec \la_1\ldots \la_{3\cdot 2^n-i}. \label{eq:inequ-2}
	\end{align}
\end{lemma}
\begin{proof}
	Since the proof of (\ref{eq:inequ-2}) is similar to that of (\ref{eq:inequ-1}),    we only prove  (\ref{eq:inequ-1}). Take $n\in\N_0$ and  $i\in\set{0,1,\ldots, 3\cdot 2^n-1}$. We will prove (\ref{eq:inequ-1}) by    the following three cases: (I) $i=3k$, (II) $i=3k+1$, and (III) $i=3k+2$ for some $0\le k<2^n$. 
	
	Case (I). $i=3k$ with $0\le k<2^n$. Then by (\ref{eq:relation-lambda-tau}) we have 
	\begin{equation}\label{eq:28-1}
		\la_{i+1}\ldots \la_{3\cdot 2^n} =(\tau_{2k+1}+N)N(\tau_{2k+2}+N)\cdots(\tau_{2^{n+1}-1}+N)N(\tau_{2^{n+1}}+N).
	\end{equation}
	Note by \cite{Komornik-Loreti-1998} that  
	\begin{equation}\label{eq:28-2}
		\overline{\tau_1\ldots\tau_{2^{n+1}-j}}\prec \tau_{j+1}\ldots \tau_{2^{n+1}}\lle \tau_1\ldots \tau_{2^{n+1}-j}\quad\forall 0\le j<2^{n+1}.
	\end{equation}
	Then by (\ref{eq:28-1}), (\ref{eq:28-2}), (\ref{eq:relation-lambda-tau}) and (\ref{eq:relation-lambda-tau-conjugacy})  it follows that 
	\[
	\la_{i+1}\ldots \la_{3\cdot 2^n}\lle(\tau_1+N)N(\tau_2+N)\cdots (\tau_{2^{n+1}-2k-1}+N)N(\tau_{2^{n+1}-2k}+N)=\la_1\ldots \la_{3\cdot 2^n-i}
	\]
	and 
	\[
	\la_{i+1}\ldots \la_{3\cdot 2^n}\succ(\overline{\tau_1}+N)N(\overline{\tau_2}+N)\cdots (\overline{\tau_{2^{n+1}-2k-1}}+N)N(\overline{\tau_{2^{n+1}-2k}}+N)=\hat\la_1\ldots \hat\la_{3\cdot 2^n-i}
	\]
	as desired. 
	
	Case (II). $i=3k+1$ with $0\le k<2^n$. If $k=2^n-1$, then 
	\[\la_{i+1}\la_{i+2}=\la_{3\cdot 2^n-1}\la_{3\cdot 2^n}=N(\tau_{2^{n+1}}+N)=N(N+1),\]
	where the last equality follows by using $\tau_{2^m}=1$ for all $m\ge 0$. Note that $\la_1\la_2=(N+1)N$ and $\hat\la_1\hat\la_2=NN$. Then it is easy to verify (\ref{eq:inequ-1}) when $i=3\cdot (2^n-1)+1$.
	
	Next we assume $0\le k<2^n-1$. Then by (\ref{eq:relation-lambda-tau}) we have $\la_{i+1}=\la_{3k+2}=N<N+1=\la_1$, which yields $\la_{i+1}\ldots \la_{3\cdot 2^n}\prec \la_1\ldots \la_{3\cdot 2^n-i}$. On the other hand, observe by (\ref{eq:relation-lambda-tau}) that 
	\begin{align*}
		\la_{i+1}\la_{i+2}\la_{i+3}&=N(\tau_{2k+2}+N)(\tau_{2k+3}+N)\\
		&=N(\tau_{k+1}+N)(1-\tau_{k+1}+N)\\
		&\lge NN(N+1)\succ NNN=\hat\la_1\hat\la_2\hat\la_3,
	\end{align*}
	where the second equality holds by the definition  of $(\tau_i)$ that $\tau_{2i}=\tau_i$ and $\tau_{2i+1}=1-\tau_i$ for all $i\ge 0$. Thus, $\la_{i+1}\ldots \la_{3\cdot 2^n}\succ \hat\la_1\ldots\hat\la_{3\cdot 2^n-i}$ as required.

	Case (III). $i=3k+2$ with $0\le k<2^n$. If $k=2^n-1$, then $\la_{i+1}=\la_{3\cdot 2^n}=N+1$. Since $\la_1=N+1$ and $\hat\la_1=N$, it is clear that (\ref{eq:inequ-1}) holds for $i=3\cdot(2^n-1)+2$. In the following we assume $0\le k<2^n-1$. Then by (\ref{eq:relation-lambda-tau}) we obtain
	\begin{align*}
		\la_{i+1}\la_{i+2}\la_{i+3}&=(\tau_{2k+2}+N)(\tau_{2k+3}+N)N\\
		&=(\tau_{k+1}+N)(1-\tau_{k+1}+N)N\\
		&\lle(N+1)NN\prec (N+1)N(N+1)=\la_1\la_2\la_3,
	\end{align*}
	and
	\[
	\la_{i+1}\la_{i+2}\la_{i+3}=(\tau_{k+1}+N)(1-\tau_{k+1}+N)N\lge N(N+1)N\succ NNN=\hat\la_1\hat\la_2\hat\la_3.
	\]
	This proves (\ref{eq:inequ-1}), completing the proof. 
\end{proof}
To prove Proposition \ref{prop:beta_c-3N+1} we also need the following result due to Mahler \cite{Mahler_1976}.
\begin{lemma}\label{lem:mahler-transcendental}
	If $z$ is an algebraic number in the open unit disc, then 
	\[Z=\sum_{i=1}^{\f}\tau_i z^i\]
	is transcendental, where $(\tau_i)$ is the classical Thue-Morse sequence. 
\end{lemma}

\begin{proof}[Proof of Proposition \ref{prop:beta_c-3N+1}]
	Note by Definition \ref{def:beta-c}  and (\ref{eq:relation-lambda-tau}) that $(\la_i)$ is an expansion of $1$ in the base $\beta_c=\beta_c(M)$. Then by Lemmas \ref{lem:char-quasi-greedy-expansion} and \ref{lem:inequality-lambda-i} it follows that $\de(\beta_c)=(\la_i)$. So, by (\ref{eq:relation-lambda-tau}) we obtain 
	\begin{align*}
		1=\sum_{i=1}^{\f}\frac{\la_i}{\beta_c^i} & =\sum_{k=0}^\f\frac{\tau_{2k+1}+N}{\beta_c^{3k+1}}+\sum_{k=0}^{\f}\frac{N}{\beta_c^{3k+2}}+\sum_{k=0}^{\f}\frac{\tau_{2k+2}+N}{\beta_c^{3k+3}} \\
		& =\sum_{i=1}^{\f}\frac{N}{\beta_c^i}+\sum_{k=0}^{\f}\frac{1-\tau_k}{\beta_c^{3k+1}}+\sum_{k=0}^{\f}\frac{\tau_{k+1}}{\beta_c^{3k+3}}\\
		&=\frac{N}{\beta_c-1}+\frac{\beta_c^2}{\beta_c^3-1}-\frac{1}{\beta_c}\sum_{k=1}^{\f}\frac{\tau_k}{\beta_c^{3k}}+\sum_{k=1}^{\f}\frac{\tau_k}{\beta_c^{3k}},
	\end{align*}
	where in the last equality we  used $\tau_0=0$. 
	Rearranging the above equation yields that
	\begin{equation}\label{eq:jan15-1}
		\sum_{k=1}^{\f}\frac{\tau_k}{\beta_c^{3k}}=\frac{1-\frac{N}{\beta_c-1}-\frac{\beta_c^2}{\beta_c^3-1}}{1-\frac{1}{\beta_c}}.
	\end{equation}
	Suppose on the contrary that $\beta_c\in(N+1, N+2)$ is algebraic. Then so is $\beta_c^{-3}\in(0,1)$. By Lemma \ref{lem:mahler-transcendental} it follows that the left-hand side of (\ref{eq:jan15-1}) is transcendental, while the right-hand side of (\ref{eq:jan15-1}) is algebraic, leading to a contradiction. So, $\beta_c$ is transcendental. 
\end{proof}

\subsection{Patterns in $\widetilde{\U}_{\beta_c,M}$}
We will construct a sequence of bases $(\beta_n)_{n=0}^\f$ strictly increasing to $\beta_c=\beta_c(M)$ such that  the difference set $\widetilde{\U}_{\beta_{n+1}, M}\setminus\widetilde{\U}_{\beta_n, M}$ is countably infinite for all $n\ge 0$. In fact we will show that $\widetilde{\U}_{\beta_{n+1}, M}\setminus\widetilde{\U}_{\beta_n, M}$ consists of only eventually periodic sequences. 

First we construct a sequence of blocks $(\t_n)$ approximating $(\la_i)=\de(\beta_c)$ componentwisely.    Recall from (\ref{eq:Sigma-1}) that $\Sigma_1=\set{(N+1)N(N+1), NNN, (N+1)NN, NN(N+1)}$. 
Let  
\begin{equation}\label{eq:tn-3N+1}
	\t_0=(N+1)NN, \quad\textrm{and}\quad\t_{n+1}:=\t_n^+\re_1(\t_n^+)\quad\forall n\ge 0,
\end{equation}
where $\re_1$ is the block map defined in (\ref{eq:substitution-3N+1}), and for a word $\om=\om_1\ldots \om_m\in\set{N, N+1}^m$ we set $\om^+=\om_1\ldots\om_{m-1}(\om_m+1)$.  Then $\t_n\in\Sigma_1^{2^n}$ for any $n\ge 0$.
For example,
\[\t_1=(N+1)N(N+1)NNN, \quad\t_2= (N+1)N(N+1)NN(N+1)\,NNN(N+1)NN.\]
Observe   that each $\t_n^+$ is a prefix of $(\la_i)$, i.e.,  $\t_n^+=\la_1\ldots \la_{3\cdot 2^{n}}$. By Lemma \ref{lem:inequality-lambda-i} it follows that 
$
\si^i(\t_n^\f)\lle \t_n^\f$ for all $ i\ge 0.
$
So, by Lemma \ref{lem:char-quasi-greedy-expansion} there exists a unique $\beta_n=\beta_n(M)\in(1, M+1]$ such that 
\begin{equation}\label{eq:beta-n-3N+1}
	\de(\beta_n)=\t_n^\f.
\end{equation} Note by Proposition \ref{prop:beta-G} that 
$
\de(\beta_G(M))=((N+1)NN)^\f=\t_0^\f=\de(\beta_0(M)).
$
Then $\beta_G(M)=\beta_0(M)$ by Lemma \ref{lem:char-quasi-greedy-expansion}. 
Observe that $\de(\beta_n)=\t_n^\f\nearrow (\la_i)=\de(\beta_c)$ as $n\to\f$. By Lemma \ref{lem:char-quasi-greedy-expansion} it follows that  
\[
\beta_G=\beta_0<\beta_1<\cdots<\beta_n<\beta_{n+1}<\cdots,\quad\textrm{and}\quad \beta_n\nearrow\beta_c\quad\textrm{as }n\to\f.\]

Note that the blocks $\t_n, n\ge 0$ are similar to the one-dimensional balanced words (see, e.g., \cite{Komornik-Loreti-1998, Komornik-Loreti-2002}). In terms of Proposition \ref{prop:lexicographical-characterization-U}, we construct a sequence of two dimensional balanced words $(\a_n)_{n\ge 0}$  appearing in sequences of $\widetilde{\U}_{\beta_c, M}$. These words $(\a_n)$ are like basic bricks to describe the set $\widetilde{\U}_{\beta, M}$ for $\beta\le\beta_c(M)$ with $M=3N+1$. Let 
\begin{align*}
	\pazocal A_1 & :=\left\{\al_{(N+1)N}\al_{N(N+1)}\al_{(N+1)N}, \quad \al_{NN}\al_{N(N+1)}\al_{NN},\right. \\
	&\hspace{1.4cm}\left.\al_{(N+1)N}\al_{N(N+1)}\al_{NN}, \quad\al_{NN}\al_{N(N+1)}\al_{(N+1)N}\right\}.
\end{align*} 
We define the substitution $\Phi_1: \pazocal A_1\to\pazocal A_1$ by
\begin{equation}
	\label{eq:substitution-Phi-3N+1}
	\begin{split}
		\al_{(N+1)N}\al_{N(N+1)}\al_{(N+1)N}\quad&\longleftrightarrow \quad  \al_{NN}\al_{N(N+1)}\al_{NN}, \\
		\al_{(N+1)N}\al_{N(N+1)}\al_{NN}\quad & \longleftrightarrow\quad \al_{NN}\al_{N(N+1)}\al_{(N+1)N}.
	\end{split}
\end{equation}
By using $\Phi_1$ we construct a sequence of blocks $\set{\a_n}_{n=0}^\f$ as follows:
\begin{equation}\label{eq:an-3N+1}
	\a_0=\al_{(N+1)N}\al_{N(N+1)}\al_{NN},\quad\textrm{and}\quad \a_{n+1}=\a_n^+\Phi_1(\a_n^+)\quad\textrm{for }n\ge 0,
\end{equation}
where $\a_n^+$ is the word in $\pazocal A_1^*$ by changing the last digit of $\a_n$ from $\al_{NN}$ to $\al_{(N+1)N}$.
Then by the definition of $\Phi_1$ each $\a_n$ has length $3\cdot 2^n$ and ends with $\al_{NN}$.  If we write $\a_n=a_1\ldots a_{3\cdot 2^n}$, then $a_{3k+2}=a_{N(N+1)}$ for all $0\le k<2^n$. By (\ref{eq:an-3N+1}) it follows that 
\[
\a_0\sim\left(\begin{array}{ccc}
N+1&N&N \\
N&N+1&N \\
N&N&N+1
\end{array}\right),\quad \a_1\sim \left(\begin{array}{cc}
(N+1)N(N+1)&NNN \\
N(N+1)N&N(N+1)N \\
NNN&(N+1)N(N+1)
\end{array}\right),
\]
and
\[
\a_2\sim\left(\begin{array}{cccc}
(N+1) N (N+1)& N N (N+1) & N N N &(N+1) N N\\
N(N+1) N &N (N+1) N & N (N+1) N &N (N+1) N \\
N N N &(N+1) N N& (N+1) N (N+1)& N N (N+1)
\end{array}\right).
\]
By comparing the definitions of $(\t_n)$ and $(\a_n)$ in  (\ref{eq:tn-3N+1}) and (\ref{eq:an-3N+1}) respectively,  it follows that 
\begin{equation}\label{eq:relation-an-tn}
\a_n^1=\t_n,\quad \a_n^2=(N(N+1)N)^n\quad\textrm{and}\quad \overline{\a_n^\+}=\re_1(\t_n)
\end{equation} 
for any $n\ge 0$, where for a word $\d=d_1\ldots d_k\in\Om_{M}^{*}$ we write $\d^1=d_1^1\ldots d_k^1$, $\d^2=d_1^2\ldots d_k^2$ and $\overline{\d^\+}=\overline{d_1^\+}\ldots \overline{d_k^\+}$ with each $d_i=(d_i^1, d_i^2)$ and $\overline{d_i^\+}=M-d_i^\+=M-(d_i^1+d_i^2)$.

The big difference between $M=1$ and $M=3N+1$ with $N\ge 1$ is that in the latter case the digit $N-1$ might occur in   $d_n^1=N-1, d_n^2=N-1$ or $\overline{d_n^\+}=N-1$ for some sequence $(d_i)\in\widetilde{\U}_{\beta_c, M}$. In other words, when $M=3N+1$ with $N\ge 1$, sequences in $\widetilde{\U}_{\beta_c, M}$ might contain digits from 
\begin{equation}\label{eq:exception-set}
\mathcal E_1:=\bigcup_{j=0}^2\set{\theta^j(\al_{(N+1)(N+1)})}=\set{\al_{(N+1)(N+1)}, \al_{(N-1)(N+1)}, \al_{(N+1)(N-1)}}.
\end{equation}
Fortunately,   we can prove   that  any sequence of  $\widetilde{\U}_{\beta_c, M}$ contains at most three digits  from $\mathcal E_1$. 

\begin{lemma}
\label{lem:exceptional-digits}
Let  $M=3N+1$ with $N\ge 1$. If $(d_i)\in\widetilde{\U}_{\beta_c, M}$, then 
\[
\#\set{i\in\N: d_i\in\mathcal E_1}\le 3.
\]
\end{lemma}
Before the proof we give an example to illustrate the main idea.
\begin{example}
Let $M=4$. Take $(d_i)\in\widetilde{\U}_{\beta_c, 4}\setminus\set{\al_{40}^\f, \al_{04}^\f, \al_{00}^\f}$. Note that $\de(\beta_c)=212112111\ldots$. By Proposition \ref{prop:lexicographical-characterization-U} there exists $n_0\in\N$ such  that for all $i\ge n_0$ we have either
\[
d_i\in\set{\al_{21}, \al_{12}, \al_{11}}\sim\set{\left(\begin{array}{c}
		2 \\
		1\\
		1 
	\end{array}\right), \left(\begin{array}{c}
		1 \\
		2\\
		1 
	\end{array}\right), \left(\begin{array}{c}
		1 \\
		1\\
		2 
	\end{array}\right)}
\]
or 
\[
d_i\in\mathcal E_1=\set{\al_{22}, \al_{02}, \al_{20}}\sim \set{\left(\begin{array}{c}
		2 \\
		2\\
		0 
	\end{array}\right), \left(\begin{array}{c}
		0 \\
		2\\
		2 
	\end{array}\right), \left(\begin{array}{c}
		2 \\
		0\\
		2 
	\end{array}\right)}.
\]
Suppose $d_n\in\mathcal E_1$ for some large integer $n>n_0$, say $d_n=\al_{22}$. Then by Proposition \ref{prop:lexicographical-characterization-U} we must have $d_{n-1}=d_{n+1}=\al_{11}$. If $d_{n+2}=\al_{22}\in\mathcal E_1$, then we can deduce by using Proposition \ref{prop:lexicographical-characterization-U} that 
\[
d_{n-1}d_n\ldots d_{n+4}\sim\left(\begin{array}{c}
	121211 \\
	121211 \\
	202022 
\end{array}\right).
\] This implies  $\overline{d_{n+3}^\+d_{n+4}^\+}=22$, which is forbidden in  sequences of $\widetilde{\U}_{\beta_c,4}$. So we can only have $d_{n+2}\in\set{\al_{21}, \al_{12}}$. Similar ideas can be applied to deduce that any digit from $\mathcal E_1$ can not occur in the tail sequence $d_{n+1}d_{n+2}\ldots$. Furthermore, we can apply this idea to prove that there are at most three digits from $\mathcal E_1$  appearing  in the prefix $d_1\ldots d_n$.  
\end{example}
\begin{proof}[Proof of Lemma \ref{lem:exceptional-digits}]
Take $(d_i)\in\widetilde{\U}_{\beta_c, M}$. Note by Proposition \ref{prop:beta-G}   that 
$\widetilde{\U}_{\beta_c, M}\supset\widetilde{\U}_{\beta_0, M}=\set{\al_{00}^\f, \al_{0M}^\f, \al_{M0}^\f}.$
If $(d_i)\in\widetilde{\U}_{\beta_0, M}$, then we are done. So in the following we assume  $(d_i)\in\widetilde{\U}_{\beta_c, M}\setminus\widetilde{\U}_{\beta_0, M}$. Note by Proposition \ref{prop:beta_c-3N+1} that $\de(\beta_c)\prec (N+1)^\f$. By Proposition \ref{prop:lexicographical-characterization-U} and (\ref{eq:sum-equal}) there exists a smallest $m\in\N$ such that 
\begin{equation}\label{eq:20-1}
	d_{m+k}^1, ~d_{m+k}^2, ~\overline{d_{m+k}^\+}\le N+1<M\quad\forall k\ge 0.
\end{equation} 

Suppose on the contrary there exist $0\le p_1<p_2<p_3<p_4$ such that 
\begin{equation}\label{eq:20-2}
	d_{m+p_j}\in\set{\al_{(N+1)(N+1)}, \al_{(N-1)(N+1)}, \al_{(N+1)(N-1)}}\quad\textrm{for all }1\le j\le  4.
\end{equation}  Write $p:=p_4$, and without loss of generality we assume $d_{m+p}=\al_{(N+1)(N+1)}$. 
Then  by (\ref{eq:20-1}),  Propositions \ref{prop:lexicographical-characterization-U} and \ref{prop:beta_c-3N+1} it follows that 
\[d_{m+p-1}d_{m+p}d_{m+p+1}=\al_{NN}\al_{(N+1)(N+1)}\al_{NN}.\] For $d_{m+p+2}$ we have three possibilities: $d_{m+p+2}\in\set{\al_{(N+1)N}, \al_{N(N+1)}, \al_{(N+1)(N+1)}}$.

{\bf Claim 1.} $d_{m+p+2}\ne \al_{(N+1)(N+1)}$.

Suppose on the contrary that $d_{m+p+2}=\al_{(N+1)(N+1)}$. Then by using $d_{m+p}d_{m+p+1}d_{m+p+2}=\al_{(N+1)(N+1)}\al_{NN}\al_{(N+1)(N+1)}$,   and Propositions \ref{prop:lexicographical-characterization-U} and \ref{prop:beta_c-3N+1} we must have \[d_{m+p+3}d_{m+p+4}=\al_{NN}\al_{NN}.\]
This implies that $\overline{d_{m+p+2}^\+}=N-1<M$ and $\overline{d_{m+p+3}^\+d_{m+p+4}^\+}=(N+1)(N+1)\succ \de_1(\beta_c)\de_2(\beta_c)$, leading to a contradiction to $(d_i)\in\widetilde{\U}_{\beta_c, M}$. So, Claim 1 holds.

By Claim 1 we have $d_{m+p+2}=\al_{(N+1)N}$ or $\al_{N(N+1)}$. By symmetry we may assume  $d_{m+p+2}=\al_{(N+1)N}$. Then by (\ref{eq:20-1}),  and Propositions \ref{prop:lexicographical-characterization-U} and \ref{prop:beta_c-3N+1} it follows that 
\[d_{m+p-2}=\al_{N(N+1)}\quad \textrm{and}\quad d_{m+p-3}\in\set{\al_{(N+1)N}, \al_{NN}, \al_{(N+1)(N-1)}}.\]

{\bf Claim 2.} $d_{m+p-3}\ne \al_{(N+1)(N-1)}$.

Suppose on the contrary that  $d_{m+p-3}=\al_{(N+1)(N-1)}$. Then 
\[
d_{m+p-3}\ldots d_{m+p+2}\quad\sim\quad\left(\begin{array}{cccccc}
	N+1&N&N&N+1&N&N+1 \\
	N-1&N+1&N&N+1&N&N \\
	N+1&N&N+1&N-1&N+1&N
\end{array}\right). 
\] 
By the definition of $p_i$ in (\ref{eq:20-2}) it follows that $p=p_4\ge 5$.
So, by (\ref{eq:20-1}), and Propositions \ref{prop:lexicographical-characterization-U} and \ref{prop:beta_c-3N+1}  it follows that $d_{m+p-4}=\al_{N(N+1)}$ and $d_{m+p-5}=\al_{(N+1)N}$ or $\al_{NN}$. However, if $d_{m+p-5}=\al_{(N+1)N}$, then by Proposition \ref{prop:beta_c-3N+1} we have
\[d_{m+p-5}^1\ldots d_{m+p+2}^1=(N+1)N(N+1)NN(N+1)N(N+1)\succ \de_1(\beta_c)\ldots \de_8(\beta_c),\]
leading to a contradiction. If $d_{m+p-5}=\al_{NN}$, then $\overline{d_{m+p-5}^\+}\ldots \overline{d_{m+p+1}^\+}=(N+1)N(N+1)N(N+1)\succ \de_1(\beta_c)\ldots \de_5(\beta_c)$, again leading to a contradiction. This proves Claim 2. 

So, by Claim 2 we have $d_{m+p-3}\in\set{\al_{(N+1)N}, \al_{NN}}$. Without loss of generality we may assume $d_{m+p-3}=\al_{NN}$. Then 
\[
d_{m+p-3}\ldots d_{m+p+2}\quad\sim\quad\left(\begin{array}{cccccc}
	N&N&N&N+1&N&N+1 \\
	N&N+1&N&N+1&N&N \\
	N+1&N&N+1&N-1&N+1&N
\end{array}\right),
\]
which implies $p=p_4\ge 6$ by (\ref{eq:20-2}). 
By   (\ref{eq:20-1}),  and Propositions \ref{prop:lexicographical-characterization-U} and \ref{prop:beta_c-3N+1} we obtain that 
$d_{m+p-5}d_{m+p-4}=\al_{N(N+1)}\al_{(N+1)N}$ and $d_{m+p-6}\in\set{\al_{(N+1)N}, \al_{NN}, \al_{(N+1)(N-1)}}$. Similar to Claim 2 we can prove that $d_{m+p-6}\ne \al_{(N+1)(N-1)}$.

Proceeding this argument we can prove that 
\[
d_{m+p-3i}d_{m+p-3i+1}d_{m+p-3i+2}\in\set{\al_{(N+1)N}\al_{N(N+1)}\al_{NN}, \al_{NN}\al_{N(N+1)}\al_{NN}} 
\]
for all $1\le i\le \lfloor\frac{p}{3}\rfloor$.
This leads to a contradiction to (\ref{eq:20-2}), completing the proof.  
\end{proof}

\begin{lemma}
\label{lem:admissible-alphabet-3N+1}
Let $M=3N+1$ with $N\in\N_0$. Then  for any sequence $(d_i)\in\widetilde{\U}_{\beta_c, M}\setminus\widetilde{\U}_{\beta_0, M}$ there exists $m_1\in\N$ such that 
\[
d_{k}\in\set{\al_{(N+1)N}, \al_{NN}, \al_{N(N+1)}}\quad\forall k>m_1.
\]
\end{lemma} 
\begin{proof}
Take $(d_i)\in\widetilde{\U}_{\beta_c, M}\setminus\widetilde{\U}_{\beta_0, M}$. If $M=3N+1$ with $N=0$, then it is clear that $d_i\in\set{\al_{(N+1)N}, \al_{NN}, \al_{N(N+1)}}$ for all $i\ge 1$. In the following we assume $N\ge 1$. Note by Proposition \ref{prop:beta-G} that $\widetilde{\U}_{\beta_0, M}=\widetilde{\U}_{\beta_G, M}=\set{\al_{0M}^\f, \al_{M0}^\f, \al_{00}^\f}$. Furthermore,  by Proposition \ref{prop:beta_c-3N+1} and Lemma \ref{lem:char-quasi-greedy-expansion} that $\de(\beta_c)\prec (N+1)^\f$. Then there exists $m_0\in\N$ such that $d_{m_0}^1<M, d_{m_0}^2<M$ and $\overline{d_{m_0}^\+}<M$. Since $M=3N+1$ with $N\ge 1$,  by (\ref{eq:sum-equal}) and Proposition \ref{prop:lexicographical-characterization-U} it follows that 
\[
d_{k}\in\bigcup_{j=0}^2\set{\theta^j(\al_{(N+1)N}), \theta^j(\al_{(N+1)(N+1)})}\quad \forall k>m_0.
\]
Hence, by Lemma \ref{lem:exceptional-digits} there exists $m_1>m_0$ such that 
\[
d_k\in\bigcup_{j=0}^2\set{\theta^j(\al_{(N+1)N})}=\set{\al_{(N+1)N}, \al_{N(N+1)}, \al_{NN}}\quad\forall k>m_1,
\]
completing the proof.
\end{proof}

Now we are ready to describe the patterns appearing in sequences of $\widetilde{\U}_{\beta_c, M}$. 
\begin{proposition}
\label{prop:patterns-3N+1}
Let $(d_i)\in\widetilde{\U}_{\beta_c,M}$ and $n\in\N_0$. Then there exists $m_1\in\N$ such that for all $m>m_1$ the following statements hold.
\begin{enumerate}[{\rm(i)}]
	\item If $d^1_m<M$ and $d_{m+1}\ldots d_{m+3\cdot 2^n}=\a_n^+$, then 
	$
	d_{m+3\cdot 2^n+1}\ldots d_{m+3\cdot 2^{n+1}}\in\set{\Phi_1(\a_n), \Phi_1(\a_n^+)}.
	$
	\item If $d^1_m>0$ and $d_{m+1}\ldots d_{m+3\cdot 2^n}=\Phi_1(\a_n^+)$, then 
	$
	d_{m+3\cdot 2^n+1}\ldots d_{m+3\cdot 2^{n+1}}\in\set{\a_n, \a_n^+}.
	$
\end{enumerate}
\end{proposition}
\begin{remark}\mbox{}
\begin{itemize}  
	\item Recall $(\a_n)$ from (\ref{eq:an-3N+1}). Then Propositin \ref{prop:patterns-3N+1} implies that if $d_m^1<M$ and $d_{m+1}\ldots d_{m+3\cdot 2^n}=\a_n^+$, then $d_{m+1}\ldots d_{m+3\cdot 2^{n+1}}=\a_{n+1}$ or $\a_{n+1}^+$. Similarly, if  $d_m^1>0$ and $d_{m+1}\ldots d_{m+3\cdot 2^n}=\Phi_1(\a_n^+)$, then $d_{m+1}\ldots d_{m+3\cdot 2^{n+1}}=\Phi_1(\a_{n+1})$ or $\Phi_1(\a_{n+1}^+)$.
	\item  Observe by (\ref{eq:cyclic-permutation}) that $(d_i)\in\widetilde{\U}_{\beta, M}$ if and only if $\theta((d_i))\in\widetilde{\U}_{\beta, M}$, which is also equivalent to $\theta^2((d_i))\in\widetilde{\U}_{\beta, M}$. So, Proposition \ref{prop:patterns-3N+1} also holds if we replace $(d_i)$ by $\theta((d_i))$ or $\theta^2((d_i))$. For example, if $\overline{d_m^\+}<M$ and $d_{m+1}\ldots d_{m+3\cdot 2^n}=\theta(\a_n^+)$, then $d_{m+3\cdot 2^n+1}\ldots d_{m+3\cdot 2^{n+1}}=\theta(\Phi_1(\a_n))$ or $\theta(\Phi_1(\a_n^+))$. 
\end{itemize}

\end{remark}

First we prove Proposition \ref{prop:patterns-3N+1} for $n=0$ and $n=1$.
\begin{lemma}
\label{lem:pattern-3N+1-0}
Let $(d_i)\in\widetilde{\U}_{\beta_c, M}$. Then there exists $m_1\in\N$ such that for all $m>m_1$ the following statements hold.
\begin{enumerate}[{\rm(i)}]
	\item If $d^1_m<M$ and $d_{m+1}d_{m+2}d_{m+3}=\a_0^+$, then 
	$
	d_{m+4}d_{m+5}d_{m+6}=\Phi_1(\a_0)$ or $\Phi_1(\a_0^+).
	$
	\item If $d^1_m>0$ and $d_{m+1}d_{m+2}d_{m+3}=\Phi_1(\a_0^+)$, then 
	$
	d_{m+4}d_{m+5} d_{m+6}=\a_0$ or $\a_0^+.
	$
\end{enumerate}
\end{lemma}
\begin{proof}Since the proof of (ii) is similar, we only prove (i). Take $(d_i)\in\widetilde{\U}_{\beta_c, M}$. If $(d_i)\in\widetilde{\U}_{\beta_0, M}$, then by Proposition \ref{prop:beta-G}  there is nothing to prove. In the following we assume $(d_i)\in\widetilde{\U}_{\beta_c, M}\setminus\widetilde{\U}_{\beta_0, M}$. Then by Lemma \ref{lem:admissible-alphabet-3N+1} there exists $m_1\in\N$ such that 
\begin{equation}\label{eq:jan15-2}
	d_j\in\set{\al_{(N+1)N}, \al_{N(N+1)}, \al_{NN}}\quad \forall j>m_1.
\end{equation}  Take $m>m_1$, and suppose
$d_{m+1}d_{m+2}d_{m+3}=\a_0^+=\al_{(N+1)N}\al_{N(N+1)}\al_{(N+1)N}$ with $d_{m}^1<M$. Then it suffices to prove 
\begin{equation}\label{eq:16-1}
	d_{m+4}d_{m+5}d_{m+6}\in\set{\al_{NN}\al_{N(N+1)}\al_{NN},\quad \al_{NN}\al_{N(N+1)}\al_{(N+1)N}}.
\end{equation}

We first prove $d_{m+4}=\al_{NN}$. Suppose on the contrary that  $d_{m+4}\ne\al_{NN}$.   Note that $d_{m+1}d_{m+2}d_{m+3}=\a_0^+$ and $d_{m}^1<M$. Moreover, observe by Proposition \ref{prop:beta_c-3N+1} that $\de(\beta_c)$ begins with $(N+1)N(N+1)NN(N+1)NNN$. Then by (\ref{eq:sum-equal}), (\ref{eq:jan15-2})  and Proposition \ref{prop:lexicographical-characterization-U} it follows that 
\[
d_{m+1}\ldots d_{m+9}\sim\left(\begin{array}{ccc}
	(N+1)N(N+1)&NN(N+1)&NNN \\
	N(N+1)N&(N+1)NN & N(N+1)(N+1)\\
	NNN&N(N+1)N&(N+1)NN
\end{array}\right)
\]
or 
\[
d_{m+1}\ldots d_{m+9}\sim\left(\begin{array}{ccc}
	(N+1)N(N+1)&NN(N+1)&NNN \\
	N(N+1)N&(N+1)NN & (N+1)NN\\
	NNN&N(N+1)N& N(N+1)(N+1)
\end{array}\right). 
\]
However, in the first case we have $d_7^2=N<M$ and $d_8^2d_9^2=(N+1)(N+1)\succ\de_1(\beta_c)\de_2(\beta_c)$; and in the second case we   have $\overline{d_7^\+}=N<M$ and $\overline{d_8^\+d_9^\+}=(N+1)(N+1)\succ\de_1(\beta_c)\de_2(\beta_c)$. Both cases lead to a contradiction to $(d_i)\in\widetilde{\U}_{\beta_c, M}$. So, $d_{m+4}=\al_{NN}$.

Note that $d_{m+1}\ldots d_{m+4}=\al_{(N+1)N}\al_{N(N+1)}\al_{(N+1)N}\al_{NN}$. Then by  (\ref{eq:sum-equal}), (\ref{eq:an-3N+1}) and Proposition  \ref{prop:lexicographical-characterization-U} it follows that $d_{m+5}d_{m+6}\in\set{\al_{N(N+1)}\al_{NN}, \al_{N(N+1)}\al_{(N+1)N}}$,  proving (\ref{eq:16-1}). 
\end{proof}

\begin{lemma}
\label{lem:pattern-3N+1-1}
Let $(d_i)\in\widetilde{\U}_{\beta_c, M}$. Then there exists $m_1\in\N$ such that for all $m>m_1$ the following statements hold.
\begin{enumerate}[{\rm(i)}]
	\item If $d^1_m<M$ and $d_{m+1}\ldots d_{m+6}=\a_1^+$, then 
	$
	d_{m+7}\ldots d_{m+12}=\Phi_1(\a_1)$ or  $\Phi_1(\a_1^+).
	$
	\item If $d^1_m>0$ and $d_{m+1}\ldots d_{m+6}=\Phi_1(\a_1^+)$, then 
	$
	d_{m+7}\ldots d_{m+12}=\a_1$ or $\a_1^+.
	$
\end{enumerate}
\end{lemma}

\begin{proof}
Since the proof of (ii) is similar, we only prove (i). Take $(d_i)\in\widetilde{\U}_{\beta_c, M}$. As in the proof of Lemma \ref{lem:pattern-3N+1-0} we may assume $(d_i)\in\widetilde{\U}_{\beta_c, M}\setminus\widetilde{\U}_{\beta_0, M}$. Then by Lemma \ref{lem:admissible-alphabet-3N+1} there exists $m_1\in\N$ such that $ d_j\in\set{\al_{(N+1)N}, \al_{N(N+1)}, \al_{NN}}$ for all $j>m_1$. Take $m>m_1$, and   suppose $d_{m+1}\ldots d_{m+6}=\a_1^+$ with $d_m^1<M$. Note by (\ref{eq:relation-an-tn}) that $d_{m+1}^1\ldots d_{m+6}^1=\t_1^+=(N+1)N(N+1)NN(N+1)$. Then by   Propositions \ref{prop:lexicographical-characterization-U} and \ref{prop:beta_c-3N+1} it follows that 
$
d_{m+7} \in \set{\al_{NN}, \al_{N(N+1)}}.
$

Suppose    $d_{m+7}=\al_{N(N+1)}$. Then by using $d_{m+1}\ldots d_{m+6}=\a_1^+$ and Proposition \ref{prop:beta_c-3N+1} it follows that 
\[
d_{m+1}\ldots d_{m+9}\quad\sim\quad\left(\begin{array}{ccc}
	(N+1)N(N+1)&NN(N+1)&N NN\\
	N(N+1)N& N(N+1)N& (N+1)NN \\
	NNN&NNN&N(N+1)(N+1)
\end{array}\right).
\]
This gives $\overline{d_{m+7}^\+}=N<M$ and $\overline{d_{m+8}^\+d_{m+9}^\+}=(N+1)(N+1)$, leading to a contradiction to $(d_i)\in\widetilde{\U}_{\beta_c, M}$. So, $d_{m+7}=\al_{NN}$.

Therefore, by Propositions \ref{prop:lexicographical-characterization-U} and \ref{prop:beta_c-3N+1} we can deduce that $d_{m+8}d_{m+9}=\al_{N(N+1)}\al_{NN}$, and thus
$d_{m+7}d_{m+8}d_{m+9}=\al_{NN}\al_{N(N+1)}\al_{NN}=\Phi(\a_0^+).$ Since $d_{m+6}^1=N+1>0$, by Lemma \ref{lem:pattern-3N+1-0} (ii) it follows that 
\[
d_{m+7}\ldots d_{m+12}=\Phi_1(\a_0^+)\a_0=\Phi_1(\a_1^+)\quad\textrm{or}\quad \Phi_1(\a_0^+)\a_0^+=\Phi_1(\a_1),
\]
completing the proof.           
\end{proof}

\begin{proof}
[Proof of Proposition \ref{prop:patterns-3N+1}]
Since the proof of (ii) is similar, we only prove (i). This will be done by induction on $n$.  Note by Lemmas \ref{lem:pattern-3N+1-0} and \ref{lem:pattern-3N+1-1} that (i) holds for $n=0, 1$. Now we assume (i) holds for all $k\le n$ with $n\in\N$, and we will prove it for $k=n+1$. Take $(d_i)\in\widetilde{\U}_{\beta_c, M}$. As in the proof of Lemma \ref{lem:pattern-3N+1-0} we may assume $(d_i)\in\widetilde{\U}_{\beta_c, M}\setminus\widetilde{\U}_{\beta_0, M}$. Then by Lemma \ref{lem:admissible-alphabet-3N+1} there exists $m_1\in\N$ such that $ d_j\in\set{\al_{(N+1)N}, \al_{N(N+1)}, \al_{NN}}$ for all $j>m_1$. Now take $m>m_1$, and suppose $d_{m+1}\ldots d_{m+3\cdot 2^{n+1}}=\a_{n+1}^+$ with $d_m^1<M$. We will prove 
\begin{equation}
	\label{eq:induction-3N+1}
	d_{m+3\cdot 2^{n+1}+1}\ldots d_{m+3\cdot 2^{n+2}}=\Phi_1(\a_{n+1}^+)\quad\textrm{or}\quad \Phi_1(\a_{n+1}).
\end{equation}

Note by (\ref{eq:relation-an-tn}) and (\ref{eq:tn-3N+1})  that 
$
d_{m+1}^1\ldots d_{m+3\cdot 2^{n+1}}^1=\t_{n+1}^+=\t_n^+\re_1(\t_n).
$
Furthermore,    $\de(\beta_c)$ begins with $\t_{n+1}^+\re_1(\t_{n+1})$ by Proposition \ref{prop:beta_c-3N+1}. Then by   (\ref{eq:tn-3N+1}) and Proposition \ref{prop:lexicographical-characterization-U} we obtain that 
\begin{equation}\label{eq:17-2}
	d_{m+3\cdot 2^{n+1}+1}^1\ldots d_{m+3\cdot 2^{n+2}}^1\lle \re_1(\t_{n+1})=\re_1(\t_n^+)\t_n^+=\re_1(\t_{n-1}^+)\t_{n-1}\t_n^+.
\end{equation}
Observe by (\ref{eq:an-3N+1}) that 
$
d_{m+1}\ldots d_{m+3\cdot 2^{n+1}}=\a_{n+1}^+=\a_n^+\Phi_1(\a_n)=\a_n^+\Phi_1(\a_{n-1}^+)\a_{n-1}^+,
$
which gives that $d^1_{m+3(2^n+2^{n-1})}=N<M$ and $d_{m+3(2^n+2^{n-1})+1}\ldots d_{m+3\cdot 2^{n+1}}=\a_{n-1}^+$. Then by the induction hypothesis for $k=n-1$ we have 
\[ 
d_{m+3\cdot 2^{n+1}+1}\ldots d_{m+3(2^{n+1}+2^{n-1})}\in\set{\Phi_1(\a_{n-1}^+), \Phi_1(\a_{n-1})}.
\]
Note that $\Phi_1(\a_{n-1})^1=\re_1(\t_{n-1})\succ \re_1(\t_{n-1}^+)$, and  $d_{m+3(2^{n}+2^{n-1})}^1=N<M$. Then by (\ref{eq:17-2}) and Proposition \ref{prop:lexicographical-characterization-U} it follows that 
\begin{equation}\label{eq:17-3}
	d_{m+3\cdot 2^{n+1}+1}\ldots d_{m+3(2^{n+1}+2^{n-1})}=\Phi_1(\a_{n-1}^+).
\end{equation}
Note that $d_{m+3\cdot 2^{n+1}}^1=N+1>0$.   Again by the induction hypothesis we obtain that 
$d_{m+3( 2^{n+1}+2^{n-1})+1}\ldots d_{m+3(2^{n+1}+2^n)}\in\set{\a_{n-1}, \a_{n-1}^+}.$
Then by (\ref{eq:17-2}) and Proposition \ref{prop:lexicographical-characterization-U} it follows that 
\begin{equation}\label{eq:17-4}
	d_{m+3(2^{n+1}+2^{n-1})+1}\ldots d_{m+3(2^{n+1}+2^n)}=\a_{n-1}.
\end{equation}
So, by (\ref{eq:17-3}) and (\ref{eq:17-4}) we have
\[d_{m+3\cdot 2^{n+1}+1}\ldots d_{m+3(2^{n+1}+2^n)}=\Phi_1(\a_{n-1}^+)\a_{n-1}=\Phi_1(\a_n^+).\] Since $d_{m+3\cdot 2^{n+1}}=N+1>0$, by the induction hypothesis for $k=n$ we obtain that 
$d_{m+3(2^{n+1}+2^n)+1}\ldots d_{m+3\cdot 2^{n+2}}=\a_n$ or $\a_n^+$. Therefore, by (\ref{eq:an-3N+1}) we conclude that 
\[
d_{m+3\cdot 2^{n+1}+1}\ldots d_{m+3\cdot 2^{n+2}}=\Phi_1(\a_n^+)\a_n^+=\Phi_1(\a_{n+1})\quad\textrm{or}\quad \Phi_1(\a_n^+)\a_n=\Phi_1(\a_{n+1}^+),
\]
establishing  (\ref{eq:induction-3N+1}). Hence, by induction this complets the proof. 
\end{proof}

\subsection{Phase transition for $\widetilde{\U}_{\beta, M}$ with $M=3N+1$}
Recall that the sequence $(\beta_n)$ is defined in (\ref{eq:beta-n-3N+1}).
Note that $\beta_n\nearrow\beta_c$ as $n\to\f$. Our strategy on proving  Theorem \ref{th:main-result} (ii)   is to show that  the difference set $\widetilde{\U}_{\beta_{n+1}, M}\setminus\widetilde{\U}_{\beta_n, M}$ is countable for any $n\ge 0$.
\begin{proposition}
\label{th:difference-Un+1-Un}
Let  $n\in\N_0$. Then any sequence in $\widetilde{\U}_{\beta_{n+1}, M}\setminus\widetilde{\U}_{\beta_n, M}$ must end in
\[
\bigcup_{j=0}^2\set{  \theta^j(\a_n^\f), \theta^j(\Phi_1(\a_n)^\f)},
\]
where $\theta$ is the cyclic permutation defined in (\ref{eq:theta}).
\end{proposition}

First we prove Proposition  \ref{th:difference-Un+1-Un} for $n=0$.
\begin{lemma}
\label{lem:difference-U1-U0}
Any sequence in $\widetilde{\U}_{\beta_1, M}\setminus\widetilde{\U}_{\beta_0, M}$ must end in $ \a_0^\f, \Phi_1(\a_0)^\f$ or their cyclic permutations under $\theta$.
\end{lemma}
\begin{proof}
Take $(d_i)\in\widetilde{\U}_{\beta_1, M}\setminus\widetilde{\U}_{\beta_0, M}$.
By Lemma \ref{lem:admissible-alphabet-3N+1} there exists $m\in\N$ such that $d_m^1<M$ and 
$
d_{m+i}\in \set{\al_{(N+1)N}, \al_{N(N+1)}, \al_{NN}}$ for all $ i\ge 1.$ 
Without loss of generality we assume $d_{m+1}=\al_{(N+1)N}$. 
Note that 
$
\de(\beta_1)=\t_1^\f=((N+1)N(N+1)NNN)^\f.
$
Then by   Proposition \ref{prop:lexicographical-characterization-U} it follows that 
\begin{equation}\label{eq:30-2}
	d_{m+1}d_{m+2}=\al_{(N+1)N}\al_{N(N+1)}\quad\textrm{and}\quad d_{m+3}\in\set{\al_{(N+1)N}, \al_{NN}},  
\end{equation}
or
\begin{equation}\label{eq:jan15-3}
	d_{m+1}d_{m+2}=\al_{(N+1)N}\al_{NN}\quad\textrm{and}\quad d_{m+3}\in\set{\al_{(N+1)N}, \al_{N(N+1)}}.
\end{equation}

Suppose (\ref{eq:30-2}) holds, and we assume  $d_{m+3}=\al_{(N+1)N}$.  Then $d_{m+1}d_{m+2}d_{m+3}=\a_0^+$. By Proposition \ref{prop:patterns-3N+1} we have $d_{m+4}d_{m+5}d_{m+6}=\Phi_1(\a_0^+)$ or $\Phi_1(\a_0)$. However, if $d_{m+4}d_{m+5}d_{m+6}=\Phi_1(\a_0)$, then we have
\[
d_m^1<M\quad\textrm{and}\quad d_{m+1}^1\ldots d_{m+6}^1=\t_1^+\succ\de_1(\beta_1)\ldots \de_6(\beta_1),
\]
leading to a contradiction to $(d_i)\in\widetilde{\U}_{\beta_1, M}$. So, $d_{m+4}d_{m+5}d_{m+6}=\Phi_1(\a_0^+)$. Again, by Proposition \ref{prop:patterns-3N+1} and the same argument as above we can deduce that 
\[d_{m+1}d_{m+2}\ldots =(\a_0^+\Phi_1(\a_0^+))^\f=\a_1^\f.\]
However, by (\ref{eq:relation-an-tn}) this implies that  $d^1_{m+1}d^1_{m+2}\ldots=\t_1^\f=\de(\beta_1)$, which again leads to a contradiction to $(d_i)\in\widetilde{\U}_{\beta_1, M}$. 

In view of (\ref{eq:30-2}), this implies    $d_{m+1}d_{m+2}d_{m+3}=\al_{(N+1)N}\al_{N(N+1)}\al_{NN}=\a_0$. By   Proposition \ref{prop:lexicographical-characterization-U} we have that $d_{m+4}=\al_{(N+1)N}$ or $\al_{N(N+1)}$. If $d_{m+4}=\al_{N(N+1)}$, then $d_{m+2}d_{m+3}d_{m+4}=\theta(\a_0^+)$. By the same argument as in the claim we can prove that $d_{m+2}d_{m+3}\ldots=\theta(\a_1^\f)$, which will lead to a contradiction. 
So, $d_{m+4}=\al_{(N+1)N}$. Then by   Proposition \ref{prop:lexicographical-characterization-U} it follows that $d_{m+5}=\al_{NN}$ or $\al_{N(N+1)}$. If $d_{m+5}=\al_{NN}$, then $d_{m+3}d_{m+4}d_{m+5}=\theta^2(\a_0^+)$, again the same argument  as in the claim yields  that $d_{m+3}d_{m+4}\ldots=\theta^2(\a_1^\f)$, which will lead to a contradiction.  
Therefore, $d_{m+5}=\al_{N(N+1)}$.  
Repeating the above argument  we conclude that  $d_{m+1}d_{m+2}\ldots=\a_0^\f$.

Similarly, if (\ref{eq:jan15-3}) holds, then by the same argument as above we could prove that \[d_{m+1}d_{m+2}\ldots=\theta(\Phi_1(\a_0))^\f, \] 
completing the proof.  
\end{proof}

\begin{proof}
[Proof of Proposition \ref{th:difference-Un+1-Un}]
We will prove the theorem by induction on $n$. Note by Lemma \ref{lem:difference-U1-U0} that it holds for $n=0$.  Now suppose the result holds for all $n\le k$ with $k\ge 0$, and we consider it for $n=k+1$. Let $(d_i)\in\widetilde{\U}_{\beta_{k+2}, M}\setminus\widetilde{\U}_{\beta_{k+1}, M}$. Then by Lemma \ref{lem:admissible-alphabet-3N+1} there exists $m_1\in\N$ such that   $d_{m_1}^1<M$ and 
$
d_{m_1+i}\in\set{\al_{(N+1)N}, \al_{N(N+1)}, \al_{NN}}$ for all $i\ge 1.
$
Note by the induction hypothesis that any sequence in $\widetilde{\U}_{\beta_{k+1}, M}\setminus\widetilde{\U}_{\beta_0, M}=\bigcup_{j=0}^k(\widetilde{\U}_{\beta_{j+1}, M}\setminus\widetilde{\U}_{\beta_j, M})$ must end with $\a_j^\f, \Phi_1(\a_j^\f)$ or their cyclic permutations under $\theta$ for some $0\le j\le k$.
Observe by (\ref{eq:relation-an-tn}) that $\a_j^1=\t_j$ for all $j\ge 0$, and $\t_0^\f\prec \t_1^\f\prec \cdots\prec \t_k^\f\prec\t_{k+1}^\f$. Since $\de(\beta_j)=\t_j^\f$ for all $j\ge 0$,  by (\ref{eq:an-3N+1})   it follows that the sequence $(d_i)\in\widetilde{\U}_{\beta_{k+2}, M}\setminus\widetilde{\U}_{\beta_{k+1}, M}$ must contain a block in $\bigcup_{j=0}^2\set{\theta^j(\a_{k}^+), \theta^j(\Phi_1(\a_{k}^+))}$. Without loss of generality we may assume 
$d_{m_2+1}\ldots d_{m_2+3\cdot 2^{k}}=\a_{k}^+$ for some $m_2\ge m_1$. Note that $d^1_{m_2}<M$. Then by Proposition \ref{prop:patterns-3N+1} it follows that 
\begin{equation}\label{eq:21-2}
	d_{m_2+1}\ldots d_{m_2+3\cdot 2^{k+1}}=\a_k^+\Phi_1(\a_k^+)=\a_{k+1}\quad\textrm{or}\quad \a_k^+\Phi(\a_k)=\a_{k+1}^+.
\end{equation}
If $d_{m_2+1}\ldots d_{m_2+3\cdot 2^{k+1}}=\a_{k+1}^+$, then by Proposition \ref{prop:patterns-3N+1} and using $(d_i)\in\widetilde{\U}_{\beta_{k+2}, M}$ it follows that $d_{m_2+1}d_{m_2+2}\ldots=\a_{k+2}^\f$. But by
(\ref{eq:relation-an-tn}) this gives 
\[
d_{m_2+1}^1d_{m_2+2}^1\ldots = \t_{k+2}^\f=\de(\beta_{k+2}),
\]
which   leads to a contradiction to $(d_i)\in\widetilde{\U}_{\beta_{k+2}, M}$.

Hence, by (\ref{eq:21-2}) we must have $d_{m_2+1}\ldots d_{m_2+3\cdot 2^{k+1}}=\a_{k+1}=\a_k^+\Phi_1(\a_k^+)$. Again, by Proposition \ref{prop:patterns-3N+1} it follows that 
\begin{equation}\label{eq:21-3}
	d_{m_2+3\cdot 2^{k+1}+1}\ldots d_{m_2+3(2^{k+1}+2^k)}=\a_k\quad\textrm{or}\quad \a_k^+.
\end{equation}
If $d_{m_2+3\cdot 2^{k+1}+1}\ldots d_{m_2+3(2^{k+1}+2^k)}=\a_k$, then $d_{m_2+3\cdot 2^k+1}\ldots d_{m_2+3\cdot(2^{k+1}+2^k)}=\Phi_1(\a_k^+)\a_k=\Phi(\a_{k+1}^+)$. By (\ref{eq:relation-an-tn}), and Proposition \ref{prop:lexicographical-characterization-U} and using $d^1_{m_2+3\cdot 2^k}>0$ it follows that 
\[
\overline{d^\+_{m_2+3\cdot 2^k+1}}\,\overline{d^\+_{m_2+3\cdot 2^k+2}}\ldots\lge \t_{k+2}^\f=\de(\beta_{k+2}),
\]
again leading to a contradiction to $(d_i)\in\widetilde{\U}_{\beta_{k+2}, M}$. 

So, by (\ref{eq:21-3}) we have $d_{m_2+3\cdot 2^{k+1}+1}\ldots d_{m_2+3(2^{k+1}+2^k)}=\a_k^+$. Continuing this argument as in (\ref{eq:21-2}) and (\ref{eq:21-3}) we can prove that 
$
d_{m_2+1}d_{m_2+2}\ldots =(\a_{k}^+\Phi_1(\a_k^+))^\f=\a_{k+1}^\f. 
$
This proves the result for $n=k+1$. 
\end{proof}

\begin{proof}
[Proof of Theorem \ref{th:main-result} for $M=3N+1$]
By Propositions \ref{prop:U-beta-coincidence}, \ref{prop:beta-G} and  \ref{prop:beta_c-3N+1} it suffices to prove (ii)--(iv) for $\widetilde{U}_{\beta, M}$. Since the map $\Pi_\beta$ from $\widetilde{\U}_{\beta, M}$ to $\widetilde{U}_{\beta, M}$ is bijective and bi-H\"older continuous, we only need to consider the symbolic intrinsic univoque set $\widetilde{\U}_{\beta, M}$. Here we equip with the symbolic space $\Om^\N_M$ the metric $\rho$ given by 
\begin{equation}\label{eq:rho}
	\rho((c_i), (d_i))=2^{-\inf\set{n: c_n\ne d_n}}
\end{equation}
for any $(c_i), (d_i)\in\Om_M^\N$.

First we prove (ii). Take $\beta\in(\beta_G, \beta_c)$. Note that $\beta_0=\beta_G$ and  $\beta_k\nearrow\beta_c$ as $k\to\f$. So there exists $k\in\N$ such that $\beta<\beta_k$. Note by Propositions \ref{prop:beta-G} and  \ref{th:difference-Un+1-Un} that $\widetilde{\U}_{\beta_k, M}=\widetilde{\U}_{\beta_0, M}\cup\bigcup_{j=0}^{k-1}\widetilde{\U}_{\beta_{j+1}, M}\setminus\widetilde{\U}_{\beta_{j}, M}$ is at most countable. This implies that $\widetilde{\U}_{\beta, M}\subset\widetilde{\U}_{\beta_k, M}$ is at most countable. On the other hand, by Proposition \ref{prop:lexicographical-characterization-U} we have
\[
\al_{00}^n(\al_{(N+1)N}\al_{N(N+1)}\al_{NN})^\f\in\widetilde{\U}_{\beta, M}\quad\forall n\in\N.
\]
So, $\widetilde{\U}_{\beta, M}$ is indeed countably infinite, establishing (ii).

Next we prove (iii). The proof is similar to that of \cite[Proposition 5.6]{Kong-Li-2020}. Note that $\beta_k\nearrow\beta_c$ as $k\to\f$. By Proposition \ref{prop:patterns-3N+1} it follows that $\widetilde{\U}_{\beta_c, M}$ contains all of the following sequences
\[
\a_1^{j_1}\a_2^{j_2}\cdots\a_k^{j_k}\cdots,\quad j_k\in\N.
\]
Note that each block $\a_k$ can not be written as concatenations of two or more blocks of the form $\a_l$ with $l<k$. This implies that $\widetilde{\U}_{\beta_c, M}$ is uncountable. On the other hand, by the recursive construction as in the proof of Proposition \ref{prop:patterns-3N+1} it follows that any sequence in $\widetilde{\U}_{\beta_c, M}\setminus\widetilde{\U}_{\beta_G, M}$ must be of the form
\[
\omega(\a_{i_1}^+\Phi_1(\a_{i_1}^+))^{j_1}\left(\a_{i_1}^+\Phi_1(\a_{i_2'}^+)\right)^{j_1'}\cdots(\a_{i_k}^+\Phi_1(\a_{i_k}^+))^{j_k}\left(\a_{i_k}^+\Phi_1(\a_{i_{k+1}'}^+)\right)^{j_k'}\cdots 
\]
or its cyclic permutations under $\theta$, where $\omega\in\Omega_M^*$, $j_k\in\N\cup\set{0,\f}$, $j_k'\in\set{0,1}$ and
\[
1\le i_1<i_2'\le i_2<\cdots<i_k'\le i_k<i_{k+1}'\le\cdots.
\]
Observe that the length of $\a_k$ grows exponentially fast. This implies that $\dim_H\widetilde{\U}_{\beta_c, M}=0$, establishing (iii).

Finally we prove (iv). For $k\in\N_0$ note that $\t_k^+=\la_1\ldots \la_{3\cdot 2^k}$. Then by Lemma \ref{lem:inequality-lambda-i} it follows that $\si^i(\t_k^+\re_1(\t_k)^\f)\lle \t_k^+\re_1(\t_k)^\f$ for all $i\ge 0$. Then by Lemma \ref{lem:char-quasi-greedy-expansion} there exists $\hat\beta_k\in(N+1, N+2)$ such that
\[
\de(\hat\beta_k)=\t_k^+\re_1(\t_k)^\f.
\]
Since $\de(\hat\beta_k)\searrow\de(\beta_c)$ as $k\to\f$, again by Lemma \ref{lem:char-quasi-greedy-expansion} it follows that
$\hat\beta_k\searrow\beta_c$ as $k\to\f$. So, for $\beta>\beta_c$ there exists a smallest $k\in\N_0$ such that $\hat\beta_k<\beta$. Then  we can show that $\widetilde{\U}_{\beta, M}$ contains a subshift of finite type $(X_k, \sigma)$ represented  by the labeled graph $G$ (see Figure \ref{figure:15}).
\begin{figure}[h!]
	\centering
	\begin{tikzpicture}[->,>=stealth',shorten >=1pt,auto,node distance=4cm,
		semithick]
		
		\tikzstyle{every state}=[minimum size=0pt,fill=black,draw=none,text=black]
		
		\node[state] (A)                    { };
		\node[state]         (B) [ right of=A] { };

		\path[->,every loop/.style={min distance=0mm, looseness=80}]
		(A) edge [loop left,->]  node {$\Phi_1(\a_k)$} (A)
		edge  [bend left]   node {$\Phi_1(\a_k^+)$} (B)
		
		(B) edge [loop right] node {$\a_k$} (B)
		edge  [bend left]            node {$\a_k^+$} (A);
	\end{tikzpicture}
	\caption{The   labeled graph $G$ presenting $(X_k, \si)$.}\label{figure:15}
\end{figure}
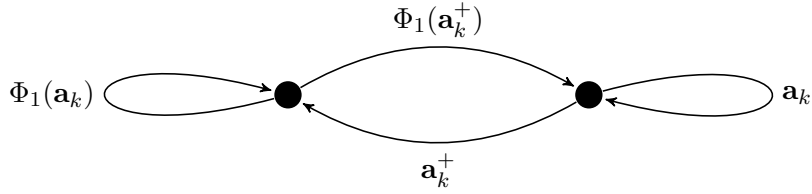
From this we can deduce that 
\[\dim_H\widetilde{\U}_{\beta, M}\ge \dim_H X_k=\frac{\log 2}{|\a_k|\log 2}=\frac{1}{|\a_k|}>0,\]
where we use the metric $\rho$ defined in (\ref{eq:rho}). This
completes  the proof.
\end{proof}
At the end of this section we give an example to illustrate that the Hausdorff dimension of $\widetilde{U}_{\beta, M}$ can be explicitly calculated for some $\beta>\beta_c(M)$.
\begin{example}\label{ex:1}
Let $\beta=\frac{1+\sqrt{5}}{2}$ and $M=1$. Then $\beta>\beta_c(1)$. Observe by Proposition \ref{prop:lexicographical-characterization-U} that   up to   countably many points the set $\widetilde{\U}_{\beta, 1}$ is a subshift of finite type in $\set{\al_{00}, \al_{10}, \al_{01}}^\N$ with the set of forbidden words 
\[
\mathcal F=\set{\al_{10}\al_{10}, \al_{01}\al_{01}, \al_{00}\al_{00}}. 
\]
Then it can be represented as a directed graph with the vertex set $V=\set{\al_{00}, \al_{10}, \al_{01}}$ and the adjacency matrix 
\[
A=\left(\begin{array}{ccc}
	0 & 1 & 1 \\
	1 & 0 & 1 \\
	1 & 1 & 0 
\end{array}\right).
\]
So, $\dim_H\widetilde{U}_{\beta, 1}=\frac{\log\rho(A)}{\log\beta}=\frac{\log 2}{\log\beta}\approx 1.44042$, where $\rho(A)$ denotes the spectral radius of $A$. Indeed, we also have $\dim_H U_{\beta, 1}=\frac{\log 2}{\log\beta}$. 
\end{example}
We remark that if $\beta$ is a multinacci number, then by a similar argument as in Example \ref{ex:1} we are able to determine the explicit Hausdorff dimension of $U_{\beta, 1}$.

\section{Generalized Komornik-Loreti constant $\beta_c(M)$ for $M=3N+2$}\label{sec:3N+2}

Let $M=3N+2$ with $N\in\N_0$. Recall $\beta_c(M)$ from Definition \ref{def:beta-c}.  
We will show in this section that $\beta_c(M)$ is the critical base for $\widetilde{\U}_{\beta, M}$, i.e., $\widetilde{\U}_{\beta, M}$ has positive Hausdorff dimension if and only if $\beta>\beta_c(M)$. The proof is similar to that for $M=3N+1$, so we only sketch the main idea in the proof of Theorem \ref{th:main-result} for $M=3N+2$.

First we introduce  a     Thue-Morse  type sequence $(\ga_i)\in\set{N, N+1}^\N$, which will be the quasi-greedy $\beta_c(M)$-expansion of $1$. In view of Definition \ref{def:beta-c}, we define the sequence $(\ga_i)$ by 
\begin{equation}\label{eq:relation-gamma-tau}
\ga_{3k+1}=\tau_{2k+1}+N,\quad \ga_{3k+2}=N+1, \quad \ga_{3k+3}=\tau_{2k+2}+N,
\end{equation}
where  $(\tau_i)_{i=0}^\f=0110100110010110\ldots\in\set{0,1}^\N$ is the classical Thue-Morse sequence. Then $(\ga_i)$ begins with $(N+1)(N+1)(N+1)N(N+1)(N+1)N(N+1)N(N+1)(N+1)(N+1)$. Similar to Lemma \ref{lem-obs-(3N+1)},  the sequence $(\ga_i)$ can also be obtained by substitutions. Let $
\varphi_2: 0\mapsto N(N+1)(N+1),~ 1\mapsto (N+1)N(N+1).
$
Then 
\[(\ga_i)_{i=1}^\f=\si(\varphi_2(\tau_0)\varphi_2(\tau_1)\varphi_2(\tau_2)\ldots),\]
where $\si$ is the left-shift map.

Similar to $(\la_i)$ defined in the previous section, we construct a sequence of blocks $(\su_n)$ approximating $(\ga_i)$ componentwisely. Let 
\[{\Sigma}_2:=\set{(N+1)(N+1)(N+1), N(N+1)N,(N+1)(N+1)N, N(N+1)(N+1)},\] and we define a block  map 
$ \re_2: {\Sigma}_2  \ra {\Sigma}_2$ by
\begin{equation*}\label{eq:substitution-3N+2}
\begin{split}
	(N+1)(N+1)(N+1)\quad &\longleftrightarrow\quad N(N+1)N, \\
	(N+1)(N+1)N   \quad&\longleftrightarrow\quad N(N+1)(N+1).
	\end{split}\end{equation*}
	Then the sequence  $(\su_n)$ is defined recursively as follows: 
	\begin{equation}\label{eq:un-3N+2}
\su_0=(N+1)(N+1)N,\quad\textrm{and}\quad \su_{n+1}:=\su_n^+\re_2(\su_n^+)\quad\forall n\ge 0.
\end{equation}
For example,
\begin{align*}
\su_1 & =(N+1)(N+1)(N+1)N(N+1)N, \\
\su_2 & = (N+1)(N+1)(N+1)N(N+1)(N+1)\,N(N+1)N(N+1)(N+1)N.
\end{align*} 
Note that $(\ga_i)$ begins with $\su_n^+$ for any $n\ge 1$. 
Then the sequence $(\su_n)$ converges to $(\ga_i)$ componentwisely, i.e., 
\begin{equation*}\label{eq:gamma-i}
(\ga_i)=\lim_{n\to\f}\su_n \in\Sigma_2^\N.
\end{equation*}

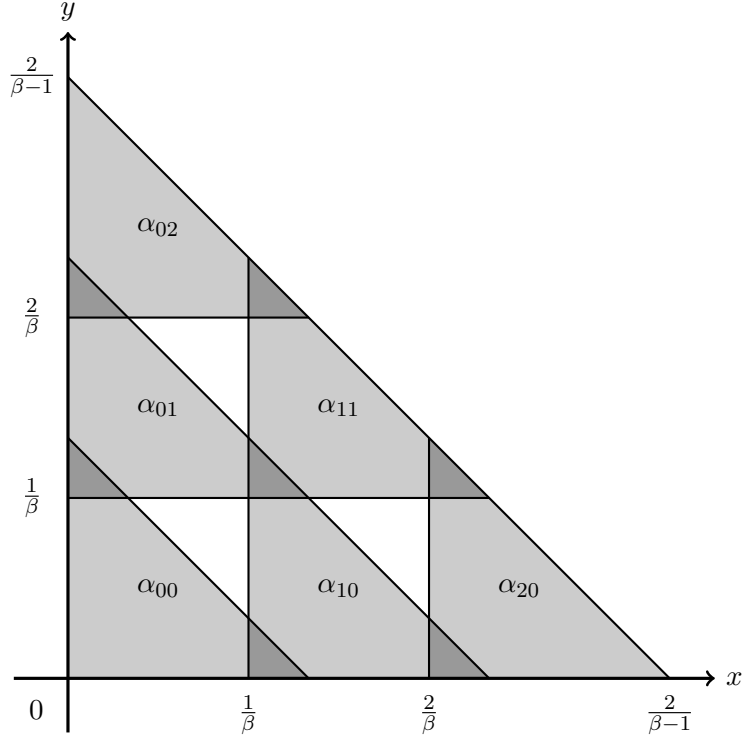
\begin{figure}[h!]
\begin{center}
	\begin{tikzpicture}[
		scale=6,  
		axis/.style={very thick, ->},
		important line/.style={thick},
		dashed line/.style={dashed, thin},
		pile/.style={thick, ->, >=stealth', shorten <=2pt, shorten
			>=2pt},
		every node/.style={color=black}
		]
		
		\fill[black!20]({0},0)--({2/(2.5-1)-2/2.5},0)--(0,{2/(2.5-1)-2/2.5})--cycle; 
		\fill[black!20]({1/2.5},0)--({2/(2.5-1)-1/2.5},0)--({1/2.5},{2/(2.5-1)-2/2.5})--cycle; 
		\fill[black!20]({2/2.5},0)--({2/(2.5-1)},0)--({2/2.5},{2/(2.5-1)-2/2.5})--cycle;
		
		\fill[black!20]({0},{1/2.5})--({2/(2.5-1)-2/2.5},{1/2.5})--(0,{2/(2.5-1)-1/2.5})--cycle; 
		\fill[black!20]({0},{2/2.5})--({2/(2.5-1)-2/2.5},{2/2.5})--(0,{2/(2.5-1)})--cycle; 
		\fill[black!20]({1/2.5},{1/2.5})--({2/(2.5-1)-1/2.5},{1/2.5})--({1/2.5},{2/(2.5-1)-1/2.5})--cycle;

		\fill[black!40]({1/2.5},0)--({2/(2.5-1)-2/2.5},0)--({1/2.5},{2/(2.5-1)-3/2.5})--cycle;    
		\fill[black!40]({2/2.5},0)--({2/(2.5-1)-1/2.5},0)--({2/2.5},{2/(2.5-1)-3/2.5})--cycle;
		
		\fill[black!40]({1/2.5},{1/2.5})--({2/(2.5-1)-2/2.5},{1/2.5})--({1/2.5},{2/(2.5-1)-2/2.5})--cycle;    
		\fill[black!40]({2/2.5},{1/2.5})--({2/(2.5-1)-1/2.5},{1/2.5})--({2/2.5},{2/(2.5-1)-2/2.5})--cycle;
		
		\fill[black!40](0,{1/2.5})--(0,{2/(2.5-1)-2/2.5})--({2/(2.5-1)-3/2.5}, {1/2.5})--cycle;    
		\fill[black!40](0,{2/2.5})--(0,{2/(2.5-1)-1/2.5})--({2/(2.5-1)-3/2.5}, {2/2.5})--cycle;
		
		\fill[black!40]({1/2.5},{2/2.5})--({2/(2.5-1)-2/2.5},{2/2.5})--({1/2.5},{2/(2.5-1)-1/2.5})--cycle;

		\draw[axis] (-0.12,0)  -- ({(2/(2.5-1))+0.1},0) node(xline)[right]
		{$x$};
		\draw[axis] (0,-0.12) -- (0,{2/(2.5-1)+0.1}) node(yline)[above] {$y$};
		\node[] at (-0.07,-0.07){$0$};
		\draw[important line] (0,{2/(2.5-1)})--({2/(2.5-1)},0);
		
		\draw[important line] ({1/2.5},0)--({1/2.5},{2/(2.5-1)-1/2.5});
		\draw[important line] ({2/2.5},0)--({2/2.5},{2/(2.5-1)-2/2.5});
		
		\draw[important line] (0,{1/2.5})--({2/(2.5-1)-1/2.5},{1/2.5});
		\draw[important line] (0,{2/2.5})--({2/(2.5-1)-2/2.5},{2/2.5});
		
		\draw[important line] ({2/(2.5-1)-2/2.5},0)--(0,{2/(2.5-1)-2/2.5});
		\draw[important line] ({2/(2.5-1)-1/2.5},0)--(0,{2/(2.5-1)-1/2.5});
		
		\node[] at ({1/(2*2.5)}, {1/(2*2.5)}){$\al_{00}$};
		
		\node[] at ({3/(2*2.5)}, {1/(2*2.5)}){$\al_{10}$};
		
		\node[] at ({5/(2*2.5)}, {1/(2*2.5)}){$\al_{20}$};
		
		\node[] at ({1/(2*2.5)}, {3/(2*2.5)}){$\al_{01}$};
		
		\node[] at ({1/(2*2.5)}, {5/(2*2.5)}){$\al_{02}$};
		
		\node[] at ({3/(2*2.5)}, {3/(2*2.5)}){$\al_{11}$};
		
		\node[] at ({1/(2.5)}, {-0.08}){$\frac{1}{\beta}$};
		
		\node[] at ({2/(2.5)}, {-0.08}){$\frac{2}{\beta}$};
		\node[] at ({2/(2.5-1)}, {-0.08}){$\frac{2}{\beta-1}$};
		
		\node[] at ( {-0.08}, {1/(2.5)}){$\frac{1}{\beta}$};
		
		\node[] at ( {-0.08}, {2/(2.5)}){$\frac{2}{\beta}$};
		\node[] at ( {-0.08}, {2/(2.5-1)}){$\frac{2}{\beta-1}$};

	\end{tikzpicture}
\end{center}
\caption{The graph for the first generation of $S_{\beta, M}$ with $M=2$ and $\beta=2.5$.   The convex hull $\Delta_{\beta, M}$ is the triangle with vertices $(0, 0), (2/(\beta-1), 0)$ and $(0, 2/(\beta-1))$. Each light grey triangle corresponds to a $f_{\beta, \al_{ij}}(\Delta_{\beta,M})$ for some $\al_{ij}\in\Omega_M$.  Furthermore, the overlap region $O_{\beta, M}$ is the union of   seven  small dark grey triangles.}\label{Fig:S-2}
\end{figure}

We also define the \emph{conjugate} of $(\ga_i)$ by 
\[
(\hat\ga_i)_{i=1}^\f:=\re_2(\ga_1\ga_2\ga_3)\re_2(\ga_4\ga_5\ga_6)\ldots\re_2(\ga_{3k+1}\ga_{3k+2}\ga_{3k+3})\ldots.
\]
Then  by  (\ref{eq:relation-gamma-tau}) and (\ref{eq:un-3N+2}) it follows that for any $k\in\N_0$,
\begin{equation}\label{eq:relation-gamma-tau-conjugacy}
\hat\ga_{3k+1}=\overline{\tau_{2k+1}}+N,\quad \hat\ga_{3k+2}=N+1,\quad \hat\ga_{3k+3}=\overline{\tau_{2k+2}}+N,
\end{equation}
where $\overline{\tau_i}=1-\tau_i$ for any $i\ge 0$.

By a similar argument as in the proof of Lemma \ref{lem:inequality-lambda-i} we have the following result.
\begin{lemma}
\label{lem:inequality-gamma-i}
For any $n\in\N_0$ and any $0\le i<3\cdot 2^n$, the following inequalities hold:
\begin{align*}
	\hat\gamma_1\ldots \hat\gamma_{3\cdot 2^n-i}&\prec \gamma_{i+1}\ldots \gamma_{3\cdot 2^n}\lle \gamma_1\ldots\gamma_{3\cdot 2^n-i},\\
	\hat\gamma_1\ldots\hat\gamma_{3\cdot 2^n-i}&\lle \hat\gamma_{i+1}\ldots\hat\gamma_{3\cdot 2^n}\prec \gamma_1\ldots \gamma_{3\cdot 2^n-i}.
\end{align*}
\end{lemma}
By Lemmas  \ref{lem:inequality-gamma-i} and \ref{lem:char-quasi-greedy-expansion} it follows that $(\ga_i)$ is the quasi-greedy $\beta_c(M)$-expansion of $1$, i.e., \begin{equation*}\label{eq:quasi-greedy-betac-3N+2}
\begin{split}
	\de(\beta_c(M))&=\ga_1\ga_2\ldots =(N+1)(N+1)(N+1)N(N+1)(N+1)\;N(N+1)N \ldots,
\end{split}
\end{equation*}
By  Lemma \ref{lem:mahler-transcendental},  (\ref{eq:relation-gamma-tau}) and the same argument as in the proof of Proposition \ref{prop:beta_c-3N+1}  we can prove that $\beta_c(3N+2)$ is transcendental.  

Note by (\ref{eq:un-3N+2}) that $\su_n=\ga_1\ldots \ga_{3\cdot 2^{n}}^-$. Then by Lemma \ref{lem:inequality-gamma-i} we have 
$
\si^i(\su_n^\f)\lle \su_n^\f$ for all  $i\ge 0.
$
So, by Lemma \ref{lem:char-quasi-greedy-expansion} there exists a unique $\beta_n=\beta_n(M)\in(1, M+1]$ such that 
\begin{equation}\label{eq:beta-n-3N+2}
\de(\beta_n)=\su_n^\f.
\end{equation}
Note by Proposition \ref{prop:beta-G} that 
$
\de(\beta_G(M))=((N+1)(N+1)N)^\f=\su_0^\f=\de(\beta_0(M)).
$
Then $\beta_G(M)=\beta_0(M)$ by Lemma \ref{lem:char-quasi-greedy-expansion}.     
Observe that $\de(\beta_n)=\su_n^\f\nearrow (\ga_i)=\de(\beta_c)$ as $n\to\f$. By Lemma \ref{lem:char-quasi-greedy-expansion} it follows that 
\[
\beta_G=\beta_0<\beta_1<\cdots<\beta_n<\beta_{n+1}<\cdots,\quad\textrm{and}\quad\beta_n\nearrow\beta_c\quad\textrm{as }n\to\f. 
\]

To prove Theorem \ref{th:main-result} for $M=3N+2$ we will construct a sequence of blocks $\set{\b_n}_{n=0}^\f$ in $\Om_M^*$. Similar to the words $(\a_n)$ in the previous section, these words $(\b_n)$ are like basic bricks to describe the set $\widetilde{\U}_{\beta, M}$ for $\beta\le\beta_c(M)$ with $M=3N+2$.
Let 
\begin{align*}
\pazocal A_2 & :=\left\{\al_{(N+1)(N+1)}\al_{(N+1)N}\al_{N(N+1)}, \quad \al_{N(N+1)}\al_{(N+1)N}\al_{(N+1)(N+1)},\right. \\
&\hspace{0.75cm}\left.\al_{(N+1)(N+1)}\al_{(N+1)N}\al_{(N+1)(N+1)}, \quad\al_{N(N+1)}\al_{(N+1)N}\al_{N(N+1)}\right\}.
\end{align*} 
We define the substitution $\Phi_2: \pazocal A_2\to\pazocal A_2$ by
\begin{equation}
\label{eq:substitution-Phi-3N+2}
\begin{split}
	\al_{(N+1)(N+1)}\al_{(N+1)N}\al_{N(N+1)}\quad&\longleftrightarrow \quad  \al_{N(N+1)}\al_{(N+1)N}\al_{(N+1)(N+1)}, \\
	\al_{(N+1)(N+1)}\al_{(N+1)N}\al_{(N+1)(N+1)}\quad & \longleftrightarrow\quad \al_{N(N+1)}\al_{(N+1)N}\al_{N(N+1)}.
\end{split}
\end{equation}
Then the sequence $\set{\b_n}_{n=0}^\f\subset\pazocal A_2^*$ is defined as follows (see Figure \ref{Fig:S-2}): 
\begin{equation}\label{eq:bn-3N+2}
\b_0=\al_{(N+1)(N+1)}\al_{(N+1)N}\al_{N(N+1)},\quad\textrm{and}\quad  \b_{n+1}=\b_n^+\Phi_2(\b_n^+)\quad\forall n\ge 0,
\end{equation}
where $\b_n^+$ is the word in $\pazocal A_2^*$ by changing the last digit of $\b_n$ from $\al_{N(N+1)}$ to $\al_{(N+1)(N+1)}$.
Then each $\b_n$ has length $3\cdot 2^n$ and ends with $\al_{N(N+1)}$.  If we write $\b_n=b_1\ldots b_{3\cdot 2^n}$, then $b_{3k+2}=\al_{(N+1)N}$ for all $0\le k<2^n$. So, 
\[
\b_0\sim\left(\begin{array}{ccc}
N+1&N+1&N \\
N+1&N&N+1 \\
N&N+1&N+1
\end{array}\right), 
\]
and
\[
\b_1\sim \left(\begin{array}{cccccc}
N+1&N+1&N+1&N&N+1&N \\
N+1&N&N+1&N+1&N&N+1 \\
N&N+1&N&N+1&N+1&N+1
\end{array}\right).
\]
Comparing the definitions of $(\su_n)$ and $(\b_n)$ in  (\ref{eq:un-3N+2}) and (\ref{eq:bn-3N+2}) respectively,  it follows that 
\begin{equation}\label{eq:relation-bn-un}
\b_n^1=\su_n,\quad \b_n^2=((N+1)N (N+1))^n\quad\textrm{and}\quad \overline{\b_n^\+}=\re_2(\su_n)
\end{equation}
for any $n\ge 0$. Recall that for  a block $\d=d_1d_2\ldots d_k\in\Omega_M^*$ we write $\d^1=d_1^1d_2^1\ldots d_k^1, \d^2=d_1^2d_2^2\ldots d_k^2$ and $\overline{\d^\+}=\overline{d_1^\+}\,\overline{d_2^\+}\ldots\overline{d_k^\+}$, where each $d_i=(d_i^1, d_i^2)$ and $\overline{d_i^\+}=M-d_i^\+$ with $d_i^\+=d_i^1+d_i^2$.

The following result  can be proved in a similar way to that for Proposition \ref{th:difference-Un+1-Un}. 
\begin{proposition}
\label{th:difference-Un+1-Un-3N+2}
Let  $n\in\N_0$. Then any sequence in $\widetilde{\U}_{\beta_{n+1},M}\setminus\widetilde{\U}_{\beta_n,M}$ must end in
\[
\bigcup_{j=0}^2\set{  \theta^j(\b_n^\f), \theta^j(\Phi_2(\b_n)^\f)},
\]
where $\theta$ is the cyclic permutation defined in (\ref{eq:theta}).
\end{proposition}
It is worth  mentioning  that for $M=3N+2$, if      $(d_i)\in\widetilde{\U}_{\beta_c, M}\setminus\widetilde{\U}_{\beta_0,M}$, then there exists $m\in\N$ such that $d_m^1<M, d_m^2<M, \overline{d_m^\+}<M$, and    
\[
d_k\in\bigcup_{j=0}^2\set{\theta^j(\al_{(N+1)(N+1)})}\quad\forall k>m.
\]
So, in this case we don't need a lemma like Lemma \ref{lem:exceptional-digits}. All the other properties are   similar to those for $M=3N+1$. For example, we have a similar characterization on patterns appearing in sequences of $\widetilde{\U}_{\beta_c,M}$.
\begin{proposition}
\label{prop:patterns-3N+2}
Let $(d_i)\in\widetilde{\U}_{\beta_c,M}$. Then   the following statements hold.
\begin{enumerate}[{\rm(i)}]
	\item If $d^1_m<M$ and $d_{m+1}\ldots d_{m+3\cdot 2^n}=\b_n^+$, then 
	$
	d_{m+3\cdot 2^n+1}\ldots d_{m+3\cdot 2^{n+1}}\in\set{\Phi_1(\b_n), \Phi_1(\b_n^+)}.
	$
	\item If $d^1_m>0$ and $d_{m+1}\ldots d_{m+3\cdot 2^n}=\Phi_1(\b_n^+)$, then 
	$
	d_{m+3\cdot 2^n+1}\ldots d_{m+3\cdot 2^{n+1}}\in\set{\b_n, \b_n^+}.
	$
\end{enumerate}
\end{proposition}
The proof of Proposition \ref{prop:patterns-3N+2} is very similar to that of Proposition \ref{prop:patterns-3N+1}, and we leave it as an exercise for interested readers.

\begin{proof}
[Proof of Theorem \ref{th:main-result} for $M=3N+2$]
The proof is similar to that for $M=3N+1$. By using Proposition  \ref{th:difference-Un+1-Un-3N+2} we can deduce that for $\beta<\beta_c$ the set $\widetilde{\U}_{\beta, M}$ is at most countable; and by Proposition \ref{prop:patterns-3N+2} we can prove that for $\beta=\beta_c$ the set $\widetilde{\U}_{\beta, M}$ is uncountable but has zero Hausdorff dimension. Furthermore, for $\beta>\beta_c$ we can construct a subshift of finite type $Y_k\subset\widetilde{\U}_{\beta, M}$ having positive Hausdorff dimension, where $Y_k$ is analogous to $X_k$ as described in Figure \ref{figure:15} with $\a_k, \a_k^+, \Phi_1(\a_k)$ and $\Phi_1(\a_k^+)$ replaced by $\b_k, \b_k^+, \Phi_2(\b_k)$ and $\Phi_2(\b_k^+)$, respectively. 
\end{proof}

\section{Generalized Komornik-Loreti constant for $\widetilde{\U}_{\beta, M}$ with $M=3N+3$}\label{sec:3N}
Let $M=3N+3$ with $N\in\N_0$. Recall that $\beta_c(M)$ is defined in Definition \ref{def:beta-c}.  In this section we will prove Theorem \ref{th:main-result} for $M=3N+3$, i.e., $\widetilde{\U}_{\beta, M}$ has positive Hausdorff dimension if and only if $\beta>\beta_c(M)$. Completely different from the case for $M=1$, we need more effort to find the critical base $\beta_c(M)$ for $M=3N+3$.

In view of Definition \ref{def:beta-c} and the definition of the classical Thue-Morse sequence  $(\tau_i)$,
we introduce  a     Thue-Morse  type sequence $(\eta_i)\in\set{N, N+1, N+2}^\N$: 
\begin{equation}\label{eq:relation-eta-tau}
\eta_{3k+1}=N+2-2\tau_k,\quad \eta_{3k+2}=N+\tau_k, \quad \eta_{3k+3}=N+1+\tau_{k+1}\quad\forall k\in\N_0.
\end{equation}
Similar to Lemma \ref{lem-obs-(3N+1)}, the sequence $(\eta_i)$ can also be obtained  by substitutions. Let 
$
\varphi_3: 0\mapsto (N+1)(N+2)N, ~1\mapsto (N+2)N(N+1).
$
Then 
\[(\eta_i)_{i=1}^\f=\si(\varphi_3(\tau_0)\varphi_3(\tau_1)\varphi_3(\tau_2)\ldots),\]
where $\si$ is the left-shift map. 
\begin{proposition}
\label{prop:transcendental-betac-3N}
Let $M=3N+3$ with $N\in\N_0$. Then the quasi-greedy $\beta_c(M)$-expansion of $1$ is given by
\[
\de(\beta_c(M))=\eta_1\eta_2\ldots=(N+2)N(N+2)N(N+1)(N+2)N(N+1)(N+1)(N+2)N(N+1)\ldots.
\]
Furthermore, $\beta_c(M)$ is transcendental.
\end{proposition}

Let  
\[{\Sigma}_3:=\set{(N+2)N(N+2),  N(N+1)(N+1), (N+2)N(N+1), N(N+1)(N+2)},\]
and we define the block map  $\re_3: {\Sigma}_3  \ra {\Sigma}_3$ by
\begin{equation}\label{eq:substitution-3N}
\begin{split}
	(N+2)N(N+2) \quad & \longleftrightarrow\quad  N(N+1)(N+1),\\
	(N+2)N(N+1) \quad  &\longleftrightarrow\quad N(N+1)(N+2).
\end{split}
\end{equation}
Accordingly,   we   define the \emph{conjugate} of $(\eta_i)$ by 
\[
(\hat\eta_i)_{i=1}^\f:=\re_3(\eta_1\eta_2\eta_3)\re_3(\eta_4\eta_5\eta_6)\ldots\re_3(\eta_{3k+1}\eta_{3k+2}\eta_{3k+3})\ldots.
\]
Then $(\hat\eta_i)\in\Sigma_3^\N$ begins with $N(N+1)(N+1)(N+2)N(N+1)$. Furthermore, by the definition of $\phi_3$ it follows that for any $k\in\N_0$,
\begin{equation}\label{eq:relation-eta-tau-conjugacy}
\hat\eta_{3k+1}=N+2-2\overline{\tau_k},\quad \hat\eta_{3k+2}=N+\overline{\tau_k},\quad \hat\eta_{3k+3}=N+1+\overline{\tau_{k+1}},
\end{equation}
where $\overline{\tau_i}=1-\tau_i$ for all $i\ge 0$.
\begin{lemma}
\label{lem:inequality-eta-i}
For any $n\in\N_0$ and any $0\le i<3\cdot 2^n$, the following inequalities hold:
\begin{align}
	\hat\eta_1\ldots \hat\eta_{3\cdot 2^n-i}&\prec \eta_{i+1}\ldots \eta_{3\cdot 2^n}\lle \eta_1\ldots\eta_{3\cdot 2^n-i},\label{eq:inequ-1-3N}\\
	\hat\eta_1\ldots\hat\eta_{3\cdot 2^n-i}&\lle \hat\eta_{i+1}\ldots\hat\eta_{3\cdot 2^n}\prec \eta_1\ldots \eta_{3\cdot 2^n-i}.\label{eq:inequ-2-3N}
\end{align}
\end{lemma}
\begin{proof}
Since the proof of (\ref{eq:inequ-2-3N}) is similar,   we only prove (\ref{eq:inequ-1-3N}). Take $i\in\set{0,1,\ldots, 3\cdot 2^n-1}$. We consider the  following three cases: (I) $i=3k$, (II) $i=3k+1$, and (III) $i=3k+2$ for some $k\in\set{0,1,\ldots, 2^n-1}$.

Case (I). $i=3k$ with $0\le k<2^n$. If $k=2^n-1$, then by (\ref{eq:relation-eta-tau})   we have 
\begin{align*}
	\eta_{i+1}\eta_{i+2}\eta_{i+3}&=(N+2-2\tau_{2^n-1})(N+\tau_{2^n-1})(N+1+\tau_{2^n})\\
	&\lle (N+2)N(N+2)=\eta_1\eta_2\eta_3.
\end{align*}
On the other hand, by (\ref{eq:relation-eta-tau-conjugacy}) and using $\tau_{2^n}=1$ it follows that 
\[
\eta_{i+1}\eta_{i+2}\eta_{i+3}\lge N(N+1)(N+2)\succ N(N+1)(N+1)=\hat\eta_1\hat\eta_2\hat\eta_3.
\]
So, in the following we assume $0\le k<2^n-1$. Note by (\ref{eq:relation-eta-tau}) that 
\begin{equation}
	\label{eq:29-1}
	\eta_{i+1}\ldots \eta_{3\cdot 2^n}=(N+2-2\tau_k)(N+\tau_k)\mathbf B(\tau_{k+1})\cdots \mathbf B(\tau_{2^n-1})(N+1+\tau_{2^n}),
\end{equation}
where 
\[\mathbf B(\tau_j):=(N+1+\tau_j)(N+2-2\tau_j)(N+\tau_j)\quad \textrm{for }j\ge 1.\]
Observe by (\ref{eq:relation-eta-tau}) and (\ref{eq:relation-eta-tau-conjugacy}) that 
\begin{equation}\label{eq:29-2}
	\begin{split}
		\eta_1\ldots \eta_{3\cdot 2^n} & =(N+2)N\mathbf B(\tau_1)\cdots \mathbf B(\tau_{2^n-1})(N+1+\tau_{2^n}), \\
		\hat\eta_1\ldots \hat\eta_{3\cdot 2^n}&=N(N+1) \mathbf B(\overline{\tau_1})\cdots \mathbf B(\overline{\tau_{2^n-1}})(N+1+\overline{\tau_{2^n}}). 
	\end{split}
\end{equation}
Since $\mathbf B(\tau_i)\prec \mathbf B(\tau_j)$ is equivalent to $\tau_i<\tau_j$, by (\ref{eq:29-1}) and (\ref{eq:29-2}) it follows that 
\begin{align*}
	&\quad\eta_{i+1}\ldots \eta_{3\cdot 2^n}(N+2-2\tau_{2^n})(N+\tau_{2^n})\\
	& \lle (N+2)N\mathbf B(\tau_{k+1})\cdots \mathbf B(\tau_{2^n}) \\
	& \lle (N+2)N\mathbf B(\tau_1)\cdots \mathbf B(\tau_{2^n-k})=\eta_1\ldots \eta_{3\cdot 2^n-i+2},
\end{align*}
where the second inequality follows by using the property of $(\tau_i)$, i.e., $\tau_{k+1}\ldots \tau_{2^n}\lle \tau_1\ldots \tau_{2^n-k}$. This proves the second inequality in (\ref{eq:inequ-1-3N}).

On the other hand, by (\ref{eq:29-1}) and (\ref{eq:29-2}) we obtain that 
\begin{equation}\label{eq:29-3}
	\begin{split}
		&\quad \eta_{i+1}\ldots \eta_{3\cdot 2^n}(N+2-2\tau_{2^n})(N+\tau_{2^n})\\
		& \lge N(N+1) \mathbf B(\tau_{k+1})\cdots \mathbf B(\tau_{2^n}) \\
		& \succ N(N+1) \mathbf B(\overline{\tau_1})\cdots \mathbf B(\overline{\tau_{2^n-k}})=\hat\eta_1\ldots \hat\eta_{3\cdot 2^n-i+2},
	\end{split}
\end{equation}
where the second inequality follows by using the fact that $\tau_{k+1}\ldots \tau_{2^n}\succ \overline{\tau_1}\ldots \overline{\tau_{2^n-k}}$. Note that $(N+2-2\tau_{2^n})(N+\tau_{2^n})=N(N+1)\lle\hat\eta_j\hat\eta_{j+1}$ for all $j\ge 1$. So, by (\ref{eq:29-3}) it follows that 
\[
\eta_{i+1}\ldots \eta_{3\cdot 2^n}\succ \hat\eta_1\ldots \hat\eta_{3\cdot 2^n-i},
\]
establishing (\ref{eq:inequ-1-3N}) when $i=3k$.

Case (II). $i=3k+1$ with $0\le k<2^n$. Then $\eta_{i+1}=N+\tau_k\le N+1<\eta_1$, which yields the second inequality in (\ref{eq:inequ-1-3N}). On the other hand, if $k=2^n-1$ then by using $\tau_{2^n}=1$ it follows that
\[
\eta_{i+1}\eta_{i+2}=(N+\tau_{2^n-1})(N+1+\tau_{2^n})\lge N(N+2)\succ N(N+1)=\hat\eta_1\hat\eta_2.
\]
So in the following we assume $0\le k<2^n-1$. Then by (\ref{eq:relation-eta-tau}) we have
\begin{align*}
	\eta_{i+1}\eta_{i+2}\eta_{i+3}&=(N+\tau_k)(N+1+\tau_{k+1})(N+2-2\tau_{k+1})\\
	&\lge N(N+1)(N+2)\succ N(N+1)(N+1)=\hat\eta_1\hat\eta_2\hat\eta_3,
\end{align*}
proving (\ref{eq:inequ-1-3N}) when $i=3k+1$.

Case (III). $i=3k+2$ with $0\le k<2^n$. Then $\eta_{i+1}=N+1+\tau_{k+1}\ge N+1>\hat\eta_1$, which gives the first inequality of (\ref{eq:inequ-1-3N}). On the other hand, if $k=2^n-1$ then we have 
$\eta_{i+1}=N+1+\tau_{2^n}=N+2=\eta_1$; and if $0\le k<2^n-1$, then by (\ref{eq:relation-eta-tau}) we have 
\begin{align*}
	\eta_{i+1}\eta_{i+2}\eta_{i+3} & =(N+1+\tau_{k+1})(N+2-2\tau_{k+1})(N+\tau_{k+1}) \\
	& \lle (N+2)N(N+1)\prec(N+2)N(N+2)= \eta_1\eta_2\eta_3,
\end{align*}
completing the proof. 
\end{proof}

\begin{proof}[Proof of Proposition \ref{prop:transcendental-betac-3N}]
Note that the Thue-Morse sequence $(\tau_i)_{i=0}^\f$ satisfies $\tau_0=0, \tau_{2k}=\tau_k$ and $\tau_{2k+1}=1-\tau_k$ for all $k\ge 0$. Then by Definition \ref{def:beta-c} and (\ref{eq:relation-eta-tau}) it follows that $(\eta_i)$ is an expansion of $1$ in the base $\beta_c=\beta_c(M)$. So, by Lemmas \ref{lem:char-quasi-greedy-expansion} and \ref{lem:inequality-eta-i} it follows that $\de(\beta_c)=\eta_1\eta_2\ldots$. Hence, 
by (\ref{eq:relation-eta-tau}) we obtain  
\begin{align*} 
	1=\sum_{i=1}^{\f}\frac{\eta_i}{\beta_c^i} & =\sum_{i=0}^\f\frac{N+2-2\tau_{i}}{\beta_c^{3i+1}}+\sum_{k=0}^{\f}\frac{N+\tau_i}{\beta_c^{3i+2}}+\sum_{i=0}^{\f}\frac{N+1+\tau_{i+1}}{\beta_c^{3i+3}} \\
	& =\sum_{i=1}^{\f}\frac{N+1}{\beta_c^i}+\sum_{i=0}^{\f}\frac{1}{\beta_c^{3i+1}}-\sum_{i=0}^{\f}\frac{1}{\beta_c^{3i+2}}-\frac{2}{\beta_c}\sum_{i=0}^{\f}\frac{\tau_i}{\beta_c^{3i}}+\frac{1}{\beta_c^2}\sum_{i=0}^{\f}\frac{\tau_i}{\beta_c^{3i}}+\frac{1}{\beta_c^3}\sum_{i=0}^{\f}\frac{\tau_i}{\beta_c^{3i}}\\
	&=\frac{N+1}{\beta_c-1}+\frac{\beta_c^2-\beta_c}{\beta_c^3-1}+\left(\frac{-2}{\beta_c}+\frac{1}{\beta_c^2}+\frac{1}{\beta_c^3}\right)\sum_{i=0}^{\f}\frac{\tau_i}{\beta_c^{3i}}.
\end{align*}
Since $\beta_c>N+2$, we have $\frac{-2}{\beta_c}+\frac{1}{\beta_c^2}+\frac{1}{\beta_c^3}\ne 0$. Using $\tau_0=0$ and rearranging the above equation yields
\[
\sum_{i=1}^{\f}\frac{\tau_i}{\beta_c^{3i}}=\frac{ 1-\frac{N+1}{\beta_c-1}-\frac{\beta_c^2-\beta_c}{\beta_c^3-1}}{\frac{-2}{\beta_c}+\frac{1}{\beta_c^2}+\frac{1}{\beta_c^3}}.
\]
Suppose on the contrary that $\beta_c\in(N+2, N+3)$ is algebraic. Then $\beta_c^{-3}\in(0,1)$ is  also algebraic, and by Lemma \ref{lem:mahler-transcendental} it follows that the left-hand side of the above equation is transcendental, while the right-hand side is algebraic, leading to a contradiction. So, $\beta_c$ is transcendental. 
\end{proof} 

\subsection{Patterns in $\widetilde{\U}_{\beta_c, M}$}
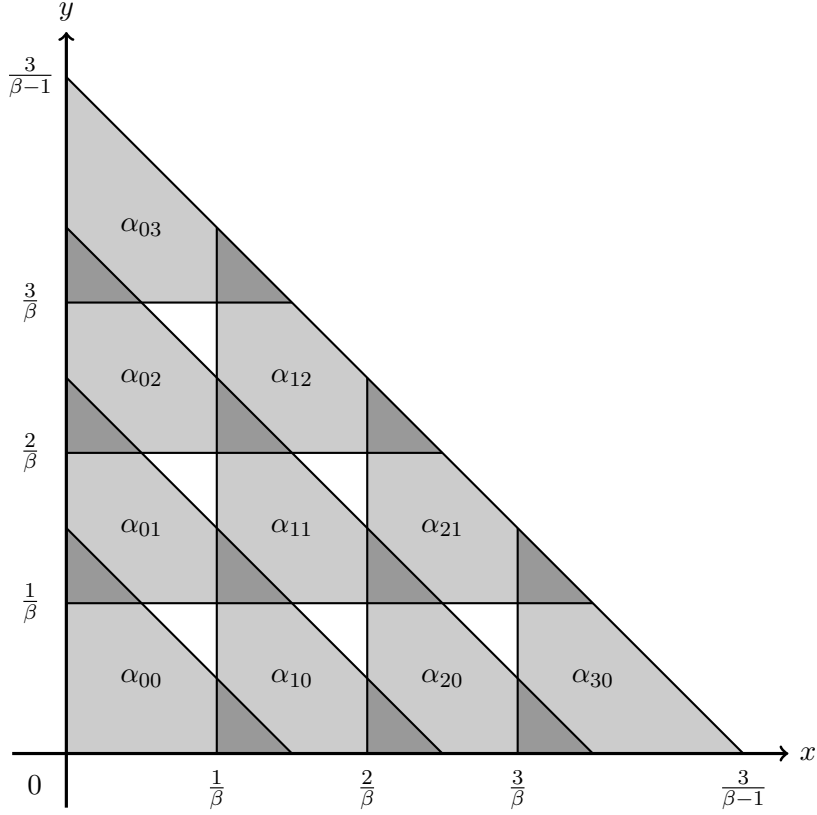
\begin{figure}[h!]
\begin{center}
	\begin{tikzpicture}[
		scale=6,  
		axis/.style={very thick, ->},
		important line/.style={thick},
		dashed line/.style={dashed, thin},
		pile/.style={thick, ->, >=stealth', shorten <=2pt, shorten
			>=2pt},
		every node/.style={color=black}
		]
		
		\fill[black!20]({0},0)--({3/(3-1)-3/3},0)--(0,{3/(3-1)-3/3})--cycle; 
		\fill[black!20]({1/3},0)--({3/(3-1)-2/3},0)--({1/3},{3/(3-1)-3/3})--cycle; 
		\fill[black!20]({2/3},0)--({3/(3-1)-1/3},0)--({2/3},{3/(3-1)-3/3})--cycle; 
		\fill[black!20]({3/3},0)--({3/(3-1)-0},0)--({3/3},{3/(3-1)-3/3})--cycle; 
		
		\fill[black!20]({0},{1/3})--({3/(3-1)-3/3},{1/3})--(0,{3/(3-1)-2/3})--cycle; 
		\fill[black!20]({1/3},{1/3})--({3/(3-1)-2/3},{1/3})--({1/3},{3/(3-1)-2/3})--cycle; 
		\fill[black!20]({2/3},{1/3})--({3/(3-1)-1/3},{1/3})--({2/3},{3/(3-1)-2/3})--cycle; 
		
		\fill[black!20]({0},{2/3})--({3/(3-1)-3/3},{2/3})--(0,{3/(3-1)-1/3})--cycle; 
		\fill[black!20]({1/3},{2/3})--({3/(3-1)-2/3},{2/3})--({1/3},{3/(3-1)-1/3})--cycle; 
		
		\fill[black!20]({0},{3/3})--({3/(3-1)-3/3},{3/3})--(0,{3/(3-1)})--cycle; 
		
		\fill[black!40]({1/3},0)--({3/(3-1)-3/3},0)--({1/3},{3/(3-1)-4/3})--cycle;  
		\fill[black!40]({2/3},0)--({3/(3-1)-2/3},0)--({2/3},{3/(3-1)-4/3})--cycle;  
		\fill[black!40]({3/3},0)--({3/(3-1)-1/3},0)--({3/3},{3/(3-1)-4/3})--cycle;  
		
		\fill[black!40]({1/3},{1/3})--({3/(3-1)-3/3},{1/3})--({1/3},{3/(3-1)-3/3})--cycle;  
		\fill[black!40]({2/3},{1/3})--({3/(3-1)-2/3},{1/3})--({2/3},{3/(3-1)-3/3})--cycle;  
		\fill[black!40]({3/3},{1/3})--({3/(3-1)-1/3},{1/3})--({3/3},{3/(3-1)-3/3})--cycle;  
		
		\fill[black!40](0,{1/3})--(0,{3/(3-1)-3/3})--({3/(3-1)-4/3},{1/3})--cycle;  
		\fill[black!40](0,{2/3})--(0,{3/(3-1)-2/3})--({3/(3-1)-4/3},{2/3})--cycle;  
		\fill[black!40](0,{3/3})--(0,{3/(3-1)-1/3})--({3/(3-1)-4/3},{3/3})--cycle;  
		
		\fill[black!40]({1/3},{2/3})--({1/3},{3/(3-1)-2/3})--({3/(3-1)-3/3},{2/3})--cycle;  
		\fill[black!40]({1/3},{3/3})--({1/3},{3/(3-1)-1/3})--({3/(3-1)-3/3},{3/3})--cycle;  
		
		\fill[black!40]({2/3},{2/3})--({2/3},{3/(3-1)-2/3})--({3/(3-1)-2/3},{2/3})--cycle; 
		
		\draw[axis] (-0.12,0)  -- ({(3/(3-1))+0.1},0) node(xline)[right]
		{$x$};
		\draw[axis] (0,-0.12) -- (0,{3/(3-1)+0.1}) node(yline)[above] {$y$};
		\node[] at (-0.07,-0.07){$0$};
		\draw[important line] (0,{3/(3-1)})--({3/(3-1)},0);
		
		\draw[important line] ({1/3},0)--({1/3},{3/(3-1)-1/3});
		\draw[important line] ({2/3},0)--({2/3},{3/(3-1)-2/3});
		\draw[important line] ({3/3},0)--({3/3},{3/(3-1)-3/3});
		
		\draw[important line] (0,{1/3})--({3/(3-1)-1/3},{1/3});
		\draw[important line] (0,{2/3})--({3/(3-1)-2/3},{2/3});
		\draw[important line] (0,{3/3})--({3/(3-1)-3/3},{3/3});
		
		\draw[important line] ({3/(3-1)-3/3},0)--(0,{3/(3-1)-3/3});
		\draw[important line] ({3/(3-1)-2/3},0)--(0,{3/(3-1)-2/3});
		\draw[important line] ({3/(3-1)-1/3},0)--(0,{3/(3-1)-1/3});
		
		\node[] at ({1/(2*3)}, {1/(2*3)}){$\al_{00}$};
		\node[] at ({3/(2*3)}, {1/(2*3)}){$\al_{10}$};
		\node[] at ({5/(2*3)}, {1/(2*3)}){$\al_{20}$};
		\node[] at ({7/(2*3)}, {1/(2*3)}){$\al_{30}$};
		
		\node[] at ({1/(2*3)}, {3/(2*3)}){$\al_{01}$};
		
		\node[] at ({1/(2*3)}, {5/(2*3)}){$\al_{02}$};
		
		\node[] at ({1/(2*3)}, {7/(2*3)}){$\al_{03}$};
		
		\node[] at ({3/(2*3)}, {3/(2*3)}){$\al_{11}$};
		
		\node[] at ({5/(2*3)}, {3/(2*3)}){$\al_{21}$};
		
		\node[] at ({3/(2*3)}, {5/(2*3)}){$\al_{12}$};
		
		\node[] at ({1/(3)}, {-0.08}){$\frac{1}{\beta}$};
		
		\node[] at ({2/(3)}, {-0.08}){$\frac{2}{\beta}$};
		\node[] at ({3/(3)}, {-0.08}){$\frac{3}{\beta}$};
		
		\node[] at ({3/(3-1)}, {-0.08}){$\frac{3}{\beta-1}$};
		
		\node[] at ( {-0.08}, {1/(3)}){$\frac{1}{\beta}$};
		
		\node[] at ( {-0.08}, {2/(3)}){$\frac{2}{\beta}$};
		
		\node[] at ( {-0.08}, {3/(3)}){$\frac{3}{\beta}$};
		\node[] at ( {-0.08}, {3/(3-1)}){$\frac{3}{\beta-1}$};

		%
		%
		%

	\end{tikzpicture}
\end{center}
\caption{The graph for the first generation of $S_{\beta, M}$ with $M=3$ and $\beta=3$.   The convex hull $\Delta_{\beta, M}$ is the triangle with vertices $(0, 0), (3/(\beta-1), 0)$ and $(0, 3/(\beta-1))$. Each light grey triangle corresponds to a $f_{\beta, \al_{ij}}(\Delta_{\beta,M})$ for some $\al_{ij}\in\Omega_M$.  Furthermore, the overlap region $O_{\beta, M}$ is the union of   eighteen  small dark grey triangles.}\label{Fig:S-3}
\end{figure}
In the following we will study possible block patterns appearing in  sequences of $\widetilde{\U}_{\beta_c, M}$. To describe these patterns  we  construct a sequence of blocks $\set{\c_n}_{n=0}^\f$ in $\Om_M^*$. Similar to the words $(\a_n)$ and $(\b_n)$ for the cases $M=3N+1$ and $M=3N+2$ respectively, these words $(\c_n)$ play an important role to describe $\widetilde{\U}_{\beta,M}$ for $\beta\le \beta_c(M)$ with $M=3N+3$.
Let 
\begin{align*}
\pazocal A_3 & :=\left\{\al_{(N+2)(N+1)}\al_{N(N+2)}\al_{(N+2)N}, \quad \al_{N(N+1)}\al_{(N+1)(N+2)}\al_{(N+1)N},\right. \\
&\hspace{0.75cm}\left.\al_{(N+2)(N+1)}\al_{N(N+2)}\al_{(N+1)N}, \quad\al_{N(N+1)}\al_{(N+1)(N+2)}\al_{(N+2)N}\right\}.
\end{align*} 
We define a substitution $\Phi_3: \pazocal A_3\to\pazocal A_3$   by
\begin{equation}
\label{eq:substitution-Phi-3N}
\begin{split}
	\al_{(N+2)(N+1)}\al_{N(N+2)}\al_{(N+2)N}\quad&\longleftrightarrow \quad  \al_{N(N+1)}\al_{(N+1)(N+2)}\al_{(N+1)N}, \\
	\al_{(N+2)(N+1)}\al_{N(N+2)}\al_{(N+1)N}\quad & \longleftrightarrow\quad \al_{N(N+1)}\al_{(N+1)(N+2)}\al_{(N+2)N}.
\end{split}
\end{equation}
Then the sequence $\set{\c_n}_{n=0}^\f\subset\pazocal A_3^*$ is defined recursively by (see Figure \ref{Fig:S-3})
\begin{equation}\label{eq:cn-3N}
\c_0=\al_{(N+1)(N+1)},\quad \c_1=\al_{(N+2)(N+1)}\al_{N(N+2)}\al_{(N+1)N},\quad\textrm{and}\quad \c_{n+1}=\c_n^+\Phi_3(\c_n^+) 
\end{equation}for all $n\ge 1$,
where $\c_n^+$ is a word in $\pazocal A_3^*$ by changing the last digit of $\c_n$ from $\al_{(N+1)N}$ to $\al_{(N+2)N}$.
Then each $\c_n$ has length $3\cdot 2^n$ and ends with $\al_{(N+1)N}$.     So, 
\[
\c_0\sim\left(\begin{array}{c}
N+1\\
N+1\\
N+1
\end{array}\right),\quad \c_1\sim\left(\begin{array}{ccc}
N+2&N&N+1 \\
N+1&N+2&N \\
N&N+1&N+2
\end{array}\right), 
\]
and
\[
\c_2\sim \left(\begin{array}{cccccc}
N+2&N&N+2&N&N+1&N+1 \\
N+1&N+2&N&N+1&N+2&N \\
N&N+1&N+1&N+2&N&N+2
\end{array}\right).
\]

Next  we recursively define a sequence of blocks $(\sv_n)$ which converges to $(\eta_i)$ componentwisely: 
\begin{equation}\label{eq:vn-3N}
\sv_0=N+1,\quad \sv_1=(N+2)N(N+1),\quad\textrm{and}\quad \sv_{n+1}:=\sv_n^+\re_3(\sv_n^+)\quad\forall n\ge 1,
\end{equation}
where $\re_3$ is the block map defined as in (\ref{eq:substitution-3N}).
So, for example,
\begin{align*}\sv_2&=(N+2)N(N+2)N(N+1)(N+1),\\
\sv_3&= (N+2)N(N+2)N(N+1)(N+2)\,N(N+1) (N+1)(N+2)N(N+1).\end{align*}
Then by (\ref{eq:relation-eta-tau}) it follows that $(\eta_i)$ is the componentwise limit of the sequence $(\sv_n)$, i.e., 
\begin{equation*}\label{eq:eta-i}
(\eta_i)=\lim_{n\to\f}\sv_n\in \Sigma_3^\N.
\end{equation*} 
One can check that $(\eta_i)$ begins with $\sv_n^+$ for any $n\in\N$.

By comparing the definitions of $(\sv_n)$ and $(\c_n)$ in  (\ref{eq:vn-3N}) and (\ref{eq:cn-3N}) respectively,  it follows that 
\begin{equation}\label{eq:relation-cn-vn}
\c_n^1=\sv_n,\quad \c_n^2=((N+1)(N+2)N)^n\quad\textrm{and}\quad \overline{\c_n^\+}=\re_3(\sv_n)
\end{equation}
for any $n\ge 0$. Remember  that for  a block $\d=d_1d_2\ldots d_k\in\Omega_M^*$ we write $\d^1=d_1^1d_2^1\ldots d_k^1, \d^2=d_1^2d_2^2\ldots d_k^2$ and $\overline{\d^\+}=\overline{d_1^\+}\,\overline{d_2^\+}\ldots\overline{d_k^\+}$, where each $d_i=(d_i^1, d_i^2)$ and $\overline{d_i^\+}=M-d_i^\+$ with $d_i^\+=d_i^1+d_i^2$.

We aim at proving  the following pattern result for $\widetilde{\U}_{\beta_{c}, M}$.

\begin{proposition}
\label{prop:patterns-3N}
Let $(d_i)\in\widetilde{\U}_{\beta_c, M}$ and $n\in\N_0$. Then  there exists $m_*\in\N$ such that for all $m>m_*$ the following statements hold.
\begin{enumerate}[{\rm(i)}]
	\item If $d^1_m<M$ and $d_{m+1}\ldots d_{m+3\cdot 2^n}=\c_{n+1}^+$, then 
	$
	d_{m+3\cdot 2^n+1}\ldots d_{m+3\cdot 2^{n+1}}=\Phi_3(\c_{n+1})$ or $\Phi_3(\c_{n+1}^+).
	$
	\item If $d^1_m>0$ and $d_{m+1}\ldots d_{m+3\cdot 2^n}=\Phi_3(\c_{n+1}^+)$, then 
	$
	d_{m+3\cdot 2^n+1}\ldots d_{m+3\cdot 2^{n+1}}=\c_{n+1}$ or  $\c_{n+1}^+.
	$
\end{enumerate} 
\end{proposition}
To prove Proposition \ref{prop:patterns-3N} we need the following lemma.
\begin{lemma}\label{beta-3N-obs}
For any $(d_i)\in\widetilde{\U}_{\beta_c, M}\setminus\widetilde{\U}_{\beta_0, M}$ there exists $m_*\in\N$ such that 
\[d_{m_*+j}^1,\quad d_{m_*+j}^2,\quad \overline{d_{m_*+j}^\+}\in\set{N,  N+1, N+2}\quad\forall  j\ge 1.
\]
\end{lemma}
\begin{proof}
Take $(d_i)\in\widetilde{\U}_{\beta_c, M}\setminus\widetilde{\U}_{\beta_0, M}$. Note by Proposition \ref{prop:beta-G} that $(d_i^1)\in\set{0,1,\ldots, M}^\f$ does not end with $M^\f$. Then there exists $m_1\in\N$ such that $d_{m_1}^1<M$. Since $\de(\beta_c)\lle (N+2)^\f$ by Proposition \ref{prop:transcendental-betac-3N}, we obtain by Proposition \ref{prop:lexicographical-characterization-U}   that $d_{m_1+1}^1\le N+2<3N+3=M$; and again  we can deduce that $d_{m_1+j}^1<M$ for all $j\ge 1$. Similarly, since $(d_i)\notin\widetilde{\U}_{\beta_0, M}$, by Proposition \ref{prop:beta-G} there exists $m_2>m_1$ such that $d_{m_2}^2<M$. The same argument as above yields that $d_{m_2+j}^2<M$ for all $j\ge 1$. Also, we can find $m_*=m_3>m_2$ such that $\overline{d_{m_*+j}^\+}<M$ for all $j\ge 1$. So, for all $j\ge 1$ we have 
\begin{equation}\label{eq:nov8-1}
	d_{m_*+j}^1<M, \quad d_{m_*+j}^2<M\quad\textrm{and}\quad \overline{d_{m_*+j}^\+}<M.
\end{equation}
Since $\de(\beta_c)\lle(N+2)^\f$, by (\ref{eq:sum-equal}) and (\ref{eq:nov8-1}) it follows that 
\[d_{m_*+j}^1, d_{m_*+j}^2, \overline{d_{m_*+j}^\+}\in\set{N, N+1, N+2}\] for all $j\ge 1$.
\end{proof}

First we prove Proposition \ref{prop:patterns-3N} for $n=0$ and $n=1$.

\begin{lemma}
\label{lem:pattern-3N-1}
For any $(d_i)\in\widetilde{\U}_{\beta_c, M}$ there exists   $m_*\in\N$ such that for all $m>m_*$ the following  statements hold.
\begin{enumerate}[{\rm(i)}]
	\item If $d^1_m<M$ and $d_{m+1}d_{m+2}d_{m+3}=\c_1^+$, then 
	$
	d_{m+4}d_{m+5}d_{m+6}=\Phi_3(\c_1)$ or $\Phi_3(\c_1^+).
	$
	\item If $d^1_m>0$ and $d_{m+1}d_{m+2}d_{m+3}=\Phi_3(\c_1^+)$, then 
	$
	d_{m+4}d_{m+5} d_{m+6}=\c_1$ or $\c_1^+.
	$
\end{enumerate}
\end{lemma}

\begin{proof}
Since the proof of (ii) is similar, we only prove (i). Take $(d_i)\in\widetilde{\U}_{\beta_{c},M}$. If $(d_i)\in\widetilde{\U}_{\beta_{0},M}$, then by Proposition \ref{prop:beta-G} it follows that $d_{i+1}d_{i+2}d_{i+3}\notin\set{\c_1^+, \Phi_3(\c_1^+)}$ for all $i\ge 0$, and we are done. In the following we assume $(d_i)\in\widetilde{\U}_{\beta_{c},M}\setminus\widetilde{\U}_{\beta_{0},M}$. Then by Lemma \ref{beta-3N-obs} there exists $m_*\in\N$ such that $d_{m_*+j}^1, d_{m_*+j}^2, \overline{d_{m_*+j}^\+}\in\set{N, N+1, N+2}$ for all $j\ge 1$. Take $m>m_*$. Suppose $d_m^1<M$ and $d_{m+1}d_{m+2}d_{m+3}=\c_1^+=\al_{(N+2)(N+1)}\al_{N(N+2)}\al_{(N+2)N}$. 
Then by Propositions \ref{prop:lexicographical-characterization-U} and \ref{prop:transcendental-betac-3N} it follows that $d_{m+4}\in\set{\al_{N(N+1)}, \al_{N(N+2)}}$. 

If $d_{m+4}=\al_{N(N+2)}$, then by Propositions \ref{prop:lexicographical-characterization-U} and \ref{prop:transcendental-betac-3N} we obtain that 
\[
d_{m+1}\ldots d_{m+9}\sim\begin{pmatrix}
	(N+2)N(N+2)&N(N+1)(N+2)&N(N+1)(N+1)\\
	(N+1)(N+2)N&(N+2)N(N+1)&(N+2)N(N+1)\\
	N(N+1)(N+1)&(N+1)(N+2)N&(N+1)(N+2)(N+1)
\end{pmatrix}
\]or
\[ 
d_{m+1}\ldots d_{m+9}\sim\begin{pmatrix}
	(N+2)N(N+2)&N(N+1)(N+2)&N(N+1)(N+1)\\
	(N+1)(N+2)N&(N+2)N(N+1)&(N+1)(N+2)(N+1)\\
	N(N+1)(N+1)&(N+1)(N+2)N&(N+2)N(N+1)
\end{pmatrix}.
\]
However, the forbidden block $(N+1)(N+2)(N+1)$ occurs in $\overline{d_{m+1}^\+}\ldots \overline{d_{m+9}^\+}$ in the first case, and occurs in $d_{m+1}^2\ldots d_{m+9}^2$ in the second case, leading to a contradiction to $(d_i)\in\widetilde{\U}_{\beta_c, M}$. So, $d_{m+4}=\al_{N(N+1)}$.

By (\ref{eq:sum-equal}), Propositions \ref{prop:lexicographical-characterization-U} and \ref{prop:transcendental-betac-3N} it follows that $d_{m+5}=\al_{(N+1)(N+2)}$ and $d_{m+6}=\al_{(N+2)N}$ or $\al_{(N+1)N}$. This implies that $d_{m+4}d_{m+5}d_{m+6}=\Phi_3(\c_1)$ or $\Phi_3(\c_1^+)$ as desired. 
\end{proof}


\begin{lemma}
\label{lem:pattern-3N-2}
For any $(d_i)\in\widetilde{\U}_{\beta_c, M}$ there exists   $m_*\in\N$ such that for all $m>m_*$ the following  statements hold.
\begin{enumerate}[{\rm(i)}]
	\item If $d^1_m<M$ and $d_{m+1}\ldots d_{m+6}=\c_2^+$, then 
	$
	d_{m+7}\ldots d_{m+12}=\Phi_3(\c_2)$ or $\Phi_3(\c_2^+).
	$
	\item If $d^1_m>0$ and $d_{m+1}\ldots d_{m+6}=\Phi_3(\c_2^+)$, then 
	$
	d_{m+7}\ldots  d_{m+12}=\c_2$ or $\c_2^+.
	$
\end{enumerate}
\end{lemma}

\begin{proof}
Since the proof of (ii) is similar, here we only prove (i). Take $(d_i)\in\widetilde{\U}_{\beta_c, M}$. By the same argument as in the proof of Lemma \ref{lem:pattern-3N-1} we may assume $(d_i)\in\widetilde{\U}_{\beta_{c},M}\setminus\widetilde{\U}_{\beta_{0},M}$. Then by Lemma \ref{beta-3N-obs} there exists $m_*\in\N$ such that each component of $d_{m_*+j}$ belongs to $\set{N, N+1, N+2}$ for all $j\ge 1$. Take $m>m_*$. Suppose $d_m^1<M$ and $d_{m+1}\ldots d_{m+6}=\c_2^+=\al_{(N+2)(N+1)}\al_{N(N+2)}\al_{(N+2)(N+1)}\al_{N(N+1)}\al_{(N+1)(N+2)}\al_{(N+2)N}$ for some $m>m_*$. Then by  (\ref{eq:sum-equal}), and Propositions \ref{prop:lexicographical-characterization-U} and \ref{prop:transcendental-betac-3N} it follows that $d_{m+7}\in\set{\al_{N(N+1)}, \al_{N(N+2)}}$.

If $d_{m+7}=\al_{N(N+2)}$, then by  (\ref{eq:sum-equal}), and Propositions \ref{prop:lexicographical-characterization-U} and \ref{prop:transcendental-betac-3N} we obtain 
\[
d_{m+1}\ldots d_{m+9}\sim\begin{pmatrix}
	(N+2)N(N+2)&N(N+1)(N+2)&N(N+1)(N+1)\\
	(N+1)(N+2)N&(N+1)(N+2)N&(N+2)N(N+1)\\
	N(N+1)(N+1)&(N+2)N(N+1)&(N+1)(N+2)(N+1)
\end{pmatrix}.
\]
Thus, the forbidden block $(N+1)(N+2)(N+1)$ appears in $\overline{d_{m+1}^\+}\ldots \overline{d_{m+9}^\+}$, which leads to a contradiction to $(d_i)\in\widetilde{\U}_{\beta_c, M}$. So, $d_{m+7}=\al_{N(N+1)}$. 

Hence, by (\ref{eq:sum-equal}), Propositions \ref{prop:lexicographical-characterization-U} and \ref{prop:transcendental-betac-3N} it follows that \[d_{m+7}d_{m+8}d_{m+9}=\al_{N(N+1)}\al_{(N+1)(N+2)}\al_{(N+1)N}=\Phi_3(\c_1^+).\]
Since $d_{m+6}^1=N+2>0$, by Lemma \ref{lem:pattern-3N-1} (ii) we obtain $d_{m+10}d_{m+11}d_{m+12}=\c_1$ or $\c_1^+$. So,
\[
d_{m+7}\ldots d_{m+12}=\Phi_3(\c_1^+)\c_1=\Phi_3(\c_2^+)\quad\textrm{or}\quad \Phi_3(\c_1^+)\c_1^+=\Phi_3(\c_2),
\]
completing the proof.
\end{proof}

\begin{proof}
[Proof of Proposition \ref{prop:patterns-3N}]
Since the proof of (ii) is similar, we only prove (i). Take $(d_i)\in\widetilde{\U}_{\beta_c, M}$. By Proposition \ref{prop:beta-G} it suffices to consider $(d_i)\in\widetilde{\U}_{\beta_c, M}\setminus\widetilde{\U}_{\beta_0, M}$. In view of Lemma \ref{beta-3N-obs}, there exists $m_*\in\N$ such that each component of $d_{m_*+j}$ belongs to $\set{N, N+1, N+2}$ for all $j\ge 1$. Take $m>m_*$. Suppose $d_m^1<M$ and $d_{m+1}\ldots d_{m+3\cdot 2^n}=\c_{n+1}^+$ for some $m>m_*$. We will prove 
\begin{equation}\label{eq:nov9-1}
	d_{m+3\cdot 2^n+1}\ldots d_{m+3\cdot 2^{n+1}}=\Phi_3(\c_{n+1}^+)\quad\textrm{or}\quad \Phi_3(\c_{n+1}),
\end{equation}
which will be done  by induction on $n$. If $n=0$ or $1$, then (\ref{eq:nov9-1}) holds by Lemmas \ref{lem:pattern-3N-1} and \ref{lem:pattern-3N-2}. Now suppose (\ref{eq:nov9-1}) holds for all $n\le k$ with $k\in\N$, and we consider it for $n=k+1$.

Note that $d_{m+1}\ldots d_{m+3\cdot 2^{k+1}}=\c_{k+2}^+$. By (\ref{eq:vn-3N}) and (\ref{eq:relation-cn-vn}) it follows that \[d_{m+1}^1\ldots d_{m+3\cdot 2^{k+1}}^1=\sv_{k+2}^+=\sv_{k+1}^+\phi_3(\sv_{k+1}). \]
Since by (\ref{eq:eta-i}) the sequence $\de(\beta_c)$ begins with $\sv_{k+2}^+\phi_3(\sv_{k+2})$, by (\ref{eq:vn-3N}) and Proposition \ref{prop:lexicographical-characterization-U} it follows that 
\begin{equation}\label{eq:nov9-2}
	d_{m+3\cdot 2^{k+1}+1}^1\ldots d_{m+3\cdot 2^{k+2}}^1\lle \phi_3(\sv_{k+2})=\phi_3(\sv_{k+1}^+)\sv_{k+1}^+=\phi_3(\sv_k^+)\sv_k\sv_{k+1}^+.
\end{equation}
Note that $d_{m+1}\ldots d_{m+3\cdot 2^{k+1}}=\c_{k+2}^+=\c_{k+1}^+\Phi_3(\c_k^+)\c_k^+$. Then $d_{m+3\cdot(2^k+2^{k-1})}^1=N+1<M$ and $d_{m+3(2^k+2^{k-1})+1}\ldots d_{m+3\cdot 2^{k+1}}=\c_k^+$. Then by the induction hypothesis for $n=k-1$ we have 
\[d_{m+3\cdot 2^{k+1}+1}\ldots d_{m+3\cdot(2^{k+1}+2^{k-1})}=\Phi_3(\c_k^+)\quad \textrm{or}\quad \Phi_3(\c_k).\]
If it takes $\Phi_3(\c_k)$, then $d^1_{m+3\cdot 2^{k+1}+1}\ldots d^1_{m+3\cdot(2^{k+1}+2^{k-1})}=\phi_3(\sv_k)\succ\phi_3(\sv_k^+)$, leading to a contradiction to (\ref{eq:nov9-2}). So, $d_{m+3\cdot 2^{k+1}+1}\ldots d_{m+3\cdot(2^{k+1}+2^{k-1})}=\Phi_3(\c_k^+)$. Since $d_{m+3\cdot 2^{k+1}}^1=N+2>0$, by the induction hypothesis we obtain that $d_{m+3\cdot 2^{k+1}+1}\ldots d_{m+3\cdot(2^{k+1}+2^{k})}=\Phi_3(\c_{k+1}^+)$ or $\Phi_3(\c_{k+1})$. In terms of (\ref{eq:nov9-2}) the block $\Phi_3(\c_{k+1})$ is not allowed. So, \[d_{m+3\cdot 2^{k+1}+1}\ldots d_{m+3\cdot(2^{k+1}+2^{k})}=\Phi_3(\c_{k+1}^+).  \] Again by using $d_{m+3\cdot 2^{k+1}}^1=N+2>0$ and the induction hypothesis it follows that 
\[d_{m+3\cdot 2^{k+1}+1}\ldots d_{m+3\cdot 2^{k+2}}=\Phi_3(\c_{k+2}^+)\quad \textrm{or}\quad \Phi_3(\c_{k+2}),
\] proving (\ref{eq:nov9-1}) for $n=k+1$. This completes the proof by induction. 
\end{proof}

\subsection{Phase transition of $\widetilde{\U}_{\beta, M}$ for $M=3N+3$}

Note by (\ref{eq:vn-3N}) that $\sv_n=\eta_1\ldots \eta_{3\cdot 2^{n-1}}^-$  for all $n\ge 1$. Then by Lemma \ref{lem:inequality-eta-i} it follows that $\si^j(\sv_n^\f)\lle\sv_n^\f$ for all $j\ge 0$. So, by Lemma \ref{lem:char-quasi-greedy-expansion} there exists $\beta_n=\beta_n(M)\in(N+2, N+3)$ such that 
\[\de(\beta_n)=\sv_n^\f.\]
Since $\sv_n\nearrow (\eta_i)=\de(\beta_c)$, again by Lemma \ref{lem:char-quasi-greedy-expansion} it follows that 
\[
\beta_G=\beta_0<\beta_1<\cdots<\beta_n<\beta_{n+1}<\cdots,\quad\textrm{and}\quad \beta_n\nearrow\beta_c\quad\textrm{as }n\to\f.
\]
\begin{proposition}
\label{prop:difference-Un+1-Un-3N}
Let $M=3N+3$ and $n\in\N_0$. Then any sequence in $\widetilde{\U}_{\beta_{n+1},M}\setminus\widetilde{\U}_{\beta_n,M}$ must end in
\[
\bigcup_{j=0}^2\set{  \theta^j(\c_n^\f), \theta^j(\Phi_3(\c_n)^\f)},
\]
where $\theta$ is the cyclic permutation defined in (\ref{eq:theta}).
\end{proposition}

We will prove Proposition \ref{prop:difference-Un+1-Un-3N} by induction on $n$. First we prove it for $n=0, 1$. 

\begin{lemma}
\label{lem:difference-Un+1-Un-3N-0n1}The following statements hold.
\begin{enumerate}[{\rm(i)}]
	\item Any sequence in $\widetilde{\U}_{\beta_{1},M}\setminus\widetilde{\U}_{\beta_0, M}$ must end in $\c_0^\infty$;
	\item  Any sequence in $\widetilde{\U}_{\beta_{2},M}\setminus\widetilde{\U}_{\beta_1,M}$ must end in
	$\bigcup_{j=0}^2\set{\theta^j(\c_1^\infty), \theta^j(\Phi_3(\c_1)^\infty)}$.
\end{enumerate}
\end{lemma}
\begin{proof}
For (i) we take $(d_i)\in\widetilde{\U}_{\beta_1,M}\setminus\widetilde{\U}_{\beta_0, M}$. Note that $\de(\beta_1)=\sv_1^\f=((N+2)N(N+1))^\f$. By the same argument as in the proof of 
Lemma \ref{beta-3N-obs} there exists $m\in\N$ such that each component of $d_{m+j}$ belongs to $\set{N, N+1, N+2}$ for all $j\ge 1$. Then by (\ref{eq:sum-equal}) it follows that   
\[
d_{m+1}=\al_{(N+1)(N+1)}\quad\textrm{or}\quad d_{m+1}\in\bigcup_{j=0}^2\set{\theta^j(\al_{(N+2)(N+1)}), \theta^j(\al_{(N+2)N})}.
\]

Suppose $d_{m+1}=\al_{(N+2)(N+1)}$. By Proposition \ref{prop:lexicographical-characterization-U} and using $\de(\beta_1)=((N+2)N(N+1))^\f$ it follows that 
\[
d_{m+1}d_{m+2}\ldots =(\al_{(N+2)(N+1)}\al_{N(N+2)}\al_{(N+1)N})^\f\quad\textrm{or}\quad (\al_{(N+2)(N+1)}\al_{N(N+1)}\al_{(N+1)(N+2)})^\f.
\]
However, in either case the component sequence $d_{m+1}^1d_{m+2}^1\ldots$ ends with $((N+2)N(N+1))^\f=\sv_1^\f$, which leads to a contradiction to $(d_i)\in\widetilde{\U}_{\beta_1,M}$. 
Similarly, if $d_{m+1}$ equals any other digit in $\bigcup_{j=0}^2\set{\theta^j(\al_{(N+2)(N+1)}), \theta^j(\al_{(N+2)N})}$, we can deduce the same contradiction. So, $d_{m+1}=\al_{(N+1)(N+1)}$; and therefore by the same argument as above we can prove that $d_{m+1}d_{m+2}\ldots=\al_{(N+1)(N+1)}^\f=\c_0^\f$.

Next we prove (ii). Take   $(d_i)\in\widetilde{\U}_{\beta_{2},M}\setminus\widetilde{\U}_{\beta_1,M}$. Note that 
\begin{equation}\label{eq:nov9-3}
	\de(\beta_2)=\sv_2^\f=((N+2)N(N+2)N(N+1)(N+1))^\f.
\end{equation}
By the same argument as in the proof of Lemma \ref{beta-3N-obs} there exists $m_*\in\N$ such that each component of $d_{m_*+j}$ belongs to $\set{N, N+1, N+2}$ for all $j\ge 1$. Again by (i) there exists $m>m_*$ such that 
\[
d_{m+j}\in\bigcup_{k=0}^2\set{\theta^k(\al_{(N+2)(N+1)}), \theta^k(\al_{(N+2)N})}.
\]
Without loss of generality we may assume $d_{m+1}=\al_{(N+2)(N+1)}$. Then by (\ref{eq:nov9-3}) and Proposition \ref{prop:lexicographical-characterization-U} it follows that 
\begin{equation}\label{eq:nov9-4}
	d_{m+1}d_{m+2}d_{m+3}\in\set{\al_{(N+2)(N+1)}\al_{N(N+2)}\al_{(N+1)N}, \al_{(N+2)(N+1)}\al_{N(N+1)}\al_{(N+1)(N+2)}}
\end{equation}
or 
\begin{equation}\label{eq:nov9-5}
	d_{m+1}d_{m+2}d_{m+3}\in\set{\al_{(N+2)(N+1)}\al_{N(N+2)}\al_{(N+2)N}, \al_{(N+2)(N+1)}\al_{N(N+1)}\al_{(N+2)(N+1)}}.
\end{equation}
We will prove by contradiction  that (\ref{eq:nov9-5}) never occurs.

If $d_{m+1}d_{m+2}d_{m+3}=\al_{(N+2)(N+1)}\al_{N(N+2)}\al_{(N+2)N}$, then by (\ref{eq:nov9-3}) and Proposition \ref{prop:lexicographical-characterization-U} it follows that \[
d_{m+1}d_{m+2}\ldots=(\al_{(N+2)(N+1)}\al_{N(N+2)}\al_{(N+2)N}\al_{N(N+1)}\al_{(N+1)(N+2)}\al_{(N+1)N})^\f,
\]
which implies that $d_{m+1}^1d_{m+2}^1\ldots=\sv_2^\f=\de(\beta_2)$, leading to a contradiction to $(d_i)\in\widetilde{\U}_{\beta_2,M}$. If $d_{m+1}d_{m+2}d_{m+3}=\al_{(N+2)(N+1)}\al_{N(N+1)}\al_{(N+2)(N+1)}$, then by (\ref{eq:nov9-3}) and Proposition \ref{prop:lexicographical-characterization-U} we obtain that 
\[
d_{m+4}d_{m+5}\ldots=(\al_{N(N+2)}\al_{(N+1)N}\al_{(N+1)(N+2)}\al_{(N+2)N}\al_{N(N+1)}\al_{(N+2)(N+1)})^\f,
\]
which implies that $d_{m+4}^1d_{m+5}^1\ldots$ ends with $\sv_2^\f=\de(\beta_2)$, again leading to a contradiction to $(d_i)\in\widetilde{\U}_{\beta_2,M}$.  So, we must have (\ref{eq:nov9-4}).

By the same argument as above we can even prove that any block  $d_{m+j}d_{m+j+1}d_{m+j+2}$ can not belong to
\begin{equation*}\label{eq:nov9-6}
	\bigcup_{k=0}^2\set{\theta^k(\al_{(N+2)(N+1)}\al_{N(N+2)}\al_{(N+2)N}), \theta^k(\al_{(N+2)(N+1)}\al_{N(N+1)}\al_{(N+2)(N+1)})}.
\end{equation*}
From this we can deduce that $d_{m+1}d_{m+2}\ldots$ must end with $\c_1^\f$ as required.

If we choose $d_{m+1}$ for some other digit from $\bigcup_{k=0}^2\set{\theta^k(\al_{(N+2)(N+1)}), \theta^k(\al_{(N+2)N})}$, then by the same argument as above it follows that $d_{m+1}d_{m+2}\ldots$ ends with $\theta^k(\c_1^\f)$ for some $k=0,1,2$. This completes the proof. 
\end{proof}

\begin{proof}
[Proof of Proposition \ref{prop:difference-Un+1-Un-3N}]
We will prove the result  by induction on $n$. By Lemma  \ref{lem:difference-Un+1-Un-3N-0n1}  the result holds for $n=0, 1$.  Now suppose it holds for all $n\le k$ for some $k\ge 1$, and we consider it for $n=k+1$. Take $(d_i)\in\widetilde{\U}_{\beta_{k+2},M}\setminus\widetilde{\U}_{\beta_{k+1},M}$. By Lemma \ref{beta-3N-obs} there exists $m_*\in\N$ such that   each component of $d_{m_*+j}$ belongs to $\set{N, N+1, N+2}$ for all $j\ge 1$. Note that $\widetilde{\U}_{\beta_{k+1},M}=\widetilde{\U}_{\beta_{0},M}\cup\bigcup_{j=0}^k(\widetilde{\U}_{\beta_{j+1},M}\setminus\widetilde{\U}_{\beta_j, M})$. 
By Proposition \ref{prop:beta-G} and the induction hypothesis it follows that any sequence in $\widetilde{\U}_{\beta_{k+1},M}$ must end in 
\begin{equation}\label{eq:nov10-1}
	\bigcup_{j=0}^k\bigcup_{l=0}^2\set{\theta^l(\c_j^\f), \theta^l(\Phi_3(\c_j^\f))}.
\end{equation}	 	Observe by (\ref{eq:relation-cn-vn}) that $\c_j^1=\sv_j$ for all $j\ge 0$, and $\sv_0^\f\prec \sv_1^\f\prec \cdots\prec \sv_k^\f\prec\sv_{k+1}^\f$. Since $\de(\beta_j)=\sv_j^\f$ for all $j\ge 0$,  by (\ref{eq:cn-3N}) and (\ref{eq:nov10-1}) there exists $m>m_*$ such that $d_{m+1}\ldots d_{m+3\cdot 2^{k-1}}\in\bigcup_{j=0}^2\set{\theta^j(\c_{k}^+), \theta^j(\Phi_3(\c_{k}^+))}$. Without loss of generality we may assume 
$d_{m+1}\ldots d_{m+3\cdot 2^{k}}=\c_{k}^+$. Note that $d^1_{m}\le N+2<M$. Then by Proposition \ref{prop:patterns-3N} it follows that 
\begin{equation}\label{eq:nov10-2}
	d_{m+1}\ldots d_{m+3\cdot 2^{k}}=\c_{k+1}\quad\textrm{or}\quad \c_{k+1}^+.
\end{equation}
If $d_{m+1}\ldots d_{m+3\cdot 2^{k}}=\c_{k+1}^+$, then by Proposition \ref{prop:patterns-3N} and (\ref{eq:relation-cn-vn}) it follows that 
\[
d_{m+1}^1d_{m+2}^1\ldots \lge \sv_{k+2}^\f=\de(\beta_{k+2}), 
\] 
which   leads to a contradiction to  $(d_i)\in\widetilde{\U}_{\beta_{k+2}, M}$.

Hence, by (\ref{eq:nov10-2}) we must have $d_{m+1}\ldots d_{m+3\cdot 2^{k}}=\c_{k+1}=\c_k^+\Phi_3(\c_k^+)$. Again, by Proposition \ref{prop:patterns-3N} we obtain that 
$
d_{m+3\cdot 2^{k}+1}\ldots d_{m+3(2^{k}+2^{k-1})}=\c_k$ or $ \c_k^+.
$
If it takes $\c_k$, then $d_{m+3\cdot 2^{k-1}+1}\ldots d_{m+3\cdot(2^{k}+2^{k-1})}=\Phi_3(\c_k^+)\c_k=\Phi_3(\c_{k+1}^+)$, and the same argument as above will lead to
a contradiction. 	 	
So,   $d_{m+3\cdot 2^{k}+1}\ldots d_{m+3(2^{k}+2^{k-1})}=\c_k^+$. Continuing this argument  we can prove that 
$
d_{m+1}d_{m+2}\ldots =(\c_{k}^+\Phi_3(\c_k^+))^\f=\c_{k+1}^\f, 
$
which establishes the result for $n=k+1$. This completes the proof by induction. 
\end{proof}

\begin{proof}
[Proof of Theorem \ref{th:main-result} for $M=3N+3$]
The proof is similar to that for $M=3N+1$. By using Proposition \ref{prop:difference-Un+1-Un-3N} we can deduce that for $\beta<\beta_c(M)$ the set $\widetilde{\U}_{\beta, M}$ is at most countable; and by Proposition \ref{prop:patterns-3N} we can prove that for $\beta=\beta_c(M)$ the set $\widetilde{\U}_{\beta, M}$ is uncountable but has zero Hausdorff dimension. Furthermore, for $\beta>\beta_c(M)$ we can construct a subshift of finite type $Z_k\subset\widetilde{\U}_{\beta, M}$ having positive Hausdorff dimension, where $Z_k$ is analogous to $X_k$ as described in Figure \ref{fig:1} with $\a_k, \a_k^+, \Phi_1(\a_k)$ and $\Phi_1(\a_k^+)$ replaced by $\c_k, \c_k^+, \Phi_3(\c_k)$ and $\Phi_3(\c_k^+)$, respectively. 
\end{proof}

\section*{Acknowledgements}
{The authors thank the anonymous referee for many useful comments which greatly improve the presentation of the paper. They also thank Wolfgang Steiner for his insight on an alternate definition of the Thue-Morse type sequences. The first author was supported by the National Natural Science Foundation of China No.~12301108. The second author was supported by the Chongqing Natural Science Foundation: CQYC20220511052 and the Scientific Research Innovation Capacity Support Project for Young Faculty No.~ZYGXQNISKYCXNLZCXM-P2P. The third author was supported by the {National Natural Science Foundation of China No.~12071148 and 12471085.}}.

%


\begin{thebibliography}{10}

\bibitem{Allaart-Kong-2023}
P.~Allaart and D.~Kong.
\newblock Critical values for the {$\beta$}-transformation with a hole at 0.
\newblock \em Ergodic Theory Dyn. Syst., 43(6):1785--1828, 2023.

\bibitem{Allouche_Shallit_1999}
J.-P. Allouche and J.~Shallit.
\newblock The ubiquitous {P}rouhet-{T}hue-{M}orse sequence.
\newblock In {\em Sequences and their applications ({S}ingapore, 1998)},
  Springer Ser. Discrete Math. Theor. Comput. Sci., pages 1--16. Springer,
  London, 1999.

\bibitem{Allouche_Shallit_2003}
J.-P. Allouche and J.~Shallit.
\newblock {\em Automatic sequences: theory, applications, generalizations}.
\newblock Cambridge University Press, Cambridge, 2003.


\bibitem{Baiocchi_Komornik_2007}
C.~Baiocchi and V.~Komornik.
\newblock Greedy and quasi-greedy expansions in non-integer bases.
\newblock {\em arXiv:0710.3001v1}, 2007.

\bibitem{Bro-Mon-Sid-04}
D.~Broomhead, J.~Montaldi, and N.~Sidorov.
\newblock Golden gaskets: variations on the {S}ierpi\'nski sieve.
\newblock {\em Nonlinearity}, 17(4):1455--1480, 2004.

\bibitem{Bunimovich-Yurchenko-2011}
L.~A. Bunimovich and A.~Yurchenko.
\newblock Where to place a hole to achieve a maximal escape rate.
\newblock {\em Israel J. Math.}, 182:229--252, 2011.

\bibitem{Lyndsey-2016}
L.~Clark.
\newblock The {$\beta$}-transformation with a hole.
\newblock {\em Discrete Contin. Dyn. Syst.}, 36(3):1249--1269, 2016.

\bibitem{Daroczy_Katai_1993}
Z.~Dar{\'o}czy and I.~K{\'a}tai.
\newblock Univoque sequences.
\newblock {\em Publ. Math. Debrecen}, 42(3-4):397--407, 1993.

\bibitem{DeVries_Komornik_2008}
M.~de~Vries and V.~Komornik.
\newblock Unique expansions of real numbers.
\newblock {\em Adv. Math.}, 221(2):390--427, 2009.

\bibitem{Vries-Komornik-Loreti-2022}
M.~de~Vries, V.~Komornik, and P.~Loreti.
\newblock Topology of univoque sets in real base expansions.
\newblock {\em Topology Appl.}, 312:Paper No. 108085, 36, 2022.

\bibitem{Demers-Wright-Young-2010}
M.~Demers, P.~Wright, and L.-S. Young.
\newblock Escape rates and physically relevant measures for billiards with
  small holes.
\newblock {\em Comm. Math. Phys.}, 294(2):353--388, 2010.

\bibitem{Demmers-2005}
M.~F. Demers.
\newblock Markov extensions for dynamical systems with holes: an application to
  expanding maps of the interval.
\newblock {\em Israel J. Math.}, 146:189--221, 2005.

\bibitem{Demers-Young-2006}
M.~F. Demers and L.-S. Young.
\newblock Escape rates and conditionally invariant measures.
\newblock {\em Nonlinearity}, 19(2):377--397, 2006.

\bibitem{Erdos_Joo_Komornik_1990}
P.~Erd\H{o}s, I.~Jo\'{o}, and V.~Komornik.
\newblock Characterization of the unique expansions $1=\sum_{i=1}^\infty
  q^{-n_i}$ and related problems.
\newblock {\em Bull. Soc. Math. France}, 118:377--390, 1990.

\bibitem{Falconer_1990}
K.~Falconer.
\newblock {\em Fractal geometry: mathematical foundations and applications}.
\newblock John Wiley \& Sons Ltd., Chichester, 1990.


\bibitem{Glendinning-Sidorov-2015}
P.~Glendinning and N.~Sidorov.
\newblock The doubling map with asymmetrical holes.
\newblock \em Ergodic Theory Dyn. Syst., 35(4):1208--1228, 2015.

\bibitem{Hochman-2015}
M.~Hochman.
\newblock On self-similar sets with overlaps and inverse theorems for entropy
  in $r^d$.
\newblock {\em Mem. Amer. Math. Soc. to appear (arXiv:1503.09043)}, 2015.

\bibitem{Hutchinson_1981}
J.~E. Hutchinson.
\newblock Fractals and self-similarity.
\newblock {\em Indiana Univ. Math. J.}, 30(5):713--747, 1981.

\bibitem{Jordan-Pollicott-06}
T.~Jordan and M.~Pollicott.
\newblock Properties of measures supported on fat {S}ierpinski carpets.
\newblock \em Ergodic Theory Dyn. Syst., 26(3):739--754, 2006.

\bibitem{Kalle-Kong-Langeveld-Li-18}
C.~Kalle, D.~Kong, N.~Langeveld, and W.~Li.
\newblock The $\beta$-transformation with a hole at 0.
\newblock \em Ergodic Theory Dyn. Syst., 40(9):2482--2514, 2020.

\bibitem{Komornik-Kong-Li-17}
V.~Komornik, D.~Kong, and W.~Li.
\newblock Hausdorff dimension of univoque sets and devil's staircase.
\newblock {\em Adv. Math.}, 305:165--196, 2017.

\bibitem{Komornik-Loreti-1998}
V.~Komornik and P.~Loreti.
\newblock Unique developments in non-integer bases.
\newblock {\em Amer. Math. Monthly}, 105(7):636--639, 1998.

\bibitem{Komornik-Loreti-2002}
V.~Komornik and P.~Loreti.
\newblock {Subexpansions, superexpansions and uniqueness properties in noninteger bases},
\newblock {\em Period. Math. Hungar.}, 44(2):195--216, 2002.

\bibitem{Kong-Li-2020}
D.~Kong and W.~Li.
\newblock Critical base for the unique codings of fat {S}ierpinski gasket.
\newblock {\em Nonlinearity}, 33(9):4484--4511, 2020.

\bibitem{Lind_Marcus_1995}
D.~Lind and B.~Marcus.
\newblock {\em An introduction to symbolic dynamics and coding}.
\newblock Cambridge University Press, Cambridge, 1995.

\bibitem{Mahler_1976}
K.~Mahler.
\newblock {\em Lectures on transcendental numbers}.
\newblock Lecture Notes in Mathematics, Vol. 546. Springer-Verlag, Berlin,
  1976.

\bibitem{Parry_1960}
W.~Parry.
\newblock On the $\beta$-expansions of real numbers.
\newblock {\em Acta Math. Acad. Sci. Hungar.}, 11:401--416, 1960.

\bibitem{Pianigiani-Yorke-1979}
G.~Pianigiani and J.~A. Yorke.
\newblock Expanding maps on sets which are almost invariant. {D}ecay and chaos.
\newblock {\em Trans. Amer. Math. Soc.}, 252:351--366, 1979.

\bibitem{Sidorov_2007}
N.~Sidorov.
\newblock Combinatorics of linear iterated function systems with overlaps.
\newblock {\em Nonlinearity}, 20(5):1299--1312, 2007.

\bibitem{Simon-Solomyak-02}
K.~Simon and B.~Solomyak.
\newblock On the dimension of self-similar sets.
\newblock {\em Fractals}, 10(1):59--65, 2002.

\bibitem{Urbanski_1986}
M.~Urba{\'n}ski.
\newblock On {H}ausdorff dimension of invariant sets for expanding maps of a
  circle.
\newblock \em Ergodic Theory Dyn. Syst., 6(2):295--309, 1986.

\end{thebibliography}

\end{document}